\documentclass[journal,twoside,web]{ieeecolor}
\usepackage{generic}
\usepackage{cite}
\usepackage{amsmath,amssymb,amsfonts}
\usepackage{algorithmic}
\usepackage{graphicx}
\usepackage{algorithm,algorithmic}
\usepackage{hyperref}
\usepackage{textcomp}
\usepackage{bm}
\usepackage{multirow}
\usepackage[dvipsnames]{xcolor}

\usepackage[caption=false, font=footnotesize]{subfig}

\usepackage{cancel}
\usepackage{stmaryrd}

\usepackage{calrsfs}
\DeclareMathAlphabet{\pazocal}{OMS}{zplm}{m}{n}

\newcommand{\calH}{\pazocal{H}}

\newcommand{\calI}{\pazocal{I}}

\newcommand{\calJ}{\pazocal{J}}

\newcommand{\calL}{\pazocal{L}}

\newcommand{\calN}{\pazocal{N}}

\newcommand{\calP}{\pazocal{P}}

\newcommand{\calU}{\pazocal{U}}

\newcommand{\calV}{\pazocal{V}}

\newcommand{\calW}{\pazocal{W}}

\newcommand{\calX}{\pazocal{X}}

\usepackage{eufrak}

\newcommand{\Eb}{\mathbb{E}}
\newcommand{\Pb}{\mathbb{P}}

\newcommand{\Rb}{\mathbb{R}}
\newcommand{\Sb}{\mathbb{S}}
\newcommand{\Zb}{\mathbb{Z}}

\newcommand{\vzero}{{\bf 0}}

\newcommand{\bdiag}{\mathrm{bdiag}}
\newcommand{\tr}{\mathrm{tr}}

\newcommand{\Cov}{\mathrm{Cov}}

\newcommand{\vGamma}{{\mathbf{\Gamma}}}

\newcommand{\vG}{{\bf G}}

\newcommand{\vI}{{\bf I}}

\newcommand{\vA}{{\bf A}}
\newcommand{\vR}{{\bf R}}

\newcommand{\vK}{{\bf K}}

\newcommand{\vF}{{\bf F}}
\newcommand{\vM}{{\bf M}}
\newcommand{\vD}{{\bf D}}

\newcommand{\vH}{{\bf H}}

\newcommand{\vC}{{\bf C}}

\newcommand{\vQ}{{\bf Q}}
\newcommand{\vB}{{\bf B}}
\newcommand{\vS}{{\bf S}}
\newcommand{\vX}{{\bf X}}

\newcommand{\vP}{{\bf P}}

\newcommand{\vPhi}{{\bf \Phi}}

\newcommand{\argmin}{\operatornamewithlimits{argmin}}
\newcommand{\T}{^\mathrm{T}}

\newcommand{\ba}{{\bm a}}

\newcommand{\bd}{{\bm d}}

\newcommand{\bg}{{\bm g}}

\newcommand{\bp}{{\bm p}}

\newcommand{\bu}{{\bm u}}

\newcommand{\bw}{{\bm w}}
\newcommand{\bx}{{\bm x}}
\newcommand{\by}{{\bm y}}
\newcommand{\bz}{{\bm z}}

\newcommand{\bQ}{{\bm Q}}

\newcommand{\bdelta}{{\bm \delta}}

\newcommand{\bzeta}{{\bm \zeta}}
\newcommand{\blambda}{{\bm \lambda}}
\newcommand{\bmu}{{\bm \mu}}
\newcommand{\bnu}{{\bm \nu}}

\newcommand{\bSigma}{{\bm \Sigma}}

\DeclareFontFamily{OT1}{pzc}{}
\DeclareFontShape{OT1}{pzc}{m}{it}{<-> s * [1.10] pzcmi7t}{}
\DeclareMathAlphabet{\mathpzc}{OT1}{pzc}{m}{it}

\usepackage{amsthm}
\usepackage{soul}

\newtheorem{proposition}{Proposition}

\newtheorem{remark}{Remark}

\newtheorem{problem}{Problem}

\def\BibTeX{{\rm B\kern-.05em{\sc i\kern-.025em b}\kern-.08em
    T\kern-.1667em\lower.7ex\hbox{E}\kern-.125emX}}
    
\begin{document}
\title{Scaling Robust Optimization for Swarms: \\
A Distributed Perspective}
\author{Arshiya Taj Abdul*, Augustinos D. Saravanos* and Evangelos A. Theodorou
\thanks{This work was supported by the ARO Award $\#$W911NF2010151. Augustinos Saravanos acknowledges financial support by the A. Onassis Foundation Scholarship. \textit{*These authors contributed equally to this work.}}% <-this % stops a space
\thanks{Arshiya Taj Abdul and Augustinos D. Saravanos are with the School of Electrical and Computer Engineering, Georgia Institute of Technology, Atlanta, GA, 30332, USA.
        {\tt\footnotesize \{aabdul6, asaravanos\}@gatech.edu}}%
\thanks{Evangelos A. Theodorou is with the Daniel Guggenheim School of Aerospace Engineering, Georgia Institute of Technology, Atlanta, GA, 30332, USA.
        {\tt\footnotesize evangelos.theodorou@gatech.edu}}%
}
\maketitle

\begin{abstract}
This article introduces a decentralized robust optimization framework for safe multi-agent control under uncertainty.
Although stochastic noise has been the primary form of modeling uncertainty in such systems, these formulations might fall short in addressing uncertainties that are deterministic in nature or simply lack probabilistic data.
To ensure safety under such scenarios, we employ the concept of \textit{robust constraints} that must hold for all possible uncertainty realizations lying inside a bounded set. 
% Nevertheless, due to huge amount and non-convex nature of the constraints involved in safe multi-agent control, standard robust optimization approaches become intractable for such problems.
Nevertheless, standard robust optimization approaches become intractable due to the large number or non-convexity of the constraints involved in safe multi-agent control.
To address this, we introduce novel robust reformulations that significantly reduce complexity without compromising safety. 
The applicability of the framework is further broadened to address both deterministic and stochastic uncertainties by incorporating robust chance constraints and distribution steering techniques.
% However, deriving tractable versions of these constraints is challenging, thus we propose novel constraint approximations to address this.s
% To achieve scalability, we derive a distributed approach based on the Alternating Direction Method of Multipliers (ADMM), accompanied with a convergence study which handles the underlying non-convexity. 
To achieve scalability, we derive a distributed approach based on the Alternating Direction Method of Multipliers (ADMM), supported by a convergence study that accounts for the underlying non-convexity. 
In addition, computational complexity bounds highlighting the efficiency of the proposed frameworks against standard approaches are presented. 
Finally, the robustness and scalability of the framework is demonstrated through extensive simulation results across diverse scenarios, including environments with nonconvex obstacles and up to 246 agents.
% \augustinos{Good start but needs work. Beginning is a bit too wordy and some contributions are understated.}
\end{abstract}
\section{Introduction}
\label{sec:introduction}
%{\color{red} Multi-agent and swarm systems have a broad range of applications, from surveillance, monitoring,  fire detection, and inspection to border patrol and drone delivery services in e-Commerce.  
%Although in the past these systems operated in relatively sparse environments, swarm systems today have to operate in cluttered environments in the presence of heterogeneous sources of uncertainty. 
%For this reason, there is the need for planning and control frameworks that can handle the increased scale of swarm systems and have the capability to be robust to heterogeneous sources and forms of uncertainty. 

%The key concepts in our work are 5\textit{scale} and \textit{uncertainty}. 

%}
Multi-agent and swarm systems are increasingly prevalent, with various applications from surveillance, border patrol, and warehouse automation to autonomous vehicles operating in urban, marine, construction and mining environments. 
These applications require the swarm systems to operate in complex environments filled with obstacles.
% Although in the past these systems operated in relatively sparse environments, swarm systems today must operate in complex, cluttered environments filled with obstacles. 
%Solving these problems using centralized approaches becomes computationally prohibitive with the increased scale of swarm systems. 
As a result, there has been active research in developing scalable planning and control frameworks for large-scale swarm systems capable of handling such challenges \cite{saravanos2023distributed, 9812050, soria2021distributed}. 
However, real-world deployments involve various disturbances, which makes ensuring safe and reliable operation of swarm systems a critical problem.
To this end, the primary motivation of this work is to address these challenges of scalability and safety in the presence of heterogeneous forms of uncertainty.
% \subsection{Background}
% {\color{blue}
% \begin{enumerate}
%     \item Multi-agent systems and safety guarantee
%     \item Two ways of representing the uncertainty. definition and examples
%     \item Stochastic optimal control methods.
%     \item Robust optimization methods, existing methods. Challenges in solving the considered problem. Need to leverage problem structure
%     \item Having both types of uncertainty, why is it important? what are the challenges? Existing work (Kotsalis). Extension to mixed case, using RO techniques for chance constraints. 
% \end{enumerate}
% }

% \augustinos{This way of introducing our work is too much from the technical side. Before we go there, the reader should be convinced that 1) there is a great need for developing algorithms that can handle scale and uncertainty (as Prof. also emphasizes), 2) modeling uncertainty as  a stochastic signal can be insufficient (perhaps a couple of examples) and we need a stronger formulation.}

The growing need to develop scalable frameworks for safe multi-agent trajectory optimization has led many existing works to explicitly incorporate uncertainty, often modeling it as stochastic disturbances. 
This form of uncertainty, characterized by known probability distribution, is well-addressed by stochastic optimal control, a mature field encompassing methods such as multi-agent Linear Quadratic Gaussian (LQG) control \cite{alvergue2016consensus, nourian2012mean}, decentralized covariance steering \cite{saravanos2021distributed, saravanos2024distributed}, and distributed chance-constrained optimization \cite{ponda2012distributed, yang2018algorithm}. 
% Significant progress has been made in extending these techniques to multi-agent settings. 
% For instance, \cite{saravanos2021distributed, saravanos2024distributed} introduce distributed covariance steering approaches that steer the distribution of agents to meet specified covariance bounds while adhering to inter-agent coupling constraints. \cite{ponda2012distributed, yang2018algorithm} present chance-constrained optimization strategies to multi-agent task allocation problems, where constraints are enforced with a specified confidence level.  
While these frameworks are reliable, the stochastic representation cannot capture forms of disturbances that lack probabilistic distribution.
% that are encountered in practical environment. 
For instance, it cannot represent disturbances such as parametric uncertainties, modeling inaccuracies, and exogenous disturbances (such as wind gusts or ocean currents) which are unknown but are confined to lie inside a bounded uncertainty set \cite{green2012linear}. 
% A safety guarantee must be provided for all uncertainty realizations in such cases. 
This motivates us to consider a second form of uncertainty modeling- \textit{deterministic} uncertainty.  

Deterministic uncertainty is characterized as a disturbance lying inside a prespecified bounded uncertainty set. Examples of uncertainty sets include ellipsoidal sets, polytopic sets, norm-based sets \cite{nemirovski2009, bertsimas2011theory}. This uncertainty modeling is rooted in the field of Robust Control \cite{green2012linear, nemirovski2009}, which focuses on controller synthesis to guarantee performance and stability under all uncertainty realizations lying in the specified uncertainty set. This field has continuously evolved with distributed frameworks \cite{li2012distributed, wang2015distributed} developed to tackle robust control problems in multi-agent settings. Regardless, these methodologies are primarily based on stability theory, and extending these to the multi-agent trajectory optimization problem addressed in this work is challenging.

Robust Optimization (RO) is a field of mathematical programming that provides methodologies and computational tools for solving optimization problems under deterministic uncertainty \cite{ben2002robust}. The goal is to find solutions to these optimization problems that remain feasible across all possible realizations of uncertainty. As a result, these problems involve \textit{semi-infinite} constraints, which are an infinite number of constraints - each corresponding to a specific realization of uncertainty within the uncertainty set. Due to this, the optimization problems are intractable, and extensive research in RO has focused on reformulating these problems into a tractable form called robust counterparts \cite{ben2002robust, ben2009robust, bertsimas2011theory}. 
These RO techniques have found applications in diverse areas spanning from portfolio management \cite{FLIEGE2014422, 10.1007/978-981-15-9817-3_23}, power systems \cite{6850470, 8946764} to air traffic management \cite{yan2018robust, 9309827} and safe autonomy \cite{pilipovsky2024data, abdul2025nonlinear}.
While these techniques introduce flexibility, there are three key challenges in applying these techniques to our problem - addressing nonconvex constraints, computational complexity, and catering to heterogeneous forms of uncertainty. 

Nonconvex constraints, arising from obstacle and inter-agent collision avoidance, are prevalent in multi-agent trajectory planning and pose a significant difficulty in deriving robust counterparts \cite{leyffer2020survey}. Although prior work \cite{kotsalis2020convex} proposed a reformulation of quadratic constraints using S-lemma, this approach would lead to high-dimensional semi definite programming (SDP) constraints, which are computationally expensive. Since collision avoidance constraints must be satisfied at each timestep, these constraints are significant in number and grow rapidly with the time horizon. This makes even solving a single-agent trajectory optimization problem demanding, rendering scaling to multi-agent systems impractical. 
Consequently, there is a dire need to develop constraint reformulations/approximations that are computationally efficient while ensuring the safety under deterministic uncertainty.
% Consequently, there is a dire need to develop scalable and computationally efficient frameworks to guarantee safety under deterministic uncertainty. 
%
Additionally, it is not straightforward to combine RO techniques with stochastic optimal control theory to address the heterogeneous forms of uncertainty. It is essential to bridge this research gap to develop a unified scalable framework that could guarantee safety under both stochastic and deterministic uncertainty.

% \subsection{Contribution}
To this end, we integrate the techniques from robust optimization, distributed optimization, covariance steering and chance-constrained optimization to develop a scalable and robust framework to solve a multi-agent trajectory optimization problem under deterministic and stochastic uncertainty. Specifically, the  key contribution of our work is as follows
\begin{enumerate}
    \item We derive equivalent/tightly approximate tractable versions of the robust constraints for a deterministic case, significantly reducing the computational complexity compared to the existing reformulation methods without compromising the safety. We then propose extension of these approaches to handle both stochastic and deterministic uncertainty (referred to as mixed case) by deriving equivalent/approximate tractable constraint reformulations of the robust chance constraints. 
    \item We introduce a decentralized optimization framework based on Consensus ADMM \cite{boyd2011distributed} and a discounted dual update \cite{yang2022proximal} to solve the multi-agent robust optimization problem.
    \item 
    We analyze the convergence of the proposed framework to handle the non-convexity based on empirical assumptions. 
    \item We provide a theoretical computational complexity analysis highlighting the improved efficiency of the proposed framework over centralized approach and prior constraint reformulation methods.
    \item We validate the efficacy and the scalability of the proposed framework through simulation experiments in varied scenarios.  
\end{enumerate}

\noindent
\textit{\textbf{Organization of Paper:}}
We begin by introducing the multi-agent robust trajectory optimization (MARTO) problem in Section \ref{Problem Statement}. Sections \ref{sec: Reformulation of Robust Deterministic Constraints} and \ref{sec: Extension to Mixed Disturbance case} present tractable formulations of the robust constraints for the deterministic and mixed disturbance cases, respectively. Next, Section \ref{sec: distr section distr alg subsec} describes a distributed framework for solving the MARTO problem in a decentralized manner, along with a convergence analysis. Section \ref{sec: Computational Complexity Analysis} provides a computational complexity analysis. Finally, in Section \ref{sec: Simulation}, we demonstrate the effectiveness and scalability of the proposed framework through simulation experiments.

\noindent
\textit{\textbf{Notations:}}
% \subsection{Notations}
The integer set $[a,b] \cap \Zb$ is denoted with $\llbracket a, b \rrbracket$.
The space of matrices $\vX \in \Rb^{n \times n}$ that are symmetric positive (semi)-definite, i.e., $\vX \succ 0$ ($\vX \succeq 0$) is denoted with $\Sb_n^{++}$ ($\Sb_n^{+}$). 
% The inner product of two vectors $\bx, \by \in \Rb^n$ is denoted with $\langle \bx, \by \rangle = \bx\T \by$.
% = \sum_{i=1}^n x_i y_i$. 
The $\ell_2$-norm of a vector 
$\bx = 
\begin{bmatrix}
x_1 & \dots & x_n
\end{bmatrix} 
\in \Rb^n, 
$ is defined as $\| \bx \|_2 = \sqrt{\langle \bx, \bx \rangle}$.
% = \sqrt{\sum_{i=1}^n x_i^2}$. 
The weighted norm $\| \bx \|_{\vQ} = \| \bQ^{1/2} \bx \|_2 = \sqrt{\bx\T \bQ \bx}$ is also defined for any $\vQ \succ 0$. In addition, the Frobenius norm of a matrix $\vX \in \Rb^{m \times n}$ 
% = [x_{ij}] {\substack{{i \in \llbracket 1, m \rrbracket} \\j \in \llbracket 1,n \rrbracket}}$ 
%_{i \in \llbracket 1, m \rrbracket, j \in \llbracket 1,n \rrbracket} $ 
is given by $\| \vX \|_F = \sqrt{\tr(\vX\T \vX)}$. 
% = \sqrt{\sum_{i = 1}^m \sum_{j=1}^n X_{x_ij}^2}$.
%
With $[\bx_1 ; \dots; \bx_n]$, we denote the vertical concatenation of a series of vectors $\bx_1, \dots, \bx_n$. 
The cardinality of a set $\calX$ is denoted as $n(\calX)$. 
Finally, given a set $\calX$, the indicator function $\calI_{\calX}$ is defined as $\calI_{\calX}(x) = 0$, if $x \in \calX$, and $\calI_{\calX}(x) = + \infty$, otherwise. $\calN (\bmu, \bSigma)$ represents Gaussian distribution with mean '$\bmu$' and variance '$\bSigma$'.
\section{Problem Statement}\label{Problem Statement}
This section presents the multi-agent trajectory optimization problem addressed in this work. 
% We begin with the problem setup and the dynamical models of the agents. Next, we formally define the disturbances involved in the system. Following this, We will introduce an affine control policy for each agent. We will then outline the robust constraints considered in this work. Finally, we present the multi-agent robust optimization problem considered in the work.
We consider a multi-agent system of $N$ agents defined by the set $\calV = \{ 1, \dots, N \}$. All agents $i \in \calV$ may be subject to diverse dynamics, inter-agent interactions, and uncertainty. Each agent $i \in \calV$ has a set of neighbors defined by $\calN_i \subseteq \calV$, and the set of agents that consider $i$ as a neighbor by $\calP_i = \{ j \in \calV | i \in \calN_j \}$. 
% %
% \subsection{System Details}
% The system details are provided as follows
%
% \paragraph{Dynamics}
\subsection{Problem Setup}
\subsubsection{Dynamics}
We consider the following discrete-time linear time-varying dynamics for each agent $i \in \calV$, 
%
% \begin{equation}    
% \bx_{k+1}^i = \vA_k^i \bx_k^i + \vB_k^i \bu_k^i + \vC_k^i \bd_k^i + \vD_k^i \bw_k^i, ~ k \in \llbracket 0, T-1 \rrbracket, 
% \label{dynamics}
% \end{equation}
\begin{align}
    \bx_{k+1}^i & = \vA_k^i \bx_k^i + \vB_k^i \bu_k^i + \vC_k^i \bd_k^i + \vD_k^i \bw_k^i, ~ k \in \llbracket 0, T-1 \rrbracket, \nonumber
    \\
    \bx_0^i & = \bar{\bx}_0^i + \bar{\bd}_0^i + \bar{\bw}_0^i,
    \label{dynamics}
\end{align}
where $\bx_k^i \in \Rb^{n_{x_i}}$ is the state, $\bu_k^i \in \Rb^{n_{u_i}}$ is the control input at time $k$, $\bx_0^i$ is initial state with $\bar{\bx}_0^i$ as its known part, and $T$ is the time horizon. The terms $\bar{\bd}_0^i \in \Rb^{n_{x_i}}, \bd_k^i \in \Rb^{n_{d_i}}$ represent \textit{deterministic} uncertainty, while the terms $ \bar{\bw}_0^i \in \Rb^{n_{x_i}}, \bw_i^k \in \Rb^{n_{w_i}}$ refer to \textit{stochastic} uncertainty. 
% Each initial state $x_0^i, ~ i \in \calV$, is of the following form
%
% \begin{equation}    
% \bx_0^i = \bar{\bx}_0^i + \bar{\bd}_0^i + \bar{\bw}_0^i,
% \end{equation}
% 
% where $\bar{\bx}_0^i$ is the known part, while $\bar{\bd}_0^i$ and $\bar{\bw}_0^i$ refer to the deterministic and stochastic parts, respectively. 
Finally, $\vA_k^i, \vB_k^i, \vC_k^i, \vD_k^i$ are the dynamics matrices of appropriate dimensions.

For convenience, let us also define the sequences 
$\bx^i = [ \bx_0^i; \dots; \bx_T^i ] \in \Rb^{(T+1) n_{x_i}}$, 
$\bu^i = [\bu_0^i; \dots; \bu_{T-1}^i] $ $\in \Rb^{T n_{u_i}}$, 
$\bw^i = [\bar{\bw}_0^i; \bw_0^i; \dots; \bw_{T-1}^i] \in \Rb^{n_{x_i} + T n_{w_i}}$,
$\boldsymbol{\zeta}^i = [\bar{\bd}_0^i; \bd_0^i; \dots; \bd_{T-1}^i] \in \Rb^{n_{x_i} + T n_{d_i}}$.
%
% and the matrix $\vPhi^i(k_1,k_2) = A_{k_1 -1}^i A_{k_1-2}^i \dots A_{k_2}^i$ for $ k_1 > k_2$.
%
The dynamics \eqref{dynamics} can then be rewritten in a more compact form as 
\begin{equation}
\bx^i = \vG_0^i \bar{\bx}_0^i + \vG_u^i \bu^i + \vG_{\bw}^i \boldsymbol{\bw}^i + \vG_{\zeta}^i \boldsymbol{\zeta}^i,
\label{compact dynamics}
\end{equation} 
where, for 
the matrix $\vPhi^i(k_1,k_2)$ defined as 
$ \vPhi^i(k_1,k_2) = A_{k_1 -1}^i A_{k_1-2}^i \dots A_{k_2}^i$ for $ k_1 > k_2$, we have 
\begin{subequations}
\begin{align}
&
\vG_0^i = 
\begin{bmatrix}
\vI; & 
\vPhi^i(1,0) ; & 
\vPhi^i(2,0) ; &
\dots ; &
\vPhi^i(T,0)
\end{bmatrix}, 
%\in \Rb^{(T+1) n_x \times n_x }, 
\nonumber
\\[0.2cm]
&
\bar{\vG}^i (B^i) =  \begin{bmatrix}
\vzero & \vzero & \dots & \vzero \\ 
B_0^i & \vzero & \dots & \vzero \\ 
\vPhi^i(2,1) B_0^i & B_1^i & \dots & \vzero \\ 
\vdots & \vdots & \vdots & \vdots\\ 
\vPhi^i(T,1) B_0^i & \vPhi^i(T,2) B_1^i & \dots & B_{T-1}^i
\end{bmatrix}
\nonumber
\end{align}
\end{subequations}
such that $\vG_u^i = \bar{\vG}^i ( B^i )$, $\vG_{\zeta}^i = \begin{bmatrix}
    \vG_0^i, & \bar{\vG}^i (C^i )
\end{bmatrix}$, and
$\vG_{\bw}^i = \begin{bmatrix}
    \vG_0^i, & \bar{\vG}^i (D^i ) \end{bmatrix}$.
 
% where the matrices $\vG_0^i$, $\vG_u^i$, $\vG_{\bw}^i$, and $\vG_{\zeta}^i$ are defined in the Appendix\ref{Appendix: Compact Dynamics Matrices}.

% 
% \paragraph{Characterization of Uncertainty}
% \label{sec: char of unc}
% We consider the existence of two types of uncertainty - exogenous deterministic and stochastic uncertainties. 
% The two types of uncertainty are formally described as follows
% \paragraph{Exogenous deterministic uncertainty}
%
% \subsubsection{Exogenous deterministic uncertainty} 
 % The unknown deterministic disturbances $\boldsymbol{\zeta}^i$ typically stem from external factors that cannot be accurately represented through stochastic signals. This creates a requirement to formulate optimization problems and achieve optimal policies that are feasible in a \textit{robust} sense for all possible values of such disturbances in a bounded \textit{uncertainty set}. Examples of such uncertainty sets include ellipsoids, polytopes, and others \cite{RobustNemirovski}. 

\subsubsection{Uncertainty Characterization}
We consider the deterministic uncertainty $\boldsymbol{\zeta}^i$ to be lying inside a bounded ellipsoidal uncertainty set defined as follows
 \begin{equation}
\begin{aligned}
    \calU_i[\tau^i] = 
    \big\{
    \boldsymbol{\zeta}^i \in \Rb^{n_{\zeta_i}} &| 
    \;	\exists (\bz_i \in \Rb^{\bar{n}_i}, \tau^i \in \Rb): 
    \; 
    \\ 
    & \boldsymbol{\zeta}^i = \boldsymbol{\Gamma}_i \bz_i, 
    \;
    \bz_i^T \vS_i  \bz_i \leq \tau^i
    \big\}, 
\end{aligned}
\label{ellipsoid uncertainty set}
\end{equation}
where $n_{\zeta_i} = n_{x_i} + T n_{d_i}$, $\boldsymbol{\Gamma}_i \in \Rb^{n_{\zeta_i} \times \bar{n}_i}$, $\vS_i \in \Sb^{++}_{\bar{n}_i}$ and $\tau^i >0$ is the uncertainty level.
We consider these sets due to their extensive usage in robust control applications \cite{Petersen2000}. Further, ellipsoid sets can be used to approximate more complex uncertainty sets. 
%
% \paragraph{Stochastic uncertainty}
% \subsubsection{Stochastic uncertainty} 
% The unknown stochastic component $\bw^i$ is assumed to be a  Gaussian random vector with zero mean ($\Eb[\bw^i] = 0$), and covariance $\bSigma_{\boldsymbol{\bw}^i} \in \Sb^{+}_{n_{x_i} + T n_{w_i}}$.

The stochastic uncertainty $\bw^i$ is considered to be a Gaussian random vector with zero mean ($\Eb[\bw^i] = 0$), and covariance $\bSigma_{\boldsymbol{\bw}^i} \in \Sb^{+}_{n_{x_i} + T n_{w_i}}$.
%
%
%

% \subsubsection{Control Policies}
% \paragraph{Control Policies}
%

We consider affine \textit{purified state} feedback policies \cite{nemirovski2009, kotsalis2020convex} for all agents. Let $\{ \hat{\bx}_k^i \}_{k = 0}^{T}$ be the disturbance-free states whose dynamics are given by
\begin{align}
\hat{\bx}_{k+1}^i & = \vA_k^i \hat{\bx}_t^i + \vB_k^i \bu_k^i, \quad k \in \llbracket 0,T-1 \rrbracket,
\\
\hat{\bx}_0^i & = \bar{\bx}_0^i.
\end{align}
Thereby, the purified states $\bdelta_k^i$ are defined as $\bdelta_k^i = \bx_k^i - \hat{\bx}_k^i$, 
% Subsequently, we define the purified states $\bdelta_k^i = \bx_k^i - \hat{\bx}_k^i$, 
whose dynamics are given by
\begin{equation}
\begin{aligned}
\bdelta_{k+1}^i & = \vA_k^i \bdelta_k^i + \vC_k^i \bd_k^i + \vD_k^i \bw_k^i 
\quad k \in \llbracket 0,T-1 \rrbracket, \\
\bdelta_0^i & = \bar{\bd}_0^i + \bar{\bw}_0^i,
\end{aligned}
\label{purified state - compact dynamics}
\end{equation}
or more compactly by $\bdelta^i = \vG_{\zeta}^i \boldsymbol{\zeta}^i + \vG_{\bw}^i \boldsymbol{\bw}^i$, 
%
% \begin{equation}
% \bdelta^i = \vG_{\zeta}^i \boldsymbol{\zeta}^i + \vG_{\bw}^i \boldsymbol{\bw}^i, 
% \label{purified state - compact dynamics}
% \end{equation}
%
where $\bdelta^i = \begin{bmatrix}
\bdelta_0^i; & \dots; & \bdelta_T^i
\end{bmatrix}$. 
We consider the following affine purified state feedback control policy for each agent $i$,
\begin{equation}
\bu^i_k = \bar{\bu}^i_k + \sum_{\ell = k - \gamma_h +1}^k \vK^i_{k,\ell} \bdelta^i_\ell, \; 
\forall \; k \in \llbracket 0,T-1 \rrbracket,
\end{equation}
where $\bar{\bu}^i_k \in \Rb^{n_{u_i}}$ are the feed-forward control inputs, and $\vK^i_{k,\ell} \in \Rb^{n_{u_i} \times n_{x_i}}$ are feedback gains on the purified states,  $\gamma_h$ denotes the length of the history interval considered for the feedback. The above policies can be written in compact form (when $\gamma_h = T$) as
% We consider the following affine purified state feedback control policy for each agent $i$, 
% with $\bar{u}^i_k \in \Rb^{n_{u_i}}$ as the feed-forward control inputs, and $K^i_{k,\ell} \in \Rb^{n_{u_i} \times n_{x_i}}$ as feedback gains on the purified states.
%
\begin{equation}
\bu^i = \bar{\bu}^i + \vK^i \bdelta^i,
\label{control expression}
\end{equation}
where $\bar{\bu}^i \in \Rb^{T n_{u_i}}$, and $\vK^i \in \Rb^{ T n_{u_i} \times (T+1)n_{x_i} }$ are defined as
\begin{align}
\bar{\bu}^i = \begin{bmatrix}
\bar{\bu}_0^i; & \dots; & \bar{\bu}_{T-1}^i
\end{bmatrix}, 
\nonumber
\end{align}
\begin{align}
\vK^i = 
\begin{bmatrix}
\vK_{0,0}^i & \vzero & \dots & \vzero & \vzero \\ 
\vK_{1,0}^i & \vK_{1,1}^i & \dots & \vzero & \vzero \\  
\vdots & \vdots & \vdots & \vdots & \vdots\\ 
\vK_{T-1,0}^i & \vK_{T-1,1}^i & \dots & \vK_{T-1, T-1}^i & \vzero
\end{bmatrix}.
\nonumber
\end{align}
\subsection{Problem Formulation} \label{sec: Robust Constraints}
% In this work, we examine two cases: the \textit{deterministic case} and the \textit{mixed case}. First, we analyze the deterministic case, where the system is influenced solely by deterministic disturbances $\{ \bzeta^i \}_{i \in \calV}$ (i.e., $\vG_{\bw}^i = 0$ for all $i \in \calV$). Next, we consider the mixed case, where both deterministic and stochastic disturbances affect the system. Accordingly, we present robust constraints considered in each case. 
In this work, we consider two cases: the \textit{deterministic case}, where the system is affected solely by deterministic disturbances $\{ \bzeta^i \}_{i \in \calV}$ (i.e., $\vG_{\bw}^i = 0$ for all $i \in \calV$), and the mixed case, where both deterministic and stochastic disturbances affect the system. Accordingly, we present robust constraints considered in each case. 
We define the constraints using the following matrices and vectors: $\vH^i \in \Rb^{n_{\text{pos}} \times n_{x_i}}$, $\ba_i \in \Rb^{(T+1)n_{x_i}}$, $b_i \in \Rb$, $\bm{p}_{\text{obs}}, \bm{p}_{\text{tar}} \in \Rb^{n_{\text{pos}}}$, and $c_i, c_{i,j} \in \Rb^+$. In multi-agent trajectory optimization setting, the terms $\bm{p}_{\text{obs}}$ and $\bm{p}_{\text{tar}} $ denote position of obstacle and target respectively, while the terms $c_i$ and $c_{i,j}$ refer to distance thresholds.
% {\color{red} In the multi-agent and control of swarm fomulation   $\bm{p}_{\text{obs}} $ denotes ...,  }
%
\subsubsection{Deterministic Case} \label{sec: Robust Constraints for deterministic case} We consider the following robust constraint forms for each agent $i \in \calV$ -
\begin{enumerate}
    \item \textit{Linear constraints:}
    \begin{equation}
    \ba_i^T \bx^i \leq b_i, \quad \forall \bzeta^i \in \calU_i.
    \label{Robust Linear Constraints}
    \end{equation}
    %
    % \item Nonlinear Constraints:
    % \begin{itemize}
        \item \textit{Nonconvex constraints:} At every time step $k$,
         \begin{equation}
        \| \vH^i \bx_k^i - \bm{p}_{\text{obs}} \|_2 \geq c_i, \quad \forall \bzeta^i \in \calU_i.
        \label{Robust nonconvex norm-of-mean constraints}
        \end{equation}
        \item \textit{Nonconvex inter-agent constraints:} For every neighboring agent $j \in \calN_i$, and at every time step $k$,
        \begin{equation}
        \begin{aligned}
            & \| \vH^i \bx^i_k - \vH^j \bx^j_k \|_2 \geq c_{ij}, 
            % \\
            % &~~~~~~~~~~~~~~~~~~~~~~~~~~
            ~~~
            \forall \bzeta^i \in \calU_i, ~ \bzeta^j \in \calU_j.
        \end{aligned}
        \label{robust inter-agent collision avoidance}
        \end{equation}
        Note that the matrices $\vH^i$ and $\vH^j$ could be different.
        \item \textit{Convex norm constraints:} At some time instant $\tilde{k}$,
        \begin{align}
        \| \vH^i \bx_{\tilde{k}}^i - \bm{p}_{\text{tar}} \|_2 \leq c_i, \quad \forall \bzeta^i \in \calU_i.
        \label{Robust convex norm-of-mean constraints}
        \end{align}
        % \end{itemize}
\end{enumerate} 
\subsubsection{Mixed Case} \label{sec: robust constraints in mixed case}
The corresponding robust constraint forms considered in mixed case for each agent $i \in \calV$, with the probabilities represented by $p_i, p_{ij}$, are as follows 
\begin{enumerate}
    \item \textit{Linear chance constraints:}
    \begin{equation}
    \Pb (\ \ba_i \T \bx^i  > b_i ) \leq p_i, 
    \quad \forall \bzeta ^i \in \calU_i.
    \label{Mixed Constraints: Robust linear Chance Constraints}
    \end{equation} 
    %
    % \item Robust Nonlinear Chance Constraints:
    % \begin{itemize}
        \item \textit{Nonconvex chance constraints:} At every time step $k$,
        \begin{equation}
        \Pb ( \| \vH^i \bx_k^i - \bm{p}_{\text{obs}} \|_2 \geq c_i) \geq 1-p_i, ~~~ \forall \bzeta^i \in \calU_i.
        \label{Mixed Constraints: nonconvex norm chance}
        \end{equation}
        \item \textit{Nonconvex inter-agent chance constraints:}
        For every neighboring agent $j \in \calN_i$, and at every time step $k$,
        \begin{equation}
        \begin{aligned}
            & \Pb( \| \vH^i \bx^i_k - \vH^j \bx^j_k \|_2 \geq c_{ij}) \geq 1-p_{ij}, 
            \\
            &~~~~~~~~~~~~~~~~~~~~~~~~~~~~~~
            \forall \bzeta^i \in \calU_i, ~ \bzeta^j \in \calU_j
        \end{aligned}
        \label{Mixed Constraints: inter-agent nonconvex chance}
        \end{equation}
        \item \textit{Convex norm chance constraints:} At some time step $\tilde{k}$
        \begin{equation}
        \Pb (\| \vH^i \bx_{\tilde{k}}^i - \bm{p}_{\text{tar}} \|_2 \leq c_i) \geq 1-p_i, \quad \forall \bzeta^i \in \calU_i,
        \label{Mixed Constraints: convex norm chance}
        \end{equation}
    % \end{itemize}
    %
    \item \textit{Covariance constraints:}
    At some time step $\tilde{k}$,
    \begin{equation}
    \Cov(\bx^i_{\tilde{k}}) \preceq \bSigma^i_{\tilde{k}}, \quad \forall \bzeta^i \in \calU_i,
    \label{Mixed Constraints: Covariance}
    \end{equation}
    where $\bSigma^i_{\tilde{k}}$ is the allowed state covariance upper bound.
    % \item \textit{Convex expectation of norm constraints:} At some time step $\tilde{k}$,
    % \begin{equation}
    %     \Eb [\|  \vH^i \bx_{\tilde{k}}^i - \bm{p}_{\text{tar}} \|_2^2 ] \leq c_i^2, \quad \forall \boldsymbol{\zeta}^i \in \calU_i.
    %     \label{Mixed Constraints: Expectation of norm constraints}
    % \end{equation}
\end{enumerate}

We consider the following multi-agent robust trajectory optimization problem in this work.
% \subsection{Problem Formulation}
\begin{problem}[Multi-Agent Robust Trajectory Optimization (MARTO)]
\label{initial problem}
For all agents $i \in \calV$, find the robust optimal policies $\bar{\bu}^i,  \vK^i$, such that 
\begin{align}
&  {~~~~~} 
% \{ \bar{\bu}^i, \vK^i \}_{i \in \calV} = \argmin
% \sum_{i \in \calV} 
% J_i(\bar{\bu}^i, \vK^i)
\min_{ \{ \bar{\bu}^i, \vK^i \}_{i \in \calV} }
\quad \sum_{i \in \calV} 
J_i(\bar{\bu}^i, \vK^i)
%\quad (= J (\bu^1, \bu^2, \dots \bu^N) ) 
\nonumber
\\[0.1cm]
\mathrm{s.t.} \; 
&~
\bx_{k+1}^i = \vA_k^i \bx_k^i + \vB_k^i \bu_k^i + \vC_k^i \bd_k^i + \vD_k^i \bw_k^i,
\nonumber
\\
&~
\bar{\bx}_0^i \text{ : given}, 
  \label{robust_optimization_problem_1}
 \\ 
&
\eqref{Robust Linear Constraints}  - \eqref{Robust convex norm-of-mean constraints}  
\quad \text{or} \quad
\eqref{Mixed Constraints: Robust linear Chance Constraints} - \eqref{Mixed Constraints: Expectation of norm constraints}
\nonumber
\end{align}  
where each cost component $J_i(\bar{\bu}^i, \vK^i) = \bar{\bu}^i{}\T \vR_{\bar{u}} \bar{\bu}^i 
+ || \vR_K \vK^i||_F^2$, with $\vR_{\bar{u}}, \vR_K \in \Sb^{+}_{Tn_{u_i}}$, penalizes the control effort of agent $i \in \calV$. We consider the set of constraints (\eqref{Robust Linear Constraints} - \eqref{Robust convex norm-of-mean constraints}) for the deterministic case; or (\eqref{Mixed Constraints: Robust linear Chance Constraints} - \eqref{Mixed Constraints: Covariance}) for the mixed case. Additionally, we can combine the constraints from both cases as required. 
\end{problem}

It is not straightforward to solve the above problem in a distributed manner due to two significant issues. First, the robust constraints must be satisfied for all realizations of deterministic uncertainty, resulting in infinite number of constraints and making the problem intractable.
Second, the coupling constraints \eqref{robust inter-agent collision avoidance} (in deterministic case) and \eqref{Mixed Constraints: nonconvex norm chance} (in mixed case) are between the states of the agents, which are dependent on disturbances. This makes the exchange of state variables among agents challenging. We address these challenges in the subsequent sections.
%

% We develop a framework to solve Problem \ref{initial problem} through the following sequence of steps. First, we derive equivalent/approximate tractable versions of the robust constraints employing RO techniques, thereby transforming Problem \ref{initial problem} into a tractable problem. Second, by taking advantage of the tractable constraints, we characterize the state of the agent to enable the exchange of state information among agents. Finally, we present a distributed framework to solve Problem \ref{initial problem} in a decentralized way.
%
\section{Reformulation of Robust Constraints For Deterministic Case}
% \subsection{Reformulation of Semi-infinite Constraints}
\label{sec: Reformulation of Robust Deterministic Constraints}
In this section, we derive the equivalent or approximate tractable versions of the semi-infinite robust constraints defined in the section \ref{sec: Robust Constraints for deterministic case}. 
For that, we begin by rewriting the dynamics \eqref{compact dynamics} of each agent $i \in \calV$, in terms of the control variables $\bar{\bu}^i$ and $\vK^i$ by using \eqref{control expression} as follows
\begin{align}
    & \bx^i = \vG_0^i \bar{\bx}_0^i + \vG_u^i \bar{\bu}^i + (\vG_u^i \vK^i + \vI) \bdelta^i.
    \label{dynamics - in terms of deltai}
\end{align}
Recall that $\vG_{\bw}^i = 0$ for this case, thus we can rewrite the above using \eqref{purified state - compact dynamics} as follows
\begin{equation}
    \bx^i = \bmu_{x,\bar{u}}^i + \vM_i \boldsymbol{\zeta}^i,
    \label{dynamics - deterministic case}
\end{equation}
where $\bmu_{x,\bar{u}}^i = \vG_0^i \bar{\bx}_0^i + \vG_u^i \bar{\bu}^i$, and $\vM_i =  (\vG_u^i \vK^i + \vI) \vG_{\zeta}^i$.
% \begin{equation}
% \bmu_{x,\bar{u}}^i = \vG_0^i \bar{\bx}_0^i + \vG_u^i \bar{\bu}^i, 
% \quad
% \vM_i =  (\vG_u^i \vK^i + \vI) \vG_{\zeta}^i.
% \nonumber
% \end{equation}
%DIF > 
For simplicity, we define a matrix $\vP_k^i \in \Rb^{n_{x_i} \times (T+1)n_{x_i}}$ such that $\bx_k^i = \vP_k^i \bx^i$. Further, we also define the variables $\bmu_{x_k,\bar{u}}^i = \vP^i_k \bmu_{x,\bar{u}}^i$, $\tilde{\vM}_k^i = \vP^i_k \vM_i$, $\hat{\bmu}^i_{\text{pos},k} = \vH^i \bmu_{x_k,\bar{u}}^i$, and $\vM_{\text{det}}^{k,i} (\bzeta^i) = \vH^i \tilde{\vM}_k^i \bzeta^i $.
%   
% and the variables $\bmu_{x_k,\bar{u}}^i$, $\tilde{\vM}_k^i$ and $\vM_{\text{det}}^{k,i} (\bzeta^i) = \vH^i \tilde{\vM}_k^i \bzeta^i $ given as $\tilde{\vM}_k^i = \vP^i_k \vM_i$, and $\bmu_{x_k,\bar{u}}^i = \vP^i_k \bmu_{x,\bar{u}}^i$.
%
Subsequently, we reformulate each robust constraint presented in Section \ref{sec: Robust Constraints for deterministic case}.
\subsubsection{Robust Linear Constraints} 
%
% Let us rewrite the linear constraint \eqref{Robust Linear Constraints} as follows
% \begin{align}
%     \ba_i \T \bx^i \leq b_i, \quad \forall \bzeta^i \in \calU_i.
%     \nonumber
% \end{align}
% 
The linear constraint \eqref{Robust Linear Constraints} can be equivalently given as 
\begin{align}
    \max_{\bzeta^i \in \calU_i} \ba_i \T \bx^i \leq b_i,
    \nonumber
\end{align}
which can be simplified using \eqref{dynamics - deterministic case} to get the following
\begin{align}
    \ba_i \T \bmu_{x,\bar{u}}^i 
    + \max_{\bzeta^i \in \calU_i} \ba_i \T \vM_i \boldsymbol{\zeta}^i \leq b_i. 
    \label{single mean constraint max}
\end{align}
%DIF > 
%
We will now derive a closed form solution to the sub-problem $\max_{\bzeta^i \in \calU_i} \ba_i \T \vM_i \boldsymbol{\zeta}^i$ in the following proposition.
\begin{proposition} \label{prop: Max bound on mean} 
The maximum and the minimum values of $\tilde{\ba}_i \T \bzeta^i$ when the uncertainty vector $\bzeta^i \in \calU_i[\tau^i]$ are given by
\begin{equation}
    \max_{\bzeta^i \in \calU_i} \tilde{\ba}_i \T \bzeta^i 
    =
    - \min_{\bzeta^i \in \calU_i} \tilde{\ba}_i \T \bzeta^i
    = 
    \sqrt{\tau^i} \| \boldsymbol{\Gamma}_i \T \tilde{\ba}_i \|_{\vS_i^{-1}}.
    \label{max bound on mean problem}
\end{equation}
\end{proposition}
\proof
The proof is based on the first-order optimality conditions of the problem, and provided in Appendix\ref{Appendix: Sec1 linear constraint proof}.
%\qed{}
\endproof
Using Proposition \ref{prop: Max bound on mean}, an equivalent tractable version of the constraint \eqref{single mean constraint max} can be given as follows
\begin{equation}
    \ba_i \T \bmu_{x,\bar{u}}^i + \sqrt{\tau^i} || \boldsymbol{\Gamma}_i \T \vM_i \T \ba_i||_{\vS_i^{-1}}
    \leq b_i.
    \label{mean_constraint_form1}
\end{equation}
Similarly, the minimum value can be provided as
\begin{align}
    & \min_{\bzeta^i \in \calU_i} \ba_i \T \bx^i = \ba_i \T \bmu_{x,\bar{u}}^i - \sqrt{\tau^i} || \boldsymbol{\Gamma}_i \T \vM_i \T \ba_i||_{\vS_i^{-1}}.
    \nonumber
\end{align}
Based on the above observations, we introduce the following additional form of robust linear constraints
\begin{equation}
\begin{aligned}
     \bar{\ba}_i \T \bmu_{x,\bar{u}}^i = \bar{b}_i,
    \quad \sqrt{\tau^i} || \boldsymbol{\Gamma}_i \T \vM_i \T \bar{\ba}_i||_{\vS_i^{-1}} \leq \epsilon_i^{b},
    \label{New mean formulation}
\end{aligned}   
\end{equation}
with $\bar{\ba}_i \in \Rb^{ (T+1) n_{x_i}}, \bar{b}_i \in \Rb$ and $\epsilon_i^{b} \in \Rb^+$, which would imply the following set of constraints
\begin{align}
 \bar{\ba}_i \T \bx^i \leq \bar{b}_i + \epsilon_i^{b}, 
 \quad
 \bar{\ba}_i \T \bx^i \geq \bar{b}_i - \epsilon_i^{b}, \quad \forall \bzeta^i \in \calU_i.   \nonumber
\end{align}
The constraint form \eqref{New mean formulation} allows us to separately control the disturbance-free state component $\bmu_{x,\bar{u}}^i$ with $\bar{\bu}^i$, and the uncertain component $\vM_i \bzeta^i$ with $\vK^i$.
%
% \subsection{Robust Nonlinear Constraints}
% %
% In this subsection, we derive the equivalent/approximate tractable versions of the robust non-linear constraints defined by \eqref{Robust nonconvex norm-of-mean constraints} - \eqref{Robust convex norm-of-mean constraints}  in the section \ref{sec: Robust Constraints for deterministic case}.
\subsubsection{Nonconvex Norm Constraints}
We present the following two approaches for deriving the equivalent/approximate tractable constraints for the nonconvex constraint \eqref{Robust nonconvex norm-of-mean constraints}.
\paragraph{Initial Approach}\label{initial approach}
The robust constraint in \eqref{Robust nonconvex norm-of-mean constraints} can be rewritten as 
\begin{align}
     \min_{\bzeta^i \in \calU_i} \| \vH^i \bx_k^i - \bm{p}_{\text{obs}} \|_2^2 \geq c_i^2.
     \label{obstacle constraint initial approach initial}
\end{align}
Although the problem in the LHS is convex, the equivalent tractable constraints to the above, obtained by applying S-lemma \mbox{
\cite{kotsalis2020convex}}\hskip0pt
, might not be convex due to the non-convexity of the overall constraint. Thus, we first linearize the constraint \eqref{obstacle constraint initial approach initial} with respect to control parameters 
$\Bar{\bu}^i, \vK^i$
around some nominal values $\Bar{\bu}^{i}_{\text{lin}}, \vK^{i}_{\text{lin}}$ to obtain the following constraint.
\begin{equation}    
\begin{aligned}
\boldsymbol{\zeta}^i {}\T \mathcal{Q}^{\text{lin}}( \vK^i)
\boldsymbol{\zeta}^i 
+ 2 \bar{\mathcal{Q}}^{\text{lin}} ( \vK^i, \Bar{\bu}^i)  \boldsymbol{\zeta}^i 
% \\
% & \qquad \qquad \qquad \qquad  
+ q^{\text{lin}} (\Bar{\bu}^i)
\geq c_i^2, \; \forall \boldsymbol{\zeta}^i \in \calU_i,
\nonumber
\end{aligned}
\end{equation} 
%
% where $\mathcal{Q}^{\text{lin}} ( \vK^i), \bar{\mathcal{Q}}^{\text{lin}} ( \vK^i, \Bar{\bu}^i),$ and $q^{\text{lin}}(\Bar{\bu}^i)$ are defined in the Appendix\ref{Appendix: Sec 1 initial approach}.
where $\mathcal{Q}^{\text{lin}} ( \vK^i) = (\vH_{\text{pos}}^i \tilde{\vM}_k^{i,\text{lin}}) \T \vH_{\text{pos}}^i \tilde{\vM}_k^i + (\vH_{\text{pos}}^i (\tilde{\vM}_k^i - \tilde{\vM}_k^{i,\text{lin}})) \T \vH_{\text{pos}}^i \tilde{\vM}_k^{i,\text{lin}}$;
$\bar{\mathcal{Q}}^\text{lin} ( \vK^i, \Bar{\bu}^i)  = (\vH_{\text{pos}}^i \bmu_{x_k,\bar{u}}^{i,\text{lin}} - \bm{p}_{\text{obs}}) \T (\vH_{\text{pos}}^i \tilde{\vM}_k^i) + (\vH_{\text{pos}}^i(\bmu_{x_k,\bar{u}}^i - \bmu_{x_k,\bar{u}}^{i,\text{lin}})) \T (\vH_{\text{pos}}^i \tilde{\vM}_k^{i,\text{lin}})$; and
$q^{\text{lin}}(\Bar{\bu}^i) = (\vH_{\text{pos}}^i \bmu_{x_k,\bar{u}}^{i,\text{lin}} - \bm{p}_{\text{obs}}) \T (\vH_{\text{pos}}^i (2\bmu_{x_k,\bar{u}}^{i}-\bmu_{x,\bar{u}}^{i,\text{lin}})  - \bm{p}_{\text{obs}} )$. 
The equivalent tractable constraints to the above constraint can be given as \cite{kotsalis2020convex} -
\begin{equation}
   \begin{aligned}
    & \beta \in \Rb, \; \beta \geq 0 \\
    &
    \begin{bmatrix}
        \boldsymbol{\Gamma}_i \T \mathcal{Q}^{\text{lin}} (\vK^i) \boldsymbol{\Gamma}_i 
        + \beta \vS_i 
        & \boldsymbol{\Gamma}_i \T \bar{\mathcal{Q}}^{\text{lin}} ( \vK^i, \Bar{\bu}^i)^T  \\
       \bar{\mathcal{Q}}^{\text{lin}} ( \vK^i, \Bar{\bu}^i) \boldsymbol{\Gamma}_i
        & q^{\text{lin}} (\Bar{\bu}^i) - c_i^2 - \beta \tau^i
    \end{bmatrix}
     \succeq 0
\end{aligned}
\label{obstacle_semidef_old}
\end{equation}
The above semi-definite constraint involves a matrix of size $\bar{n}_i + 1$ (where $\boldsymbol{\Gamma}_i = I$,  $\bar{n}_i = n_{x_i} + T n_{d_i}$).  
As a result, the computational cost of the formulation increases exponentially with the number of constraints. Since the collision avoidance constraints need to be satisfied at every time step for each obstacle per agent, the above formulation is unsuitable for such cases. To address this issue, we propose a second approach to derive an approximate tractable version of the constraint \eqref{Robust nonconvex norm-of-mean constraints}.
\paragraph{Reducing computational complexity}
\label{new reformulation}
Let us first write the constraint \eqref{Robust nonconvex norm-of-mean constraints} in terms of $\bzeta^i$ as follows -
% \begin{align}
%     \| \vH^i( \bmu_{x_k,\bar{u}}^i +  \tilde{\vM}_k^i \bzeta^i) - \bm{p}_{\text{obs}} \|_2 
%     \geq c_i, \quad \forall \bzeta^i \in \calU_i.
%     \nonumber
% \end{align}
% %
% For simplicity, let us define $\hat{\bmu}^i_{\text{pos},k} = \vH^i \bmu_{x_k,\bar{u}}^i$. Using this, we can rewrite the above constraint as 
\begin{align}
    \| \hat{\bmu}^i_{\text{pos},k} - \bm{p}_{\text{obs}}
    +  \vM_{\text{det}}^{k,i} (\bzeta^i) \|_2 
    \geq c_i, \quad \forall \bzeta^i \in \calU_i. \nonumber
\end{align}
Based on the triangle inequality, we consider the following tighter approximation of the above constraint  
% \begin{align}
%     \| \hat{\bmu}^i_{\text{pos},k} - \bm{p}_{\text{obs}} \|_2 
%     - \| \vH^i \tilde{\vM}_k^i \bzeta^i \|_2
%     \geq c_i, \quad \forall \bzeta^i \in \calU_i
%     \nonumber
% \end{align}
\begin{align}
    \| \hat{\bmu}^i_{\text{pos},k} - \bm{p}_{\text{obs}} \|_2 
    - \|\vM_{\text{det}}^{k,i} (\bzeta^i) \|_2
    \geq c_i, \quad \forall \bzeta^i \in \calU_i
    \nonumber
\end{align}
which can be rewritten as 
\begin{align}
    \| \hat{\bmu}^i_{\text{pos},k} - \bm{p}_{\text{obs}} \|_2 \geq
    c_i 
    + \max_{\bzeta^i \in \calU_i} \| \vM_{\text{det}}^{k,i} (\bzeta^i) \|_2.
\end{align}
Let us now introduce a slack variable $\tilde{c}_i^k$ such that we split the above constraint into the following system of two constraints %
\begin{align}
    & \| \hat{\bmu}^i_{\text{pos},k} - \bm{p}_{\text{obs}} \|_2 
    \geq  \tilde{c}_i^k + c_i,
    \label{nonconvex obstacle final ubar}\\
    & 
    \max_{\bzeta^i \in \calU_i} \|\vM_{\text{det}}^{k,i} (\bzeta^i) \|_2 
    \leq \tilde{c}_i^k 
   \label{nonconvex obstacle K1}
\end{align}
%
% Note, however, that the constraint \eqref{nonconvex obstacle final ubar}  is non-convex. We address this by linearizing the constraint around the nominal trajectory $\bmu_{x,\bar{u}}^{i,l}$.
% In addition, 
The constraint \eqref{nonconvex obstacle K1} is still a semi-infinite constraint; thus, we need to derive its equivalent tractable constraints. Since the LHS in the constraint \eqref{nonconvex obstacle K1} involves a non-convex optimization problem, deriving the equivalent tractable constraints is NP-hard. Hence, we derive tighter approximations of the constraint \eqref{nonconvex obstacle K1} in the following proposition. 
\begin{proposition} \label{prop- obs avoidance}
    The tighter approximate constraints to the semi-infinite constraint \eqref{nonconvex obstacle K1} are given by the following second-order cone constraints (SOCP constraints)
\begin{align}
        &
        \| \bmu_{d,k}^i \|_2 \leq \tilde{c}_i^k
        \label{semi-definite obstacle new} \\
        &
        (\bmu_{d,k}^i)_{\bar{m}} \geq  
        \sqrt{\tau^i} || \boldsymbol{\Gamma}_i \T \vM_i \T \bm{h}_{k,\bar{m}}^i||_{\vS_i^{-1}}
        \label{mu_d expression}
\end{align}
where $\bm{h}_{k,\bar{m}}^i {}\T $ is $\bar{m}^{th}$ row of $ \vH^i \vP_k^i$.
\end{proposition}
\proof
% A detailed proof is provided in \cite{abdul2024scaling}. 
A detailed proof is provided in Appendix\ref{Appendix: Sec 1 nonconvex obstacle mu_d}.
% The outline of the proof is as follows. First, we obtain a tighter upper bound of the term $\| \vH^i \tilde{\vM}_k^i \bzeta^i \|_2 $ by introducing a variable $\bmu_{d,k}^i \geq 0$ which upper bounds each element of the vector $\vH^i \tilde{\vM}_k^i \bzeta^i$ i.e., $(\bmu_{d,k}^i)_{\bar{m}} \geq \max_{\bzeta^i \in \calU_i} (\vH^i \tilde{\vM}_k^i \bzeta^i)_{\bar{m}}$. Based on proposition \ref{prop: Max bound on mean}, we also have $\min_{\bzeta^i \in \calU_i} (\vH^i \tilde{\vM}_k^i \bzeta^i)_{\bar{m}} \geq -(\bmu_{d,k}^i)_{\bar{m}}$. Consequently, we consider a tighter approximation to the constraint \eqref{nonconvex obstacle K1} as $\| \bmu_{d,k}^i \|_2 \leq \tilde{c}_i^k$.
%\qed{}
\endproof
\begin{figure}[t!] 
    \centering
       \includegraphics[width= 0.75\linewidth, trim={0cm 0cm 0cm 0cm},clip]{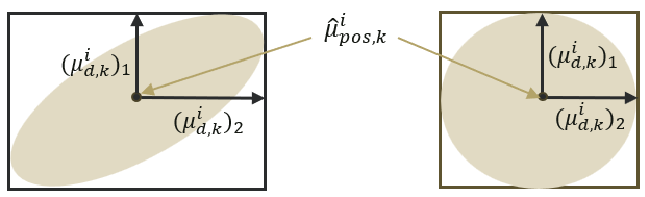}
        \hfill
  \caption{ \textbf{Characterization of State for 2D case:} The state components corresponding to position of agent $i$ at timestep $k$.  }
  \label{fig1} 
\end{figure}
\begin{remark}
    The variable $\bmu_{d,k}^i$ serves as a bound for the deviation in the state $x_k^i$ due to the deterministic disturbance, and is used to characterize the agents' state (as shown in Fig.\ref{fig1} for a 2D case).
\end{remark}
In this formulation, only the constraints \eqref{mu_d expression}, independent of $\bm{p}_{\text{obs}}$, include $\vS_i$. Hence, only the SOCP constraints \eqref{semi-definite obstacle new} of size $n_{\text{pos}}+1$ ($\ll \bar{n} +1$) increase with the number of obstacles significantly reducing the complexity compared to the SDP constraints \eqref{obstacle_semidef_old}. A detailed complexity analysis is provided in Section \ref{sec: Computational Complexity Analysis}.

% Therefore, it is only the amount of constraints of the form \eqref{semi-definite obstacle new} that would increase with the number of obstacles. Also, the size of each SOCP constraint \eqref{semi-definite obstacle new} is $n_{\text{pos}}+1$ which is significantly less than the size of the SDP constraints \eqref{obstacle_semidef_old} (i.e., $\bar{n} +1$). 
% Additionally, SOCP constraints are computationally less expensive than SDP constraints \cite{LOBO1998193}.
%
\subsubsection{Nonconvex Inter-agent Norm Constraints}
Subsequently, we proceed with considering the inter-agent robust collision avoidance constraints \eqref{robust inter-agent collision avoidance},
which can be rewritten as 
\begin{align}
    \| \tilde{\vH}_{ij}\Tilde{\bx}^{ij}_k \|_2 \geq c_{ij}, \quad \forall \bzeta^i \in \calU_i, ~ \bzeta^j \in \calU_j, ~ j \in \calN_i,
    \label{nonconvex interagent concatenated constraint}
\end{align}
where $\tilde{\vH}_{ij} = \begin{bmatrix}
    \vH^i,& -\vH^j
\end{bmatrix}$, $\Tilde{\bx}^{ij}_k = \begin{bmatrix}
    \bx^i_k; & \bx^j_k
\end{bmatrix}$. The uncertainty involved for the variable $\Tilde{\bx}^{ij}$ would be $\Tilde{\bzeta}^{ij} = \begin{bmatrix}
    \bzeta^i; & \bzeta^j
\end{bmatrix}$. Further, it should be noted that if $\bzeta^i, \bzeta^j$ belong to the defined ellipsoidal sets $\calU_i$, $\calU_j$ respectively, then the uncertainty vector $\Tilde{\bzeta}^{ij}$ can also be modeled by an ellitope (which is an intersection of ellipsoids). Thus, by using the initial approach, the constraint \eqref{nonconvex interagent concatenated constraint} would result in semi-definite constraints that are of the size $\bar{n}_i + \bar{n}_j + 1$ \mbox{
\cite{kotsalis2020convex}}\hskip0pt. This might be particularly disadvantageous in a large-scale multi-agent setup, which might involve a significant number of neighbor agents. Thus, we derive equivalent/approximate tractable constraints using the second approach. 

% For that, let us first rewrite a single constraint corresponding to a neighbor agent $j$ of the agent $i$ ($j \in \calN_i$) from the set of constraints \eqref{robust inter-agent collision avoidance} as follows
For that, let us first rewrite the constraint \eqref{robust inter-agent collision avoidance} as follows
\begin{equation}
\begin{aligned}
    & \| \hat{\bmu}^i_{\text{pos},k} + \vM_{\text{det}}^{k,i} (\bzeta^i)
    - (\hat{\bmu}^j_{\text{pos},k}  +\vM_{\text{det}}^{k,j} (\bzeta^j) ) \|_2 
    \geq c_{ij}, \nonumber
    \\
    &~~~~~~~~~~~~~~~~~~~~~~~~~~~~~~~~~~~~~~~~~~~~~~~
    \forall \bzeta^i \in \calU_i, ~ \bzeta^j \in \calU_j.
\end{aligned}
\end{equation}
We then consider the following tighter approximation to the above 
constraint-
\begin{equation}
\begin{aligned}
     & \| \hat{\bmu}^i_{\text{pos},k}  
    - \hat{\bmu}^j_{\text{pos},k}  \|_2
    - \| \vM_{\text{det}}^{k,i} (\bzeta^i)
    -  \vM_{\text{det}}^{k,j} (\bzeta^j) \|_2  
    \geq c_{ij}, \nonumber
    \\
    &~~~~~~~~~~~~~~~~~~~~~~~~~~~~~~~~~~~~~~~~~~~~
    \forall \bzeta^i \in \calU_i, ~ \bzeta^j \in \calU_j
\end{aligned}
\end{equation}
% \begin{displaymath}
% \begin{aligned}
%     \forall \bzeta^i & \in \calU_i, ~ \bzeta^j \in \calU_j,   \\
%     &
%     \begin{aligned}
%     &
%     \| \hat{\bmu}^i_{\text{pos},k}  
%     - \hat{\bmu}^j_{\text{pos},k}  \|_2
%     % \\
%     % & \qquad \qquad \qquad
%     - \| \vH^i \tilde{\vM}_k^i \bzeta^i
%     -  \vH^j \tilde{\vM}_k^j \bzeta^j \|_2  
%     \geq c_{ij}
%     \end{aligned}
% \end{aligned}
% \end{displaymath}
%
As with the previous case, we introduce a slack variable $\tilde{c}_{ij}$, through which we can rewrite the above constraint as follows
\begin{align}
    & 
    \| \hat{\bmu}^i_{\text{pos},k}  
    - \hat{\bmu}^j_{\text{pos},k} \|_2
    \geq \tilde{c}_{ij}^k + c_{ij}
    \label{inter-agent collision avoidance ubar}\\
    &
    \begin{aligned}
         \max_{\bzeta^i \in \calU_i, ~ \bzeta^j \in \calU_j}
         \| \vM_{\text{det}}^{k,i} (\bzeta^i)
    -  \vM_{\text{det}}^{k,j} (\bzeta^j) \|_2  
        \leq \tilde{c}_{ij}^k 
    \end{aligned}    
    \label{interagent collision avoidance K1}
\end{align}
%
% The non-convex constraint \eqref{inter-agent collision avoidance ubar} is addressed by linearizing it around 
% the nominal values $\bmu_{x,\bar{u}}^{i,l}, \bmu_{x,\bar{u}}^{j,l}$. 
% Similar to the previous case, the constraint \eqref{interagent collision avoidance K1} is semi-infinite, and deriving its equivalent tractable constraints is NP-hard. Thus,
% We derive the tighter approximations of the constraint \eqref{interagent collision avoidance K1}.
Similar to the previous case, derivation of equivalent tractable version of the above is NP-hard, hence we derive the tighter approximation of the constraint \eqref{interagent collision avoidance K1}.
\begin{proposition} \label{Proposition: nonconvex interagent deterministic}
    The tighter approximate constraints to the semi-infinite constraint \eqref{interagent collision avoidance K1} are given by the  following  second-order cone constraints (SOCP)
    \begin{align}
        & 
       \| \bmu_{d,k}^i +  \bmu_{d,k}^j \|_2
        \leq 
        \tilde{c}_{ij}^k,
        \label{semi-definite interagent new} \\
        &
        (\bmu_{d,k}^i)_{\bar{m}} \geq  
        \sqrt{\tau^i} || \boldsymbol{\Gamma}_i \T \vM_i \T \bm{h}_{k,\bar{m}}^i||_{\vS_i^{-1}},
        \label{interagent mu_d expression} \\
        &
        (\bmu_{d,k}^j)_{\bar{m}} \geq  
        \sqrt{\tau^j} || \boldsymbol{\Gamma}_j \T \vM_j \T \bm{h}_{k,\bar{m}}^j ||_{\vS_j^{-1}},
        \label{interagent mu_d expression 2}
    \end{align}
where $\bm{h}_{k,\bar{m}}^i {}\T, \bm{h}_{k,\bar{m}}^j {}\T $ are the $\bar{m}^{th}$ rows of the matrices $\vH^i \vP_k^i, \vH^j \vP_k^j$ respectively.
\end{proposition}
\proof
    % The proof of this proposition follows a similar intuition as the proposition \ref{prop- obs avoidance}. We derive the bounds of the vector $\vH^i \tilde{\vM}_k^i \bzeta^i
    % -  \vH^j \tilde{\vM}_k^j \bzeta^j$ and consider a tighter approximation of the constraint \eqref{interagent collision avoidance K1} based on that. 
    A detailed proof is provided in Appendix\ref{Appendix: sec1 nonconvex interagent}.
%\qed{}
\endproof
\noindent
Note that the constraints \eqref{interagent mu_d expression}, \eqref{interagent mu_d expression 2} are the same as \eqref{mu_d expression}.
% It should be noted that the constraints \eqref{interagent mu_d expression}, \eqref{interagent mu_d expression 2} will be the same as \eqref{mu_d expression}, therefore alleviating the need to introduce them as additional constraints in our framework. 
%
% \begin{remark}
%     The realizations of the state of the agents, under deterministic uncertainty, is a set of points. Thus making the sharing of the state information among the agents challenging. The proposed reformulation not only provides a tractable form of the constraints but also enables the state to be characterized using variables $\bmu_{x,\bar{u}}^i, \bmu_{d}^i$. Thereby, allowing us to develop a distributed framework where the agents do not need to share the control parameters or store any system information about the other agents. 
% \end{remark}
%
%
\subsubsection{Robust convex norm constraints} \label{Robust convex norm constraints}
% We now look into the convex norm constraints of the form \eqref{Robust convex norm-of-mean constraints}. 
Let us rewrite the convex norm constraint \eqref{Robust convex norm-of-mean constraints} as follows,
\begin{align}
    \|  \hat{\bmu}^i_{\text{pos},k} - \bm{p}_{\text{tar}} \|_2 \leq c_i, \quad \forall \bzeta^i \in \calU_i, \nonumber
\end{align}
% Using \eqref{dynamics - deterministic case}, 
which can written in terms of $ \bmu_{x,\bar{u}}^i$, and $\vM_{\text{det}}^{k,i} (\bzeta^i)$ as 
\begin{align}
     \|  \hat{\bmu}^i_{\text{pos},k}  +  \vM_{\text{det}}^{k,i} (\bzeta^i) - \bm{p}_{\text{tar}} \|_2^2 
    \leq c_i^2, \quad \forall \bzeta^i \in \calU_i, \nonumber
\end{align}
and can be equivalently given as
\begin{align}
    \max_{\bzeta^i \in \calU_i} \| \vM_{\text{det}}^{k,i} (\bzeta^i) 
    + \hat{\bmu}^i_{\text{pos},k} - \bm{p}_{\text{tar}} \|_2^2 
    \leq c_i^2.
    \label{convex norm-of-mean constraints: interm1}
\end{align}
%
%
% We now present the equivalent tractable constraints to the above constraint in the following proposition.
%
\begin{proposition} \label{Proposition: robust convex norm-mean equivalent constraints}
    The equivalent constraints to the constraint \eqref{convex norm-of-mean constraints: interm1} can be given using the S-lemma as follows
    \begin{align}
    & ~~~~~~~~~~~~~~~~~~~~~~ \alpha \in \Rb, \; \alpha \geq 0
    \label{Convex norm-of-mean constraints: final 1}
    \\[0.2cm]
    & \begin{bmatrix}
    \alpha \vS_i & 0 
    & (\vH^i \tilde{\vM}_{\tilde{k}}^i \boldsymbol{\Gamma}_i) \T \\
    0 &  c_i^2 - \alpha \tau^i 
    & (\hat{\bmu}^i_{\text{pos},k} - \bm{p}_{\text{tar}} ) \T \\
    \vH^i \tilde{\vM}_{\tilde{k}}^i \boldsymbol{\Gamma}_i
    & \hat{\bmu}^i_{\text{pos},k} - \bm{p}_{\text{tar}} 
    & \vI_{n_{\text{pos}}}
    \end{bmatrix}
    \succeq 0
    \label{Convex norm-of-mean constraints: final 2}
    \end{align}
\end{proposition}
\proof
    A detailed proof is provided in the Appendix\ref{Appendix: sec1 convex norm s-lemma}. 
    % An overview of the proof involves the following steps. 
    % To outline the proof, we first rewrite the LHS of the constraint \eqref{convex norm-of-mean constraints: interm1} as a quadratic function of $\bm{z}_i$. Then, applying the S-lemma followed by Schur Complement, we obtain the equivalent constraints \eqref{Convex norm-of-mean constraints: final 1} and \eqref{Convex norm-of-mean constraints: final 2}. 
    % Furthermore, it should be noted that $c_i$ in this constraint is a constant, therefore the power of 2 on $c_i$ does not affect the convexity of the constraint \eqref{Convex norm-of-mean constraints: final 2}.
    %\qed
\endproof
As discussed in the previous subsection, solving a problem with numerous SDP constraints is challenging. Thus, the above reformulation can be used only while dealing problems with relatively fewer constraints of the form \eqref{Robust convex norm-of-mean constraints}. Accordingly, we also provide a tighter approximation of the constraints \eqref{Robust convex norm-of-mean constraints}, which are computationally less expensive to solve compared to the equivalent tractable constraints \eqref{Convex norm-of-mean constraints: final 1} and \eqref{Convex norm-of-mean constraints: final 2}. For that, we use an intuition similar to the previous robust constraint cases. We start by rewriting the constraint \eqref{Robust convex norm-of-mean constraints} as follows
\begin{align}
    \| \hat{\bmu}^i_{\text{pos},k} +  \vM_{\text{det}}^{k,i} (\bzeta^i)  
    - \bm{p}_{\text{tar}} \|_2 \leq c_i, \quad \forall \bzeta^i \in \calU_i, 
    \nonumber
\end{align}
Using the triangle inequality, we consider the tighter approximation of the above constraint as follows 
\begin{align}
    \| \hat{\bmu}^i_{\text{pos},k} - \bm{p}_{\text{tar}} \|_2 
    +  \| \vM_{\text{det}}^{k,i} (\bzeta^i)  \|_2 \leq c_i, \quad \forall \bzeta^i \in \calU_i.
    \nonumber
\end{align}
By introducing a slack variable $\tilde{c}_i \in \Rb$, the above constraint can be equivalently given as follows
\begin{align}
    &
    \| \hat{\bmu}^i_{\text{pos},k} - \bm{p}_{\text{tar}} \|_2 
    \leq c_i - \tilde{c}_i^k
    \label{Convex norm-of-mean: tighter final 1}
    \\
    &
    \| \vM_{\text{det}}^{k,i} (\bzeta^i)  \|_2 \leq \tilde{c}_i^k, \quad \forall \bzeta^i \in \calU_i
    \label{Convex norm-of-mean: tighter interm 1}
\end{align}
%
% It can be observed that the constraint \eqref{Convex norm-of-mean: tighter final 1} is a convex constraint and does not need to be reformulated. 
% In addition, 
The constraint \eqref{Convex norm-of-mean: tighter interm 1} is the same as \eqref{nonconvex obstacle K1}, and its equivalent constraints can be given by \eqref{semi-definite obstacle new} and \eqref{mu_d expression}.

\section{Extension to Mixed Disturbance case}
\label{sec: Extension to Mixed Disturbance case}
In this section, we consider mixed disturbance case with both stochastic and deterministic disturbances. We start by rewriting the state dynamics \eqref{dynamics - deterministic case} with $\bw^i$ as follows
% For that, we rewrite the state equation \eqref{dynamics - deterministic case} by including the term corresponding to $\bw^i$ as follows
\begin{align}
    & \bx^i = \bmu_{x,\bar{u}}^i + \vM_i \boldsymbol{\zeta}^i + \vM^{w}_i \bw^i,
    \label{dynamics - mixed case}
\end{align}
where $ \vM^{w}_i =  (\vG_u^i \vK^i + \vI) \vG_{\bw}^i$.
% \begin{equation}
%     \vM^{w}_i =  (\vG_u^i \vK^i + \vI) \vG_{\bw}^i,
% \end{equation}
% and $ \bmu_{x,\bar{u}}^i$, $\vM_i$ are as defined in the section \ref{sec: Reformulation of Robust Deterministic Constraints}. 
For simplicity, we also define $\tilde{\vM}^{w,i}_k = \vP_k^i \vM^{w}_i$, and $\vM_{\text{mix}}^{k,i}(\bzeta^i, \bw^i) = \vH^i (\tilde{\vM}_k^i \bzeta^i + \tilde{\vM}^{w,i}_k \bw^i)$. Further, the expectation of the state $\Eb[\bx^i]$ is given as follows 
\begin{align}
    \Eb[\bx^i] = \bmu_{x,\bar{u}}^i + \vM_i \boldsymbol{\zeta}^i,
    \label{Mixed case: Expecation of x}
\end{align}
and covariance of the state $ \Cov(\bx^i)$ is given as 
\begin{equation}
    \Cov(\bx^i) = \vM^{w}_i \bSigma_{\bw^i} \vM^{w}_i{}\T,
\end{equation}
and matrix $\boldsymbol{\varphi}^i \in \Sb^{+}_{n_{x_i} + T n_{w_i}}$ is defined such that $\boldsymbol{\varphi}^i \boldsymbol{\varphi}^i {}\T = \bSigma_{\bw^i}$.
% and is dependent only on the deterministic disturbance. 
% This section is structured as follows: we first derive the tractable reformulations for each of the constraints (a)-(c) defined in the section \ref{sec: robust constraints in mixed case}. Then, we derive the tractable reformulations of the rest of the constraints (i.e., chance constraints) by characterizing the stochastic uncertainty as a deterministic disturbance.
%
%
Subsequently, we reformulate each robust constraint presented in Section \ref{sec: robust constraints in mixed case}.
\subsubsection{Robust Linear Chance Constraints}
Let us rewrite the chance constraint \eqref{Mixed Constraints: Robust linear Chance Constraints} as follows
\begin{align}
    \Pb \big( \ba_i \T \bx^i  > b_i \big) \leq p_i, 
\quad \forall \bzeta ^i \in \calU_i,
\nonumber
\end{align}
which can be further rewritten as
\begin{equation}
    \Pb \big( \ba_i \T \bmu_{x,\bar{u}}^i 
    + \max_{\bzeta ^i \in \calU_i} \ba_i \T \vM_i \boldsymbol{\zeta}^i 
    + \ba_i \T \vM^{w}_i \bw^i  > b_i \big) \leq p_i.
\end{equation}
Using Proposition \ref{prop: Max bound on mean}, we can rewrite the above constraint as
\begin{equation}
    \Pb \big( \ba_i \T \bmu_{x,\bar{u}}^i 
    + \sqrt{\tau^i} || \boldsymbol{\Gamma}_i \T \vM_i \T \ba_i ||_{\vS_i^{-1}} 
    + \ba_i \T \vM^{w}_i \bw^i  > b_i \big) \leq p_i \nonumber
\end{equation}
The above constraint can be equivalently given by the following set of constraints \cite{kotsalis2020convex} 
\begin{align}
    & 
    \ba_i \T \bmu_{x,\bar{u}}^i
    + 
    \sqrt{\tau^i} || \boldsymbol{\Gamma}_i \T \vM_i \T \ba_i ||_{\vS_i^{-1}} \leq b_i - \alpha, 
    \label{chance constraints 1} \\
    &
    \eta || \ba \T \vM^{w}_i \boldsymbol{\varphi}^i ||_2 \leq \alpha 
    \label{chance constraints 2},
\end{align}
where $ \eta : \Pb_{y \in \calN(0,1)}(y \geq \eta) = p_i $. 
%
%
%
% \subsection{Robust Nonlinear Chance Constraints}
% In this subsection, we derive the equivalent/approximate tractable versions of the non-linear chance constraints defined by \eqref{Mixed Constraints: nonconvex norm chance} - \eqref{Mixed Constraints: convex norm chance} in the section \ref{sec: robust constraints in mixed case}. 
% For that, we start by characterizing the stochastic disturbance as a deterministic uncertainty and convert the constraints involving two types of deterministic uncertainty. Then, we derive the reformulated constraints using the results obtained in the deterministic case (Section \ref{sec: Reformulation of Robust Deterministic Constraints}).
%
\subsubsection{Robust nonconvex norm chance constraints}
\label{subsec: Robust nonconvex norm chance constraints}
Let us rewrite the constraint \eqref{Mixed Constraints: nonconvex norm chance} here as follows
\begin{align}
    \Pb \big( \| \vH^i x_k^i - \bm{p}_{\text{obs}} \|_2 \geq c_i \big) \geq 1-p_i, \quad \forall \bzeta^i \in \calU_i.
    \label{nonconvex chance constraint eq1}
\end{align}
Since there are two disturbances involved in the above non-linear chance constraint, it is hard to derive its exact tractable version. Thus, we derive a tighter approximation extending the approaches from deterministic case. We begin by writing the following relation based on triangle inequality -
% \begin{equation}
% \begin{aligned}
%     \| \vH^i x_k^i - \bm{p}_{\text{obs}} \|_2 
%     \geq
%      & \| \hat{\bmu}^i_{\text{pos},k} - \mathbf{p}_{\text{obs}} \|_2
%      \\
%      &~~~~
%      - \| \vH^i \big( \tilde{\vM}_k^i \bzeta^i
%      + \tilde{\vM}^{w,i}_k \bw^i \big) \|_2.
%      \nonumber
% \end{aligned}
% \end{equation}
\begin{equation}
\begin{aligned}
    \| \vH^i \bx_k^i - \bm{p}_{\text{obs}} \|_2 
    \geq
     & \| \hat{\bmu}^i_{\text{pos},k} - \mathbf{p}_{\text{obs}} \|_2
     % \\
     % &~~~~
     - \| \vM_{\text{mix}}^{k,i}(\bzeta^i, \bw^i) \|_2.
     \nonumber
\end{aligned}
\end{equation}
Based on the above relation, consider the following tighter approximation of the constraint \eqref{nonconvex chance constraint eq1}
\begin{align}
    & \| \hat{\bmu}^i_{\text{pos},k} - \bm{p}_{\text{obs}} \|_2 
    \geq  \tilde{c}_i^k + c_i,
    \label{mixed: nonconvex obstacle final ubar}
    \\
    & \Pb \big( \| \vM_{\text{mix}}^{k,i}(\bzeta^i, \bw^i)  \|_2 
    \leq \tilde{c}_i^k \big) \geq 1-p_i,
    ~\forall \bzeta^i \in \calU_i.
    \label{nonconvex mixed chance obstacle K1}
\end{align}
% \begin{equation}
%     \| \hat{\bmu}^i_{\text{pos},k} - \bm{p}_{\text{obs}} \|_2 
%     \geq  \tilde{c}_i^k + c_i,
%     \label{mixed: nonconvex obstacle final ubar}
% \end{equation}
% %
% \begin{equation}
% \begin{aligned}
%     & 
%     \Pb \big( \| \vH^i ( \tilde{\vM}_k^i \bzeta^i
%      + \tilde{\vM}^{w,i}_k \bw^i ) \|_2 
%     \leq \tilde{c}_i^k \big) \geq p_i,
%     % \\
%     % &~~~~~~~~~~~~~~~~~~~~~~~~~~~~~~~~~~~~~~~~~~~~~~~~~~
%     ~\forall \bzeta^i \in \calU_i.
%     \label{nonconvex mixed chance obstacle K1}
% \end{aligned}
% \end{equation}
% \begin{align}
%     & \| \hat{\bmu}^i_{\text{pos},k} - \bm{p}_{\text{obs}} \|_2 
%     \geq  \tilde{c}_i^k + c_i,
%     \label{mixed: nonconvex obstacle final ubar}\\
%     & 
%     \Pb ( \| \vH^i \big( \tilde{\vM}_k^i \bzeta^i
%      + \tilde{\vM}^{w,i}_k \bw^i \big) \|_2 
%     \leq \tilde{c}_i^k ) \geq 1-p_i, 
%      \quad \forall \bzeta^i \in \calU_i.
%    \label{nonconvex mixed chance obstacle K1}
% \end{align}
%
%
The nonconvex constraint \eqref{mixed: nonconvex obstacle final ubar}, same as constraint \eqref{nonconvex obstacle final ubar} in the deterministic case, is linearized around a nominal trajectory. 

We now derive the tractable version of the constraint \eqref{nonconvex mixed chance obstacle K1}, dependent on both the disturbances $\bzeta^i$ and $\bw^i$, in the following two steps. 
% We derive the tractable constraints in the following two steps. 
First, we characterize the stochastic disturbance $\bw^i$ as deterministic uncertainty \cite{CINQUEMANI20112082}, \cite{VANHESSEM2006225}.
% Next, we reformulate the constraints using the similar intuition as in the deterministic case (Section \ref{sec: Reformulation of Robust Deterministic Constraints}).
%
% \paragraph{Characterizing stochastic disturbance as deterministic uncertainty}
% In this section, we use a similar method as in \cite{CINQUEMANI20112082}, \cite{VANHESSEM2006225}.
For that, we introduce a variable $\bar{\bw}^i \in \calN(0, \bSigma_{\bar{\bw}^i}) $ such that $\bar{\bw}^i = \vH^i \tilde{\vM}^{w,i}_k \bw^i$,
% \begin{equation}
%     \bar{\bw}^i = \vH^i \tilde{\vM}^{w,i}_k \bw^i, \nonumber
% \end{equation}
%
and $\bSigma_{\bar{\bw}^i} = \vH^i \tilde{\vM}^{w,i}_k \bSigma_{\boldsymbol{\bw}^i}  \tilde{\vM}^{w,i}_k{}^T \vH^i{}^T $. Subsequently, we also define $\hat{\bw}^i \in \calN(0, \vI_{n_{\text{pos}}})$ such that $\bar{\bw}^i$ can written as $\bar{\bw}^i = \bSigma_{\bar{\bw}^i}^{1/2} \hat{\bw}^i$.
% \begin{equation}
%     \bar{\bw}^i = \bSigma_{\bar{\bw}^i}^{1/2} \hat{\bw}^i. \nonumber
% \end{equation}
We can then rewrite \eqref{nonconvex mixed chance obstacle K1} in terms of $\hat{\bw}^i$ as follows
\begin{equation}
\begin{aligned}
    % \\
    % & 
    \Pb \big( \| \vM_{\text{det}}^{k,i} (\bzeta^i)
     +  \bSigma_{\bar{\bw}^i}^{1/2} \hat{\bw}^i \|_2 
    \leq \tilde{c}_i^k \big) \geq 1-p_i, 
    ~~
    \forall \bzeta^i \in \calU_i
    \label{mixed: nonconvex norm chance wbar}
\end{aligned}
\end{equation}
% \begin{equation}
%     \bw^i = \bSigma_{\boldsymbol{\bw}^i}^{1/2} \hat{\bw}^i.
% \end{equation}
% Now, we rewrite the constraint \eqref{nonconvex mixed chance obstacle K1} in terms of $\hat{\bw}^i$ as follows.
% \begin{equation}
% \begin{aligned}
%     & \forall \bzeta^i  \in \calU_i, \\
%     &~ \Pb ( \| \vH^i \big(\tilde{\vM}_k^i \bzeta^i
%      + \tilde{\vM}^{w,i}_k \bSigma_{\boldsymbol{\bw}^i}^{1/2} \hat{\bw}^i \big) \|_2 
%     \leq \tilde{c}_i^k ) \geq 1-p_i
% \end{aligned}
% \label{mixed: nonconvex norm chance wbar}
% \end{equation}
%
% We now present the tighter approximation of the above constraint in the following proposition.
% To derive a tighter approximation of the above constraint,
We will use a similar intuition as in \cite{VANHESSEM2006225}, \cite{CINQUEMANI20112082} to assert the following - if the constraint $\| \vM_{\text{det}}^{k,i} (\bzeta^i)
     +  \bSigma_{\bar{\bw}^i}^{1/2} \hat{\bw}^i \|_2 
    \leq \tilde{c}_i^k$ is satisfied for all stochastic disturbances $\hat{\bw}^i$ lying inside its $100(1-p_i)\%$ confidence ellipsoidal, then it implies that the constraint \eqref{mixed: nonconvex norm chance wbar} is satisfied. Thus, we propose the following approximation of the constraint \eqref{mixed: nonconvex norm chance wbar} - 
\begin{equation}
\begin{aligned}
    % \\
    % & \qquad \quad
    \| \vM_{\text{det}}^{k,i} (\bzeta^i)
     +  \bSigma_{\bar{\bw}^i}^{1/2} \hat{\bw}^i \|_2 
    \leq \tilde{c}_i^k, 
    ~
    \forall \bzeta^i  \in \calU_i, 
     \hat{\bw}^i \in \calW_i[\bar{\eta}_i ],
\end{aligned}
\label{mixed: nonconvex norm chance wbar interm2 - appendix}
\end{equation}
where $\calW_i[\bar{\eta}_i]$ is $100(1-p_i)\%$ confidence ellipsoid of $\hat{\bw}^i$, and is given as follows,
\begin{equation}
\begin{aligned}
    \calW_i[\bar{\eta}_i] 
    = \{ \hat{\bw} \in \Rb^{n_{\text{pos}}} \; | \; \| \hat{\bw} \|_2^2 \leq \bar{\eta}_i, 
    ~ \bar{\eta}_i = \mathcal{Q} (1-p_i) \}, 
    \label{Stochastic noise uncertainty set - confidence ellipsoid}
\end{aligned}
\end{equation}
and $\mathcal{Q}(.)$ is the quantile function of chi-square distribution with $n_{\text{pos}}$ degrees of freedom. Subsequently, we derive the tighter approximation of the constraint \eqref{mixed: nonconvex norm chance wbar interm2 - appendix}.
\begin{proposition}
\label{proposition: nonconvex chance}
The tighter approximation of the semi-infinite constraints \eqref{mixed: nonconvex norm chance wbar interm2 - appendix} is as follows-
\begin{align}
    & 
    \| \bmu_{w,k}^i \|_2 
    \leq \tilde{c}_i^k 
    \label{mixed: nonconvex chance mu_d upper bound final} \\
    &
    (\bmu_{w,k}^i)_{\bar{m}} \geq  
    \sqrt{\tau^i}\| \boldsymbol{\Gamma}_i \T \vM_i \T h_{k,\bar{m}}^i\|_{\vS_i^{-1}} 
    \nonumber \\
    & \qquad \qquad \qquad
    +
    \sqrt{\bar{\eta}_i} \| \boldsymbol{\varphi}^i \vM_w^i {}\T h_{k,\bar{m}}^i \|_2
    \label{mixed: nonconvex chance mu_d lower bound final}
\end{align}
where $h_{k,\bar{m}}^i {}\T$ is $\bar{m}^{th}$ row of the matrix $\vH^i \vP_k^i$.
    % \begin{align}
    %     & h_{k,\bar{m}}^i {}\T \text{ is } \bar{m}^{th} \text{ row of the matrix } \vH^i \vP_k^i, \\
    %     &
    %     \bar{\eta}_i = \mathcal{Q} (1-p_i), \\
    %     & \mathcal{Q}(.) -  \text{Quantile function of Chi-square distribution with } \nonumber\\
    %     & \qquad \qquad n_{x_i} + T n_{w_i} \text{ degrees of freedom.}
    %     \nonumber
    % \end{align}
\end{proposition}
\proof
    A detailed proof is provided in the Appendix\ref{Appendix: sec2 proposition 5}. 
    % To outline the proof, we use the idea that the satisfaction of the constraint $\| \vH^i \tilde{\vM}_k^i \bzeta^i
    %  +  \bSigma_{\bar{\bw}^i}^{1/2} \hat{\bw}^i \|_2 
    % \leq \tilde{c}_i^k $ for all stochastic disturbances lying inside the $100(1-p_i)\%$ confidence ellipsoid of $\hat{\bw}^i$ implies the satisfaction of constraint \eqref{mixed: nonconvex norm chance wbar} \cite{CINQUEMANI20112082}, \cite{VANHESSEM2006225}. Thus, we consider the tighter approximation of the constraint \eqref{mixed: nonconvex norm chance wbar} by characterizing the stochastic disturbance by a deterministic uncertainty lying inside the $100(1-p_i)\%$ confidence ellipsoid of $\hat{\bw}^i$ (an uncertainty set of a similar form as \eqref{ellipsoid uncertainty set}). We then use the results derived in the deterministic case (Section \ref{sec: Robust Constraints for deterministic case}) to derive the approximate tractable constraints \eqref{mixed: nonconvex chance mu_d upper bound final}, \eqref{mixed: nonconvex chance mu_d lower bound final}.
\endproof
\subsubsection{Robust inter-agent nonconvex chance constraints}
Let us now rewrite the constraint \eqref{Mixed Constraints: inter-agent nonconvex chance} as
\begin{equation}
\Pb( \| \vH^i \bx^i_k - \vH^j \bx^j_k \|_2 \geq c_{ij}) \geq 1-p_{ij}, ~~ \forall \bzeta^i \in \calU_i, ~ \bzeta^j \in \calU_j
\label{mixed: nonconvex interagent chance eq1}
\end{equation}
Similar to the earlier non-convex constraints, we consider the following tighter approximation of the constraint \eqref{mixed: nonconvex interagent chance eq1} based on the triangle inequality,
\begin{align}
    & \| \hat{\bmu}^i_{\text{pos},k} 
    - \hat{\bmu}^j_{\text{pos},k} \|_2
    \geq \tilde{c}_{ij} + c_{ij}^k
    \label{mixed: inter-agent collision avoidance ubar}
    \\
    &
    \Pb \bigg( \| \vM_{\text{mix}}^{k,i}(\bzeta^i, \bw^i)
    -  \vM_{\text{mix}}^{k,j}(\bzeta^j, \bw^j) \|_2  
        \leq \tilde{c}_{ij}^k \bigg) \geq 1- p_{ij} \nonumber
    \\
    & \qquad \qquad \qquad \qquad \qquad \qquad 
    \forall ~ \bzeta^i \in \calU_i, ~ \bzeta^j \in \calU_j
    \label{mixed: interagent collision avoidance K1}
\end{align}
% \begin{equation}
%      \| \hat{\bmu}^i_{\text{pos},k} 
%     - \hat{\bmu}^j_{\text{pos},k} \|_2
%     \geq \tilde{c}_{ij} + c_{ij}^k
%     \label{mixed: inter-agent collision avoidance ubar}
% \end{equation}
% \begin{equation}
% \begin{aligned}
%     &
%     \Pb \bigg( \| \vM_{\text{mix}}^{k,i}(\bzeta^i, \bw^i)
%     -  \vM_{\text{mix}}^{k,j}(\bzeta^j, \bw^j) \|_2  
%         \leq \tilde{c}_{ij}^k \bigg) \geq 1- p_{ij}
%     \\
%     & \forall ~ \bzeta^i \in \calU_i, ~ \bzeta^j \in \calU_j
% \end{aligned}
% \label{mixed: interagent collision avoidance K1}
% \end{equation}
% \begin{align}
%     & 
%     \| \hat{\bmu}^i_{\text{pos},k} 
%     - \hat{\bmu}^j_{\text{pos},k} \|_2
%     \geq \tilde{c}_{ij} + c_{ij}^k
%     \label{mixed: inter-agent collision avoidance ubar}\\
%     &
%     \begin{aligned}
%     & \forall ~ \bzeta^i \in \calU_i, ~ \bzeta^j \in \calU_j, \\
%     & ~
%     {\color{red}
%     \| \vH^i \big( \tilde{\vM}_k^i \bzeta^i
%      + \tilde{\vM}^{w,i}_k \bw^i \big)
%     -  \vH^j \big( \tilde{\vM}_k^j \bzeta^j
%      + \tilde{\vM}^{w,j}_k \bw^j \big)\|_2  
%         \leq \tilde{c}_{ij}^k }
%     \end{aligned}    
%     \label{mixed: interagent collision avoidance K1}
% \end{align}
%
The nonconvex constraint \eqref{mixed: nonconvex obstacle final ubar} is of the same form as \eqref{inter-agent collision avoidance ubar} and is addressed by linearizing around a nominal trajectory. The constraint \eqref{mixed: interagent collision avoidance K1} is a semi-infinite constraint and we use a similar intuition as in the subsection \ref{subsec: Robust nonconvex norm chance constraints} to derive its approximate tractable version. The disturbances $\bar{\bw}^i$, $\bar{\bw}^j$, $\hat{\bw}^i$, and $\hat{\bw}^j$ are considered as defined in the subsection \ref{subsec: Robust nonconvex norm chance constraints}. We will now rewrite the constraint \eqref{mixed: interagent collision avoidance K1} in terms of $\hat{\bw}^i$, and $\hat{\bw}^j$ as follows
\begin{equation}
\begin{aligned}
    & \forall ~ \bzeta^i \in \calU_i, ~ \bzeta^j \in \calU_j, \\
    & \quad 
    \Pb \big( \| \vM_{\text{det}}^{k,i} (\bzeta^i)
    - \vM_{\text{det}}^{k,j} (\bzeta^j)
    \\
    & \qquad \quad
    + \bSigma_{\bar{\bw}^i}^{1/2} \hat{\bw}^i 
    - \bSigma_{\bar{\bw}^j}^{1/2} \hat{\bw}^j \|_2  
    \leq \tilde{c}_{ij}^k \big) \geq 1-p_{ij}
\end{aligned}
\label{mixed: interagent collision avoidance eq2}
\end{equation}
%
%
% For that, we again start by introducing a variable  $\bar{\bw}^{i,j} \in \calN(0, \bSigma_{\bar{\bw}^{i,j}}) $  such that
% \begin{equation}
%     \bar{\bw}^{i,j} = \vH^i \tilde{\vM}^{w,i}_k \bw^i
%     - \vH^j \tilde{\vM}^{w,j}_k \bw^j, \nonumber
% \end{equation}
% with $\Eb[\bar{\bw}^{i,j}] = 0$, and $\bSigma_{\bar{\bw}^{i,j}} = \vH^i \tilde{\vM}^{w,i}_k \bSigma_{\boldsymbol{\bw}^i}  \tilde{\vM}^{w,i}_k{}^T \vH^i{}^T 
% + \vH^j \tilde{\vM}^{w,j}_k \bSigma_{\boldsymbol{\bw}^j}  \tilde{\vM}^{w,j}_k{}^T \vH^j{}^T$. Subsequently, we define $\hat{\bw}^{i,j} \in \calN(0, \vI_{n_{\text{pos}}})$ such that $\bar{\bw}^{i,j}$ can written as 
% \begin{equation}
%     \bar{\bw}^{i,j} = \bSigma_{\bar{\bw}^{i,j}}^{1/2} \hat{\bw}^{i,j}. \nonumber
% \end{equation}
% We can then rewrite \eqref{mixed: interagent collision avoidance K1} in terms of $\hat{\bw}^{i,j}$ as follows
% \begin{equation}
% \begin{aligned}
%     \forall ~ \bzeta^i & \in \calU_i, ~ \bzeta^j \in \calU_j, \\
%     &
%     \| \vH^i \tilde{\vM}_k^i \bzeta^i
%     -  \vH^j \tilde{\vM}_k^j \bzeta^j
%     + \bSigma_{\bar{\bw}^{i,j}}^{1/2} \hat{\bw}^{i,j} \|_2  
%     \leq \tilde{c}_{ij}^k 
% \end{aligned}
% \label{mixed: interagent collision avoidance eq2}
% \end{equation}
%
% We present the tighter approximation of the above constraints in the following proposition.
\begin{proposition} \label{Proposition: nonconvex interagent chance constraints}
    The tighter approximation of the semi-infinite constraints \eqref{mixed: interagent collision avoidance eq2} is given as
    \begin{align}
        & 
       \| \bar{\bmu}_{w,k}^i +  \bar{\bmu}_{w,k}^j \|_2
        \leq 
        \tilde{c}_{ij}^k,
        \label{mixed: nonconvex interagent mu_d upperbound} \\
        &
        (\bar{\bmu}_{w,k}^i)_{\bar{m}} \geq  
        \sqrt{\tau^i} || \boldsymbol{\Gamma}_i \T \vM_i \T \bm{h}_{k,\bar{m}}^i||_{\vS_i^{-1}}, \nonumber \\
        & \qquad \qquad \qquad
        + \sqrt{\bar{\eta}_{ij}} \|  \boldsymbol{\varphi}^i \vM_w^i {}\T \bm{h}_{k,\bar{m}}^i \|_2
        \label{mixed: nonconvex interagent mu_di lowerbound} \\
        &
        (\bar{\bmu}_{w,k}^j)_{\bar{m}} \geq  
        \sqrt{\tau^j} || \boldsymbol{\Gamma}_j \T \vM_j \T \bm{h}_{k,\bar{m}}^j ||_{\vS_j^{-1}}, \nonumber \\
        & \qquad \qquad \qquad
        + \sqrt{\bar{\eta}_{ij}} \|  \boldsymbol{\varphi}^j \vM_w^j {}\T \bm{h}_{k,\bar{m}}^j \|_2
         \label{mixed: nonconvex interagent mu_dj lowerbound} 
    \end{align}
    where $\bm{h}_{k,\bar{m}}^i {}\T$, $\bm{h}_{k,\bar{m}}^j {}\T$ are $\bar{m}^{th}$ rows of the matrices $\vH^i \vP_k^i$ and $\vH^j \vP_k^j$ respectively; $\bar{\eta}_{ij} = \mathcal{Q} (\sqrt{1-p_{ij}})$.
    % and $\mathcal{Q}(.)$ is the Quantile function of Chi-square distribution with $n_{\text{pos}}$ degrees of freedom.
\end{proposition}
\proof
    A detailed proof is provided in the Appendix\ref{Appendix: Proof of Proposition 6 nonconvex interagent chance}.
\endproof
\subsubsection{Robust convex norm chance constraints}
We can rewrite the constraint \eqref{Mixed Constraints: convex norm chance} as follows
% We now derive the equivalent/approximate tractable constraints to the convex norm chance constraints \eqref{Mixed Constraints: convex norm chance} given as follows
\begin{equation}
\Pb (\| \vH^i \bx_k^i - \bm{p}_{\text{tar}} \|_2 \leq c_i) \geq 1-p_i, \quad \forall \bzeta^i \in \calU_i.
\label{mixed: convex norm chance eq1}
\end{equation}
Even if the above constraint is convex, it is still challenging to derive exact tractable form. 
% Even though the norm constraint involved is a convex constraint in this case, it is still challenging to directly reformulate the above constraint to a tractable form. 
Thus, we use a similar methodology as in Subsection \ref{subsec: Robust nonconvex norm chance constraints} to derive its approximate tractable constraints. For that, we again use the triangle inequality to derive the following relation -
\begin{equation}
\begin{aligned}
    \| \vH^i x_k^i - \bm{p}_{\text{tar}} \|_2 
    \leq
     \| \hat{\bmu}^i_{\text{pos},k} - \mathbf{p}_{\text{tar}} \|_2
     % \\
     % &~~ 
     + \| \vM_{\text{mix}}^{k,i}(\bzeta^i, \bw^i) \|_2 \nonumber
\end{aligned}
\end{equation}
Using this relation, we now consider the tighter approximation of the constraint \eqref{mixed: convex norm chance eq1} as follows
\begin{align}
    & \| \hat{\bmu}^i_{\text{pos},k}  - \bm{p}_{\text{tar}} \|_2 
    \leq  c_i - \tilde{c}_i^k ,
    \label{mixed: convex norm final ubar}
    \\
    & 
    \Pb ( \| \vM_{\text{mix}}^{k,i}(\bzeta^i, \bw^i) \|_2 
    \leq \tilde{c}_i^k ) \geq 1-p_i, 
    ~~ \forall \bzeta^i \in \calU_i
   \label{convex norm mixed chance K1}
\end{align}
% \begin{equation}
% \begin{aligned}
%     \| \hat{\bmu}^i_{\text{pos},k}  - \bm{p}_{\text{tar}} \|_2 
%     \leq  c_i - \tilde{c}_i^k ,
%     \label{mixed: convex norm final ubar}
% \end{aligned}
% \end{equation}
% \begin{equation}
% \begin{aligned}
%     \forall \bzeta^i \in & \calU_i, \\
%     &
%     \Pb ( \| \vM_{\text{mix}}^{k,i}(\bzeta^i, \bw^i) \|_2 
%     \leq \tilde{c}_i^k ) \geq 1-p_i.
%    \label{convex norm mixed chance K1}
% \end{aligned}
% \end{equation}
% \begin{align}
%  & \| \hat{\bmu}^i_{\text{pos},k}  - \bm{p}_{\text{tar}} \|_2 
%     \leq  c_i - \tilde{c}_i^k ,
%     \label{mixed: convex norm final ubar}\\
%     & 
%     \Pb ( \| \vH^i \big( \tilde{\vM}_k^i \bzeta^i
%      + \tilde{\vM}^{w,i}_k \bw^i \big) \|_2 
%     \leq \tilde{c}_i^k ) \geq 1-p_i, 
%      \quad \forall \bzeta^i \in \calU_i.
%    \label{convex norm mixed chance K1}
% \end{align}
The constraint \eqref{mixed: convex norm final ubar} is convex and tractable. The semi-infinite constraint \eqref{convex norm mixed chance K1} has the same form as constraint \eqref{nonconvex mixed chance obstacle K1}, and its tractable approximation is obtained using Proposition \ref{proposition: nonconvex chance}.
% Thus, the approximate tractable constraints can be derived using the proposition \ref{proposition: nonconvex chance}.
%
\subsubsection{Robust Covariance Constraints}
We now consider the covariance constraints of the form \eqref{Mixed Constraints: Covariance}, which is rewritten here as follows
\begin{equation}
    \Cov(\bx^i_{\tilde{k}}) \preceq \bSigma^i_{\tilde{k}}, \quad \forall \bzeta^i \in \calU_i.
    \label{covariance constraint: interm1}
\end{equation}
Using \eqref{dynamics - mixed case}, the covariance 
 $\Cov(\bx^i_{\tilde{k}})$ is given as 
\begin{align}
    \Cov(\bx^i_{\tilde{k}}) 
    = 
    \vP_{\tilde{k}}^i (\vG_u^i \vK^i + \vI) \vG_{\bw}^i
    \bSigma_{\bw^i}
    \vG_{\bw}^i {}\T (\vG_u^i \vK^i + \vI) \T \vP_{\tilde{k}}^i {}^T
    \nonumber
\end{align}
%
% The derivation of the above expression can be found in the Appendix\ref{Appendix: Sec2 Derivation of State Covariance}. 
The covariance $\Cov(\bx^i_{\tilde{k}})$ does not depend on the deterministic disturbance $\bzeta^i$, and we can rewrite the constraint \eqref{covariance constraint: interm1} as
\begin{align}
    \vP_{\tilde{k}}^i (\vG_u^i \vK^i + \vI) \vG_{\bw}^i
    \bSigma_{\bw^i}
    \vG_{\bw}^i {}\T (\vG_u^i \vK^i + \vI) \T \vP_{\tilde{k}}^i {}^T
    \preceq \bSigma^i_{\tilde{k}} \nonumber
\end{align}
We can reformulate the above non-convex constraint using Schur-Complement to the following convex constraint
% The above constraint is a non-convex constraint, and can be reformulated to the following convex constraint using Schur-Complement
% which can be equivalently given by the following convex constraint using Schur-Complement -
\begin{equation}
    \begin{bmatrix}
        \bSigma^i_{\tilde{k}} &  (\vP_{\tilde{k}}^i \vG_u^i \vK^i + \vI) \vG_{\bw}^i
        \boldsymbol{\varphi}^i \\
        (\vP_{\tilde{k}}^i(\vG_u^i \vK^i + \vI) \vG_{\bw}^i
        \boldsymbol{\varphi}^i) \T & \vI_{(n_x + N n_w)}
    \end{bmatrix} \succeq 0. \label{covariance_constraint_convex}
\end{equation}
\section{Distributed Robust Optimization Framework}
\label{sec: distr section distr alg subsec}
In this section, we introduce a distributed framework to solve the tractable version of Problem \ref{initial problem} with the reformulated constraints in a decentralized manner. Throughout this section, we consider the deterministic case for the analysis, though the framework applies to mixed case as well. 
%

% \subsection{Tractable Version of the Problem}
We begin by writing the tractable version of Problem \ref{initial problem} for the deterministic case with reformulated constraints.
\begin{problem}[Deterministic Case: MARTO with Reformulated Constraints]
\label{Transformed problem with reformulated constraints}
For all agents $i \in \calV$, find the robust optimal policies $\bar{\bu}^i,  \vK^i$, such that 
%
% \begin{equation}
% \begin{aligned}
% % \{ \bar{\bu}^i, \vK^i, \{ \bmu_{d,k}^i \}_{k = 0}^{T} \}_{i \in \calV} = & \argmin
% % \sum_{i \in \calV} 
% % J_i(\bar{\bu}^i, \vK^i)
% \min_{\{ \bar{\bu}^i, \vK^i, \{ \bmu_{d,k}^i \}_{k = 0}^{T} \}_{i \in \calV}} &
% \quad \sum_{i \in \calV} 
% J_i(\bar{\bu}^i, \vK^i)
% \nonumber
% \\[0.1cm]
% \mathrm{s.t.} \; \forall \; i \in \calV, \quad
% & (a): \; \eqref{chance constraints 1}, \; \eqref{chance constraints 2},
%  \\ 
% & (b): \; \eqref{covariance_constraint_convex}, 
% \\
% & (c): \; \eqref{Expectation of norm constraints Reformulated: mixed alpha}, \; 
% \eqref{Expectation of norm constraints Reformulated: mixed SDP}, 
% \\
% & (d): \; \eqref{mixed: nonconvex obstacle final ubar}, \; \eqref{mixed: nonconvex chance mu_d upper bound final}, \; \eqref{mixed: nonconvex chance mu_d lower bound final},
% \\
% & (e): \; \eqref{mixed: inter-agent collision avoidance ubar}, \; \eqref{mixed: nonconvex interagent mu_d upperbound} \quad \forall \; j \in \calN_i.
% \end{aligned}
% \end{equation}
% \end{problem}
%
\begin{equation}
\begin{aligned}
\min_{\{ \bar{\bu}^i, \vK^i, \{ \bmu_{d,k}^i \}_{k = 0}^{T} \}_{i \in \calV}} &
\quad \sum_{i \in \calV} 
J_i(\bar{\bu}^i, \vK^i)
\nonumber
\\[0.1cm]
\mathrm{s.t.} \; \forall \; i \in \calV, \quad
& (i): \; \eqref{mean_constraint_form1}, \; \eqref{New mean formulation}
 \\ 
& (ii): \; \eqref{nonconvex obstacle final ubar}, \; \eqref{semi-definite obstacle new},\; \eqref{mu_d expression}
\\
& (iii): \; \eqref{inter-agent collision avoidance ubar},\; \eqref{semi-definite interagent new}, ~~ \forall \; j \in \calN_i
\\
& (iv): \; \eqref{Convex norm-of-mean constraints: final 1}, \; \eqref{Convex norm-of-mean constraints: final 2}
\end{aligned}
\end{equation}
\end{problem}
Recall that we consider the linearized versions of the nonconvex constraints \eqref{nonconvex obstacle final ubar} and \eqref{inter-agent collision avoidance ubar}, and the nominal trajectory around which they are linearized will be disclosed later in this section. We will now transform the above problem into a form that allows us to solve it in decentralized way.
\subsection{Problem Transformation}
To develop a distributed framework, we need to focus on the inter-agent constraints \eqref{inter-agent collision avoidance ubar} and \eqref{semi-definite interagent new} involved in the problem.
To enforce these inter-agent constraints, each agent $i \in \calV $ needs $\{ \hat{\bmu}^i_{\text{pos},k}, \bmu_{d,k}^j \}_{k=0}^T$ from each of its neighbor agents $j \in \calN_i$.  
For convenience, let us define the sequences 
$ \hat{\bmu}^{i}_{\text{pos}} = [
    \hat{\bmu}^i_{\text{pos},1}; 
    \hat{\bmu}^i_{\text{pos},2}; 
    \dots;
    \hat{\bmu}^i_{\text{pos},T}]$, 
$\bmu^{i}_d = [
\bmu_{d,0}^i; 
\bmu_{d,1}^i; 
\dots; 
\bmu_{d,T}^i]$. 
We can now rewrite Problem \ref{Transformed problem with reformulated constraints} in a simplified form as follows
% in a simplified way as 
%
\begin{equation}
\begin{aligned}
    % \{ \bar{\bu}^i, \vK^i, \bmu_{d}^i \}_{i \in \calV} = \argmin
    % & \;
    % \sum_{i \in \calV} 
    % \hat{J}_i(\bar{\bu}^i,  \vK^i, \bmu_{d}^i)  
    \min_{\{ \bar{\bu}^i, \vK^i, \bmu_{d}^i \}_{i \in \calV} }
    \;
    \sum_{i \in \calV} 
    & \hat{J}_i(\bar{\bu}^i,  \vK^i, \bmu_{d}^i) 
    \\
    \mathrm{s.t.} \quad  
    \forall  \; i \in \calV, \;
    &\forall \; j \in \calN_i,  \;
    % \\
    % & 
    \calH_{i,j}( \hat{\bmu}^i_{\text{pos}} , \hat{\bmu}^j_{\text{pos}}) \geq 0 
    \\
    & \qquad \qquad \quad \Tilde{\calH}_{i,j}(\bmu_{d}^i, \bmu_{d}^j ) \geq 0
\end{aligned}
\label{Problem: Reformulated Simplified form 1}
\end{equation}
where $\hat{J}_i(\bar{\bu}^i,  \vK^i, \bmu_{d}^i) =  J_i(\bar{\bu}^i, \vK^i) + \calI_{X^i}$; $X^i$ represents the feasible region for the constraints  \eqref{mean_constraint_form1}, \eqref{New mean formulation}, \eqref{semi-definite obstacle new}, \eqref{mu_d expression}, \eqref{Convex norm-of-mean constraints: final 1}, \eqref{Convex norm-of-mean constraints: final 2}, and linearized version of \eqref{nonconvex obstacle final ubar}.
The constraints $\Tilde{\calH}_{i,j} \geq 0$ and $\calH_{i,j} \geq 0$ represent the constraints \eqref{semi-definite interagent new}, and linearized version of \eqref{inter-agent collision avoidance ubar} respectively.
%
%
% We then solve this problem using the global consensus technique \cite{boyd2011distributed} in a decentralized way. 
We then introduce, at each agent $i$, the copy variables $ \Bar{\bmu}^i_{j, \text{pos}}, \Bar{\bmu}^i_{j,d}$ corresponding to the variables $\hat{\bmu}^j_{\text{pos}}, \bmu_d^{j}$ of each of its neighbor agents $j \in \calN_i$ to enforce the inter-agent constraints. This necessitates the introduction of consensus constraints between these copy variables and their corresponding actual state variables. Thus, the problem of the form \eqref{Problem: Reformulated Simplified form 1} can be equivalently given as 
\begin{subequations}
\begin{align}
    % \{ \bar{\bu}^i, \vK^i, \bmu_{d}^i \}_{i \in \calV} = \argmin
    % & \;
    % \sum_{i \in \calV} 
    % \hat{J}_i(\bar{\bu}^i,  \vK^i, \bmu_{d}^i)  
    \min_{ \{ \bar{\bu}^i, \vK^i, \bmu_{d}^i \}_{i \in \calV} }
    & \quad
    \sum_{i \in \calV} 
    \hat{J}_i(\bar{\bu}^i,  \vK^i, \bmu_{d}^i)  \nonumber
    \\
    \mathrm{s.t.} \quad  
    \forall & \; i \in \calV, \;  \nonumber
    \\
    & 
    \calH_{i,j}( \hat{\bmu}^i_{\text{pos}} , \Bar{\bmu}^i_{j, \text{pos}} ) \geq 0,
    \; \forall \; j \in \calN_i, 
    \label{Problem simp2 const: hij}
    \\
    & \Tilde{\calH}_{i,j}(\bmu_{d}^i, \Bar{\bmu}^i_{j,d} ) \geq 0, 
    \; \forall \; j \in \calN_i, 
    \label{Problem simp2 const: htildeij}
    \\
    &
    \hat{\bmu}^i_{\text{pos}} = \Bar{\bmu}^{\hat{j}}_{i, \text{pos}}, 
    \;
    % \forall \; \hat{j} \in \calP_i,
    % \label{Problem simp2 const: consensus pos}
    % \\
    % &
    \bmu_d^{i} = \Bar{\bmu}^{\hat{j}}_{i,d},
    \;
    \forall \; \hat{j} \in \calP_i
    \label{Problem simp2 const: consensus mu_d}
\end{align}
\label{Problem: Reformulated Simplified form 2}
\end{subequations}
%
% This necessitates the introduction of consensus constraints between the actual state variables $\hat{\bmu}^{i}_{\bar{u}}, \bmu^{i}_{d} $ and their corresponding copy variables $\Bar{\bmu}^{\hat{j}}_{i, \bar{u}} $, $\Bar{\bmu}^{\hat{j}}_{i, d} $ respectively at each agent $\hat{j}$ to which agent $i$ is a neighbor of ($\hat{j} \in \calP_i$). 
% We are interested in solving the above problem in a decentralized manner. For that,
% 
To solve the above problem in a decentralized way, we convert it into a 'global consensus problem' \cite{boyd2011distributed} by introducing the global variables $\{ \bnu^i_{\text{pos}}, \bnu^i_{d} \}_{i \in \calV}$ such that the consensus constraints
% \eqref{Problem simp2 const: consensus pos}, and 
\eqref{Problem simp2 const: consensus mu_d} for each $i \in \calV$ are given as  
\begin{equation}
\begin{aligned}
    \forall \hat{j} \in \calP_i, ~~ 
    \hat{\bmu}^i_{\text{pos}} = \bnu^i_{\text{pos}}, 
    ~ \Bar{\bmu}^{\hat{j}}_{i, \text{pos}} = \bnu^i_{\text{pos}}, 
    ~ \bmu^{i}_d =  \bnu^i_{d} , 
    ~ \Bar{\bmu}^{\hat{j}}_{i,d}  =  \bnu^i_{d} 
\end{aligned}
\label{consensus relation -1}
\end{equation}
%
% For convenience, we define the global variable sequences as $\bnu_{\text{pos}} = [\bnu^1_{\text{pos}}; \bnu^2_{\text{pos}}; \dots; \bnu^N_{\text{pos}}]$, 
% $\bnu_{d} = [\bnu^1_{d}; \bnu^2_{d}; \dots; \bnu^N_{d}]$, and the local variables of each agent  $i \in \calV$ as  
% $ \Tilde{\bmu}_{\text{pos}}^i  = [
%      \hat{\bmu}^{i}_{\text{pos}}; 
%     \{  \Bar{\bmu}^i_{j, \text{pos} } \}_{j \in \calN_i} ]$, 
% $\Tilde{\bmu}_{d}^i = [
%         \bmu^{i}_d; 
%         \{ \Bar{\bmu}^i_{j,d} \}_{j \in \calN_i}]$. 
% We also define the matrices $\{ \mathrm{\vC}_g^i \}_{i \in \calV}$ such that the consensus constraints \eqref{consensus relation -1} can be equivalently given as
% \begin{equation}
%     \Tilde{\bmu}_{\text{pos}}^i =  \mathrm{\vC}_g^i \bnu_{\text{pos}},
%     \quad 
%     \Tilde{\bmu}_{d}^i = \mathrm{\vC}_g^i \bnu_{d} 
%     \quad 
%     \forall i \in \calV.
% \end{equation}
Let us define the local variable of each agent  $i \in \calV$ as  
% $ \Tilde{\bmu}^i_{\text{pos}}  = [
%      \hat{\bmu}^{i}_{\text{pos}}; 
%     \{  \Bar{\bmu}^i_{j, \text{pos} } \}_{j \in \calN_i} ]$, 
% $\Tilde{\bmu}_{d}^i = [
%         \bmu^{i}_d; 
%         \{ \Bar{\bmu}^i_{j,d} \}_{j \in \calN_i}]$, and the global variables sequences as
% $\Tilde{\bnu}_{\text{pos}}^i  = [
%      \bnu^{i}_{\text{pos}}; 
%     \{  \bnu^j_{\text{pos} } \}_{j \in \calN_i} ]$, 
% $\Tilde{\bnu}_{d}^i = [
%         \bnu^{i}_d; 
%         \{ \bnu^j_{d} \}_{j \in \calN_i}]$,
$ \Tilde{\bmu}^i  = [
     \hat{\bmu}^{i}_{\text{pos}}; 
    \{  \Bar{\bmu}^i_{j, \text{pos} } \}_{j \in \calN_i};
    \bmu^{i}_d; 
    \{ \Bar{\bmu}^i_{j,d} \}_{j \in \calN_i} ]$, and the global variable as 
    $\bnu^i = [\bnu^i_{\text{pos}}; \bnu^i_{d} ]$. Further, by defining the global variable sequences
$\Tilde{\bnu}^i  = [
     \bnu^{i}_{\text{pos}}; 
    \{  \bnu^j_{\text{pos} } \}_{j \in \calN_i};
    \bnu^{i}_d; 
    \{ \bnu^j_{d} \}_{j \in \calN_i}]$, 
    the consensus constraints \eqref{consensus relation -1} can be equivalently given as 
\begin{equation}
% \forall i \in \calV, \quad
    \mathrm{\vC} \Tilde{\bmu}^i = \mathrm{\vC} \tilde{\bnu}^i
    % \Tilde{\bmu}^i_{\text{pos}} = \Tilde{\bnu}_{\text{pos}}^i,
    % ~~~
    % \Tilde{\bmu}_{d}^i = \Tilde{\bnu}_{d}^i.
\end{equation}
where $\mathrm{\vC} = \bdiag (\vI, \rho_r \vI )$ and $\rho_r >0$.
%
% For convenience, let us also define $\{ \bnu^i \}_{i \in \calV}$ as $\bnu^i = [\bnu^i_{\text{pos}}; \bnu^i_{d} ]$.
Thereby, we can rewrite the problem in \eqref{Problem: Reformulated Simplified form 2} as the following global consensus problem.
\begin{problem}[Deterministic Case: MARTO Consensus Problem]
\label{Problem 1 distributed formulation}
For all agents $i \in \calV$, find the robust optimal policies $\bar{\bu}^i,  \vK^i$, such that 
\begin{equation}
\begin{aligned}
    & \min_{ \{ \bar{\bu}^i, \vK^i, \Tilde{\bmu}^i, \bnu^i \}_{i \in \calV} } 
    \quad 
    \sum_{i \in \calV}
    \hat{\calJ}_i(\bar{\bu}^i, \vK^i,  \Tilde{\bmu}^i)
    \\ 
    &~~~~~~~~
    \mathrm{s.t.} \; \forall  \; i \in \calV, 
    \quad  \mathrm{\vC} \Tilde{\bmu}^i = \mathrm{\vC} \tilde{\bnu}^i
\end{aligned}
\end{equation}
where $\hat{\calJ}_i ( \bar{\bu}^i, \vK^i,  \Tilde{\bmu}^i) = \hat{J}_i(\bar{\bu}^i,  \vK^i, \bmu_{d}^i)  + \sum_{j \in \calN_i} \calI_{ \hat{X}^{ij}  } $;  $\hat{X}^{ij} $ represent the feasible region for the constraints \eqref{Problem simp2 const: hij}, \eqref{Problem simp2 const: htildeij}.
\end{problem}
\subsection{Distributed Framework}
We now provide a distributed framework in Algorithm \ref{alg1} to solve Problem \ref{Problem 1 distributed formulation} based on Consensus Alternating Direction Method of Multipliers (CADMM) \cite{boyd2011distributed}, which is a variant of Augmented Lagrangian (AL) methods. 
% We now use CADMM to solve the above problem in a distributed way as disclosed in the algorithm \ref{alg1}. 
Thus, we begin by writing the AL of Problem \ref{Problem 1 distributed formulation} as follows
\begin{equation}
\begin{aligned}
     \calL_{\rho} ( 
     \{ \Tilde{\bmu}^i, 
        \bar{\bu}^i, \vK^i,
        \bnu^i,
        \blambda^{i} \}_{i \in \calV} )= 
     \sum_{i \in \calV} \calL_{\rho}^i ( \Tilde{\bmu}^i, \bar{\bu}^i, \vK^i,
     \Tilde{\bnu}^i , \blambda^{i}) \nonumber
    \end{aligned}
    % \label{Distributed Algorithm: Lagrangian}
\end{equation}
where $ \{ \blambda^i \}_{i \in \calV}$ are dual variables, $\rho > 0$ is penalty parameter, % 
\begin{equation}
\begin{aligned}
    \calL_{\rho}^i ( \Tilde{\bmu}^i, \bar{\bu}^i, \vK^i,
     \Tilde{\bnu}^i , \blambda^{i} )
    & =
    \hat{\calJ}_i(\bar{\bu}^i, \vK^i, \Tilde{\bmu}^i )
    + \blambda^i {}\T \mathrm{\vC} ( \Tilde{\bmu}^i - \Tilde{\bnu}^i)
    \\
    &~~~~~~~
    + \frac{\rho}{2} || \mathrm{\vC}( \Tilde{\bmu}^i - \Tilde{\bnu}^i) ||_2^2.  \nonumber
\end{aligned}
\end{equation}
% and $ \{ \blambda^i \}_{i \in \calV}$ are dual variables, $\rho > 0$ is penalty parameter. 

Let the local variables $\{\Tilde{\bmu}^i \}_{i \in \calV}$ constitute the first block, and the global variables $\{ \bnu^i \}_{i \in \calV}$ constitute the second block of CADMM. 
Each iteration of CADMM involves sequential minimization of $\calL_{\rho}$ with respect to each block, followed by a dual update. 
Since our problem involves non-convex constraints \eqref{nonconvex obstacle final ubar} and \eqref{inter-agent collision avoidance ubar}, 
% standard implementation of CADMM might not guarantee convergence. Thus, 
we consider a discounted dual update step \cite{yang2022proximal} to ensure convergence (the details of which are disclosed in Appendix\ref{Appendix sec: Convergence Analysis}). 
Each $(l+1)^{th}$ CADMM iteration of Algorithm \ref{alg1} involves the following steps
\begin{enumerate}
    \item Local Update: 
    \begin{equation}
    \begin{aligned}
         \{ \Tilde{\bmu}^{i} \}^{l+1}_{i \in \calV} = \argmin \calL_{\rho} ( 
        \{ \Tilde{\bmu}^i, \bar{\bu}^i, \vK^i,
        \bnu^{i,l},
        \blambda^{i,l} \}_{i \in \calV} ). \nonumber
    \end{aligned}
    \end{equation}
    There is no coupling the agents to carryout the above update, thus it is done in a decentralized manner with each agent $i$ solving the following problem 
    \begin{equation}
    \begin{aligned}
        \Tilde{\bmu}^{i,l+1}  = \argmin \calL_{\rho}^i ( \Tilde{\bmu}^i, \bar{\bu}^i, \vK^i,
     \Tilde{\bnu}^{i,l} , \blambda^{i,l}),
    \end{aligned}
    \label{Algorithm: ADMM local update}
    \end{equation}
    by receiving the global variable $\bnu^{j,l}$ from all $j \in \calN_i$. 
    \item Global Update:
    \begin{equation}
        \{ \bnu^{i,l+1} \}_{i \in \calV} = \argmin \calL_{\rho} ( 
        \{ \Tilde{\bmu}^{i,l+1},
        \bnu^{i},
        \blambda^{i,l} \}_{i \in \calV} ). \nonumber
    \end{equation}
    The above update can also be carried out in a decentralized manner with each agent $i$ carrying out the following updates in parallel
    \begin{equation}
    \begin{aligned}
    &  \bnu^{i, l+1}_{\text{pos}}
    =
    \frac{1}{1 + n(\calP_i)} \bigg( \hat{\bmu}^{i, l+1}_{\text{pos}} 
    + \sum_{j \in \calP_i} \Bar{\bmu}^{\hat{j},l+1}_{i, \text{pos}} \bigg),  
    \\
    &
     \bnu^{i, l+1}_{d}
    =
    \frac{1}{1 + n(\calP_i)} \bigg( \hat{\bmu}^{i, l+1}_{d} 
    + \sum_{j \in \calP_i} \Bar{\bmu}^{\hat{j},l+1}_{i, d} \bigg),
    \end{aligned}
    \label{Algorithm: ADMM global update}
    \end{equation}
    by receiving the local variables updates $\bar{\bmu}_{i,d}^{\hat{j}, l+1}$ from all the agents $ \hat{j}  \in \calP_i$ (derivation in Appendix\ref{Appendix: sec D Derivation of Global Updates}). 
    % The details of derivation of the above update are provided in the Appendix\ref{Appendix: sec D Derivation of Global Updates}.
    %
    %
    \item Discounted Dual Update: 
    Each agent $i$ carries out the following discounted dual update by receiving the global variable updates $  \bnu_{\text{pos}}^{j, l+1}, \bnu_{d}^{j,l+1}$ from all $j \in \calN_i$,
    \begin{equation}
        \blambda^{i,l+1} 
        = (1- \delta) \blambda^{i,l}
        + \rho \mathrm{\vC} ( \Tilde{\bmu}^{i,l+1} - \Tilde{\bnu}^{i,l+1} ),
        \label{Algorithm: Discounted Dual update}
    \end{equation}
    where $\delta$ is the discounting parameter.
\end{enumerate}
The above iterations are carried out until the following convergence criteria, where both primal and dual residuals fall below the corresponding set thresholds $\epsilon_{\text{p}}$ and $\epsilon_{\text{d}}$, is fulfilled. 
\begin{equation}
\begin{aligned}
    & \bigg[
    \frac{ \sum_{i=1}^N  \| \Tilde{\bmu}^{\bar{u}, l}_i - \Tilde{\bnu}^{\bar{u},l}_i \|^2 }
    { N + \sum_{i = 1}^{N} n(\calN_i) } \bigg]^{1/2} \leq \epsilon_{\text{p}},
    \\
    &
    \frac{ \big( \sum_{i=1}^N  \rho^2 \| \bnu^{l}_i - \bnu^{l-1}_i \|^2 \big)^{1/2} }{ N }  \leq \epsilon_{\text{d}}. 
\end{aligned}   
\label{Convergence conditions}
\end{equation}
\begin{remark}
    The agents need not store system information of their neighbor agents'  including dynamics and uncertainty model parameters, to carry out the above updates; thereby, this framework preserves security.
\end{remark}
%
%
%
% Further, we define $\blambda^i_{\text{pos}} $, $\blambda^i_{d} $ as sequence of dual vectors given as 
% %
% \mbox{$\blambda^i_{\text{pos}}  = \begin{bmatrix}
%     \bar{\blambda}^i_{\text{pos}}; &
%     \{ \bar{\blambda}^{i}_{j, \text{pos}} \}_{j \in \calN_i}
% \end{bmatrix} $},
% %
% \mbox{
% $ \blambda^i_{d}  = \begin{bmatrix}
%      \bar{\blambda}^i_{d} ; &
%     \{ \bar{\blambda}^{i}_{j, d} \}_{j \in \calN_i}
% \end{bmatrix}$}.
%
%
% \begin{remark}
%     As the considered problem involves non-convex constraints, standard implementation of Consensus ADMM might not guarantee convergence. Thus, we consider a discounted dual update step to ensure convergence (the details of which are disclosed in the section \ref{Convergence Analysis}).
%     % Consensus ADMM is guaranteed to be converging for the convex problems. Thus, to establish convergence guarantees in the presence of non-convex constraints, we consider a discounted dual update step. 
% \end{remark}
%
%
\begin{algorithm}[H]
\caption{\strut Distributed Robust Trajectory Optimization}
\begin{algorithmic}[1]
\STATE
Set $\delta \in (0, 0.5]$, $\rho, \rho_r > 0$,
max ADMM iterations $L_{max}$.
\\
\FOR{$l = 0$ to $L_{max}$} 
\STATE \textit{Local update:} $ \forall \; i \in \calV$ (in parallel), \eqref{Algorithm: ADMM local update}
\vskip -2.5pt
\STATE 
\textit{Global update:} $ \forall \; i \in \calV$ (in parallel), \eqref{Algorithm: ADMM global update}
\vskip -2.5pt
\STATE
\textit{Discounted Dual update: } $ \forall \; i \in \calV$ (in parallel), \eqref{Algorithm: Discounted Dual update}
\vskip -3.5pt
\STATE \textbf{if} \eqref{Convergence conditions} \textbf{then} break
\ENDFOR
\end{algorithmic}
\label{alg1}
\end{algorithm}
Assuming that the local sub-problems (in \eqref{Algorithm: ADMM local update}) are feasible and $\delta$ lies in $(0,0.5]$, we study the convergence properties of Algorithm \ref{alg1} in Appendix\ref{Appendix sec: Convergence Analysis}. 
To outline the analysis, we construct a merit function in a systematic manner. We start by deriving upper-bound of the difference of a regularized Lagrangian function in each iteration, which depends on the change in dual variables ($\| \blambda^{l+1} - \blambda^l \|_2^2$). Then, we derive an upper bound of dual variables in terms of local and global variables. By combining the above results, we derive a descent relation and propose a merit function.  
Based on an empirical assumption, we prove Algorithm \ref{alg1} convergence to $\delta \rho^{-1} \| \blambda^* \|_2$ approximate stationary points (as defined in \cite{yang2022proximal}) of Problem \ref{Problem 1 distributed formulation} where $\blambda^*$ is the limit point of the algorithm.
%
% we first upper bound the change in dual variables in each iteration (i.e., $\| \blambda^{l+1} - \blambda^l \|_2^2$), which is achieved due to discounted dual update. 
% We propose a merit function and based on an empirical assumption that there exists a finite ADMM iteration $\bar{l}$ such that $\alpha_{\delta} r_{c}^l + r_h^l \leq 0$ for all $l \geq \bar{l}$ (defined in Appendix\ref{Appendix sec: Convergence Analysis}), we show that Algorithm \ref{alg1} converges to $\delta \rho^{-1} \| \blambda^* \|_2$ approximate stationary points (as defined in \cite{yang2022proximal}) of Problem \ref{Problem 1 distributed formulation} where $\blambda^*$ is the limit point of the algorithm.
% \begin{enumerate}
%     \item 
%     \item 
% \end{enumerate}
%
\section{Computational Complexity Analysis}
\label{sec: Computational Complexity Analysis}
In this section, we provide an analysis of the notable computational performance of the proposed distributed framework through two key comparisons as provided below. For simplicity of presentation, we consider the deterministic case, and  define $\hat{n}_{\text{obs}} = n_{u_i} \gamma_h n_{x_i} + n_{\text{obs}}, \hat{n}_{\text{com}} = n_{u_i} \gamma_h n_{x_i} + n_{\text{obs}} + n_{\text{inter}} $.
\subsubsection{Efficacy of the Proposed Constraint Reformulation}
\label{Complexity subsec: Efficacy of the Proposed Constraint Reformulation}
First, we highlight the efficacy of the proposed constraint reformulation compared to the SDP approach (as disclosed in section \ref{initial approach}). 
\begin{proposition}
    The computational complexity of solving a single-agent robust optimization problem using the proposed constraint reformulation is given as 
% \begin{equation}
% \begin{aligned}
%     & O \bigg( T^{7/2} n_{\text{obs}}^{1/2}
%     (n_{u_i} \gamma_h n_{x_i} + n_{\text{obs}})
%     \\
%     & \qquad \qquad \quad
%     \big[ (n_{u_i} \gamma_h n_{x_i} + n_{\text{obs}})^2 
%      + n_{d_i}^2 T
%     \big] \bigg), \nonumber
% \end{aligned}
% \end{equation}
\begin{equation}
\begin{aligned}
     O \big( T^{7/2} n_{\text{obs}}^{1/2} \hat{n}_{\text{obs}}
    \big[ \hat{n}_{\text{obs}}^2 
     + n_{d_i}^2 T
    \big] \big), \nonumber
\end{aligned}
\end{equation}
which is significantly better compared to the computational complexity of the SDP approach given as  
% \begin{equation}
% \begin{aligned}
%     O \bigg( 
%     T^6 n_{d_i}^{7/2} n_{\text{obs}}^{1/2} \big[ n_{u_i} \gamma_h n_{x_i} + n_{obs} \big]^3 
%     \bigg).  \nonumber
% \end{aligned}
% \end{equation}
\begin{equation}
\begin{aligned}
    O \big( 
    T^6 n_{d_i}^{7/2} n_{\text{obs}}^{1/2} \hat{n}_{\text{obs}}^3 
    \big).  \nonumber
\end{aligned}
\end{equation}
Furthermore, our proposed constraint reformulation is particularly useful in a multi-agent setting with the computational complexity of the proposed distributed framework given as 
% \begin{equation}
% \begin{aligned}
%     & O \bigg( L_{\text{ADMM}} T^{7/2} ( n_{\text{obs}} + n_{\text{inter}})^{1/2}
%     (n_{u_i} \gamma_h n_{x_i} + n_{\text{obs}} + n_{\text{inter}})
%     \\
%     & \qquad \quad
%     \big[ (n_{u_i} \gamma_h n_{x_i} + n_{\text{obs}} + n_{\text{inter}})^2 
%      + n_{d_i}^2 T
%     \big] \bigg), \nonumber
% \end{aligned}
% \end{equation}
\begin{equation}
\begin{aligned}
    & O \big( L_{\text{ADMM}} T^{7/2} ( n_{\text{obs}} + n_{\text{inter}})^{1/2}
    \hat{n}_{\text{com}}
    \big[ \hat{n}_{\text{com}}^2 
     + n_{d_i}^2 T
    \big] \big), \nonumber
\end{aligned}
\end{equation}
in contrast to the computational complexity of the distributed SDP approach given as 
% \begin{equation}
% \begin{aligned}
%     O \bigg( L_{\text{ADMM}} 
%     T^6 n_{d_i}^{7/2} \big[ n_{u_i} \gamma_h n_{x_i} n_{\text{inter}} + n_{obs} \big]^3 ( n_{\text{obs}} + n_{\text{inter}} )^{1/2}
%     \bigg)  \nonumber
% \end{aligned}
% \end{equation}
\begin{equation}
\begin{aligned}
    O \big( L_{\text{ADMM}} 
    T^6 n_{d_i}^{7/2} \big[ n_{u_i} \gamma_h n_{x_i} n_{\text{inter}} + n_{obs} \big]^3 ( n_{\text{obs}} + n_{\text{inter}} )^{1/2}
    \big)  \nonumber
\end{aligned}
\end{equation}
\end{proposition}
\proof The proof provided in the Appendix\ref{Appendix: sec computational complexity}.
\endproof
%
%
% It should be noted that the above complexity bounds provided are conservative, and in practice, our proposed framework has better convergence speed.
% Recall that the proposed constraint reformulation also enables the state of the agents to be characterized using $\bmu_{x,\bar{u}}^i, \bmu_{d}^i$, which is not the case for the 
%
\subsubsection{Scalability of the Proposed Distributed Framework}
Next, we emphasize the scalability of the proposed distributed framework by providing a computational complexity bound for solving Problem \ref{Transformed problem with reformulated constraints} in a centralized manner. 
\begin{proposition}
    The computational complexity for solving Problem \ref{Transformed problem with reformulated constraints} in a centralized manner, is given as 
% \begin{equation}
% \begin{aligned}
%     O \bigg( (NT)^{7/2} ( n_{\text{obs}} + n_{\text{inter}})^{1/2} \big[ n_{u_i} \gamma_h n_{x_i} 
%     + n_{\text{obs}} 
%     + n_{\text{inter}} \big]^3 \bigg) \nonumber
% \end{aligned}
% \end{equation}
\begin{equation}
\begin{aligned}
    O \big( (NT)^{7/2} ( n_{\text{obs}} + n_{\text{inter}})^{1/2} \hat{n}_{\text{com}}^3 \big) \nonumber
\end{aligned}
\end{equation}
\end{proposition}
\proof
    The proof is provided in the Appendix\ref{Appendix: sec computational complexity}.
\endproof
The above bound brings attention to a significant challenge in solving Problem \ref{Transformed problem with reformulated constraints} in centralized manner, wherein the computational complexity rapidly increases with the number of agents ($N$) at the rate of $N^{7/2}$, rendering it to be computationally intractable for large-scale multi-agent systems. 

Note that these bounds are conservative, and in practice, the proposed method performs considerably faster as demonstrated in the simulation results.
% \section{Convergence Analysis} \label{Convergence Analysis}
%
% Assuming that the local subproblem is feasible and the solution is bounded. $\delta_{\text{pos}} \in (0, 0.5]$.

\section{Simulation} \label{sec: Simulation}
In this section, we demonstrate the effectiveness and scalability of the proposed framework through experimental results. We set the time horizon to $T = 20$ unless stated otherwise. Each agent has 2D double integrator dynamics with matrices $\vC_k^i = \texttt{randn}(4, 2)$ normalized to $\| \vC_k^i \|_F = 1$, and $\vD_k^i = \vI$. For all agents, the uncertainty parameters are $\vS_i = \vI, \vGamma_i = \vI$, and different $\tau_i \in [0.01, 0.075]$ for each scenario. We set the ADMM parameters to $\delta = 0.001$, $\rho = 100$, $\rho_r = 0.15$, $\epsilon_{\text{p}} = 0.05$ and $\epsilon_{\text{d}} = 0.1$. 
Further, the probability for satisfying nonlinear inter-agent chance constraints is set to 0.81, and to 0.9 for the rest. For each agent $i$, the set $\calN_i$ is considered to be a set of specified number of agents closest to the agent at the initial state. 
% A supplementary video \footnote{}is provided showcasing all the experiments.
%
%
\subsubsection{Single Agent Computational Complexity Comparison}
We compare the optimizer time for solving a single agent problem using SDP approach and the proposed constraint reformulation (say "SOCP approach") for the deterministic case in Table \ref{Table-1}. A time horizon of $T= 15$ is considered. The following key observations are made.
First, the change in the optimizer time with the disturbance dimension $n_d$ or the number of obstacles is significantly more for SDP approach compared to the SOCP approach, which is inline with the complexity bounds derived in Section \ref{Complexity subsec: Efficacy of the Proposed Constraint Reformulation}. For more than two obstacles, the solver fails to solve using SDP approach. 
Second, for a set number of obstacles, the optimizer time varies greatly across different uncertainty levels for the SDP approach compared to the SOCP approach. Since the SDP approach fails to scale with the obstacles even in a single agent case, we omit the comparison of frameworks for multi-agent scenarios.
% First, the change in the optimizer time with the disturbance dimension $n_d$ is significantly more for SDP framework which is inline with the complexity bounds derived in Section \ref{Complexity subsec: Efficacy of the Proposed Constraint Reformulation}. Second, for a set number of obstacles, the optimizer time varies greatly across different uncertainty levels for the SDP framework compared to the SOCP framework. Third, the increase in the optimizer time for SDP framework wrt obstacles is significantly more compared to the SOCP framework. For more than two obstacles, the solver fails to solve using SDP approach. Since the SDP approach fails to scale wrt obstacles even in a single agent case, we omit the comparison of frameworks for multi-agent scenarios.
%
\begin{table}[t!]
    \centering
    \begin{tabular}{ |c|c|c|c|c|c| } 
\hline
 \multirow{2}{*}{$n_{\text{obs}}$} & \multirow{2}{*}{$\tau$} & \multicolumn{4}{|c|}{Optimizer time (sec)}\\ \cline{3-6}
 & & \multicolumn{2}{|c|}{$n_{d} = 2$}&  \multicolumn{2}{|c|}{$n_{d} = 4$} \\ \cline{3-6}
& & SOCP & SDP & SOCP & SDP \\
\hline
\multirow{6}{*}{$1$} & $0.01$ & $0.17$ & $3.37$ & $0.26$ & $6.71$\\
& $0.02$ & $0.18$ & $2.37$ & $0.24$ & $7.84$ \\
& $0.03$ & $0.17$ & $2.84$ & $0.2$ & $8.21$\\
& $0.04$ & $0.27$ & $2.7$ & $0.28$ & $16.65$ \\
& $0.05$ & $0.15$ & $3.43$ & $0.34$ & $9.16$ \\ \cline{2-6}
& avg. & $0.186$ & $2.942$ & $0.264$ & $9.714$\\
\hline
\multirow{6}{*}{$2$} & $0.01$ & $0.16$ & $6.93$ & $0.24$ & $16.44$ \\
& $0.02$ & $0.2$ & $8.74$ & $0.3$ & $31.3$ \\
& $0.03$ & $0.19$ & $4.78$ & $0.28$ & $32.53$ \\
& $0.04$ & $0.21$ & $4.64$ & $0.27$ & $14.30$ \\
& $0.05$ & $0.19$ & $4.5$ & $0.26$ & $15.58$ \\ \cline{2-6}
& avg. & $0.19$ & $5.918$ & $0.27$ & $22.03$\\
\hline 
$3$ & $0.01$ & 0.17 & - & 0.27 & - \\
\hline
\end{tabular}
\caption{Optimizer time comparison between the proposed SOCP approach (in Section \ref{new reformulation}) and the SDP approach (in Section \ref{initial approach}) for a single agent problem}
\label{Table-1}
\vspace{-0.5cm}
\end{table}
\subsubsection{Effectiveness of Proposed Constraint Reformulations}
\begin{figure}[t!] 
    \centering
  \subfloat[ Without robust constraints \label{2a}]{
       \includegraphics[width= 0.45\linewidth, trim={0.95cm 4.5cm 1.75cm 5cm},clip]{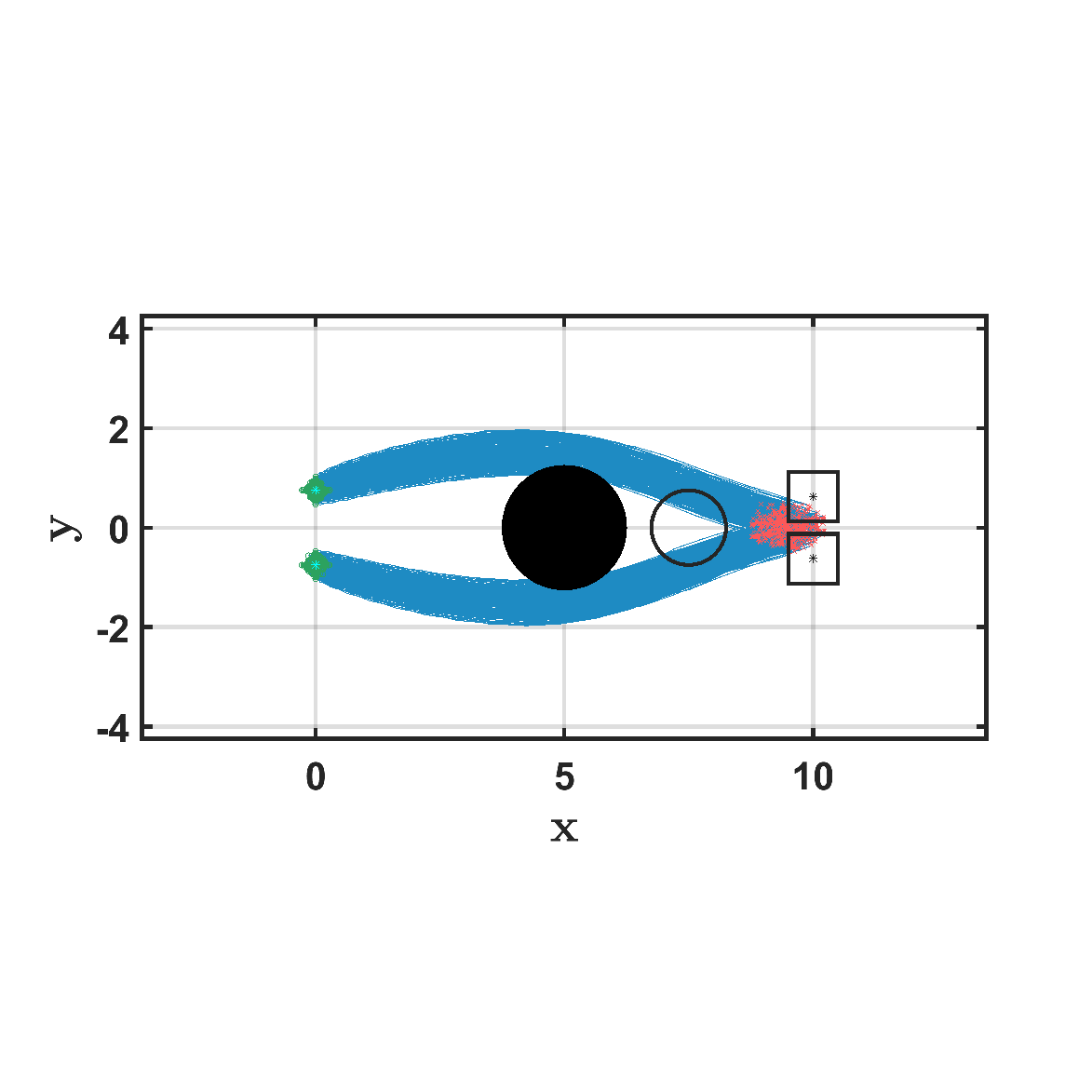}}
  \subfloat[ At $k= 15$ \label{2b}]{
        \includegraphics[width=0.45\linewidth, trim={0cm 0cm 0cm 0cm},clip]{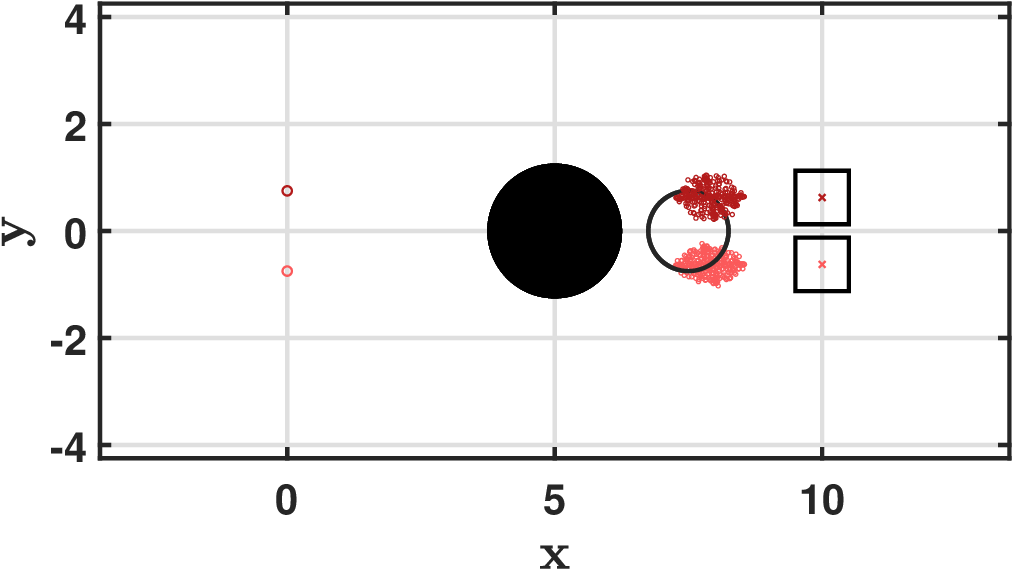}}
        \\
        \subfloat[ With robust constraints \label{2c}]{
       \includegraphics[width= 0.45\linewidth, trim={0.95cm 4.5cm 1.75cm 5cm},clip]{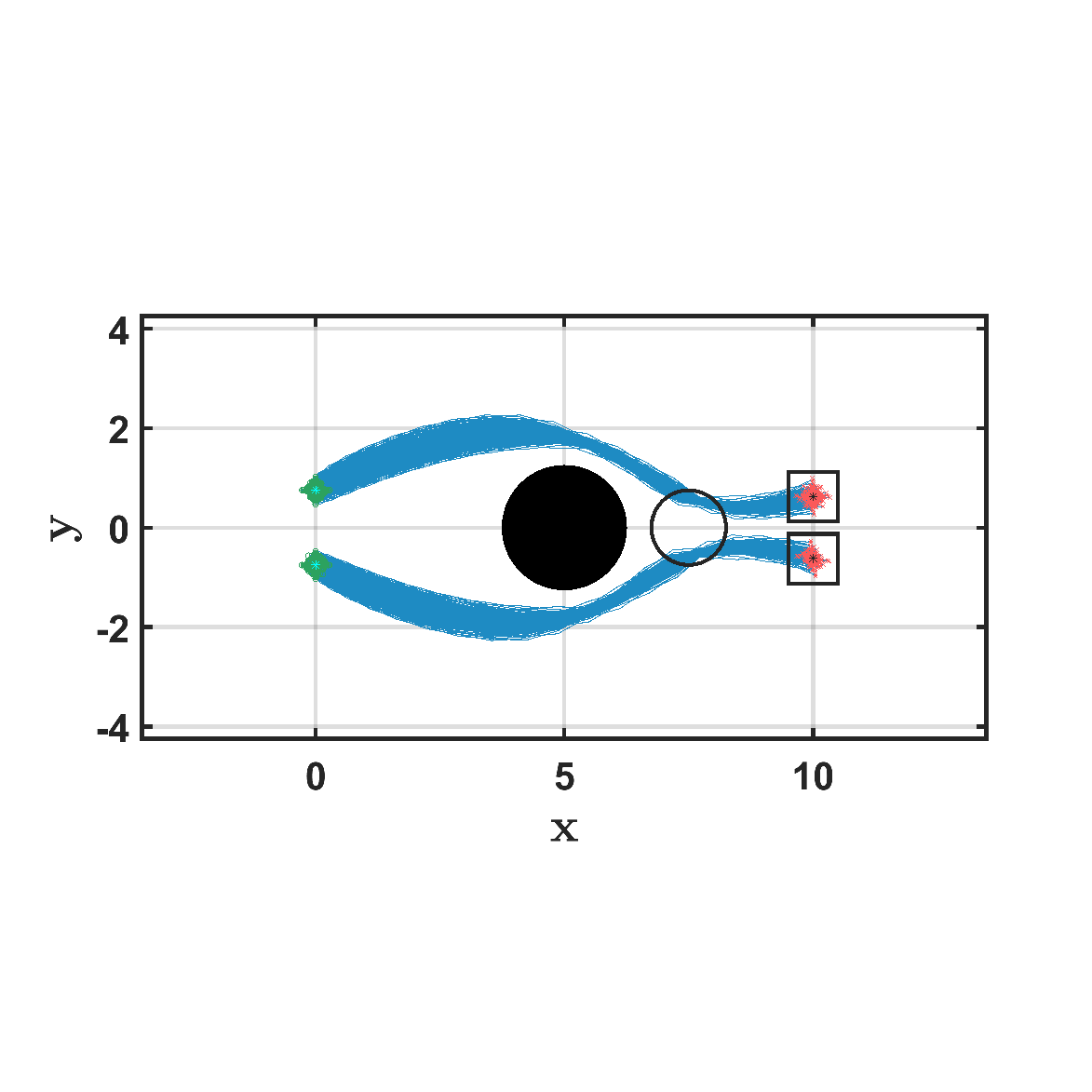}}
  \subfloat[At $k= 15$  \label{2d}]{
        \includegraphics[width=0.45\linewidth, trim={0cm 0cm 0cm 0cm},clip]{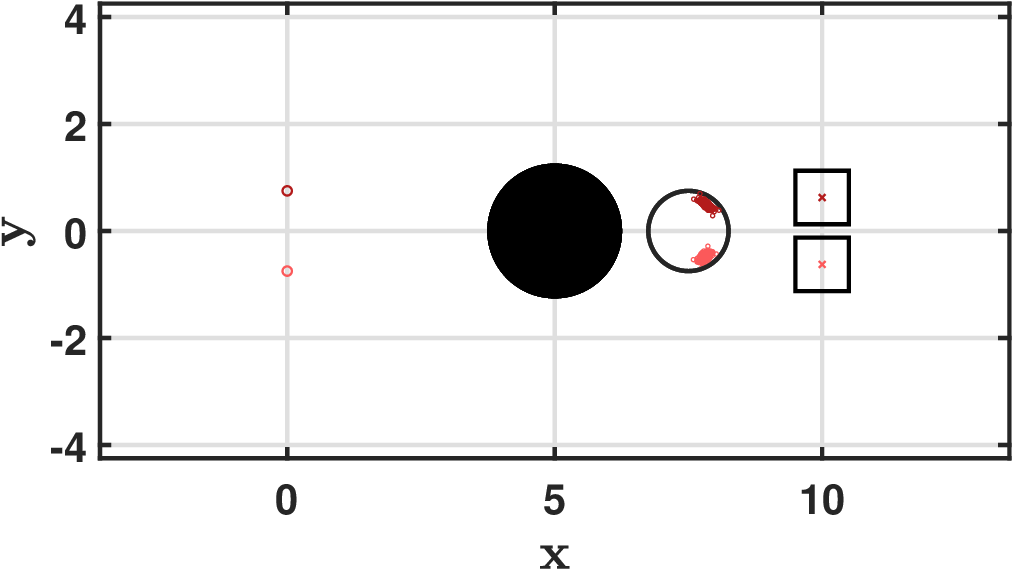}}
  \hfill
  \caption{ \textbf{Deterministic Case:} With terminal linear \eqref{Robust Linear Constraints}, collision avoidance \eqref{Robust nonconvex norm-of-mean constraints} \eqref{robust inter-agent collision avoidance} and convex norm constraints \eqref{Robust convex norm-of-mean constraints}. }
  \label{fig2} 
  \vspace{-0.25cm}
\end{figure}

We demonstrate the effectiveness of the robust constraints for the deterministic case (defined in Section \ref{sec: Robust Constraints for deterministic case}) in Fig. \ref{fig2}. In Fig. \ref{fig2}, two agents are tasked to reach the target bounds (shown as black boxes) while avoiding the circular obstacle, and remaining inside the black circle at a time step $k=15$. The Fig. \ref{2a}-\ref{2b} show the trajectory realizations of the agents without considering the robust constraints (i.e., finding the control policy without accounting for the deterministic uncertainty). The agents violate all the constraints. The Fig. \ref{2c}- \ref{2d} show that with the robust constraints, all the constraints are satisfied for all realizations of uncertainty. 

\begin{figure}[t!] 
    \centering
  \subfloat[ With constraints on $\Eb(\bx_i)$ \label{3a}]{
       \includegraphics[width= 0.45\linewidth, trim={0.95cm 4cm 1.75cm 5cm},clip]{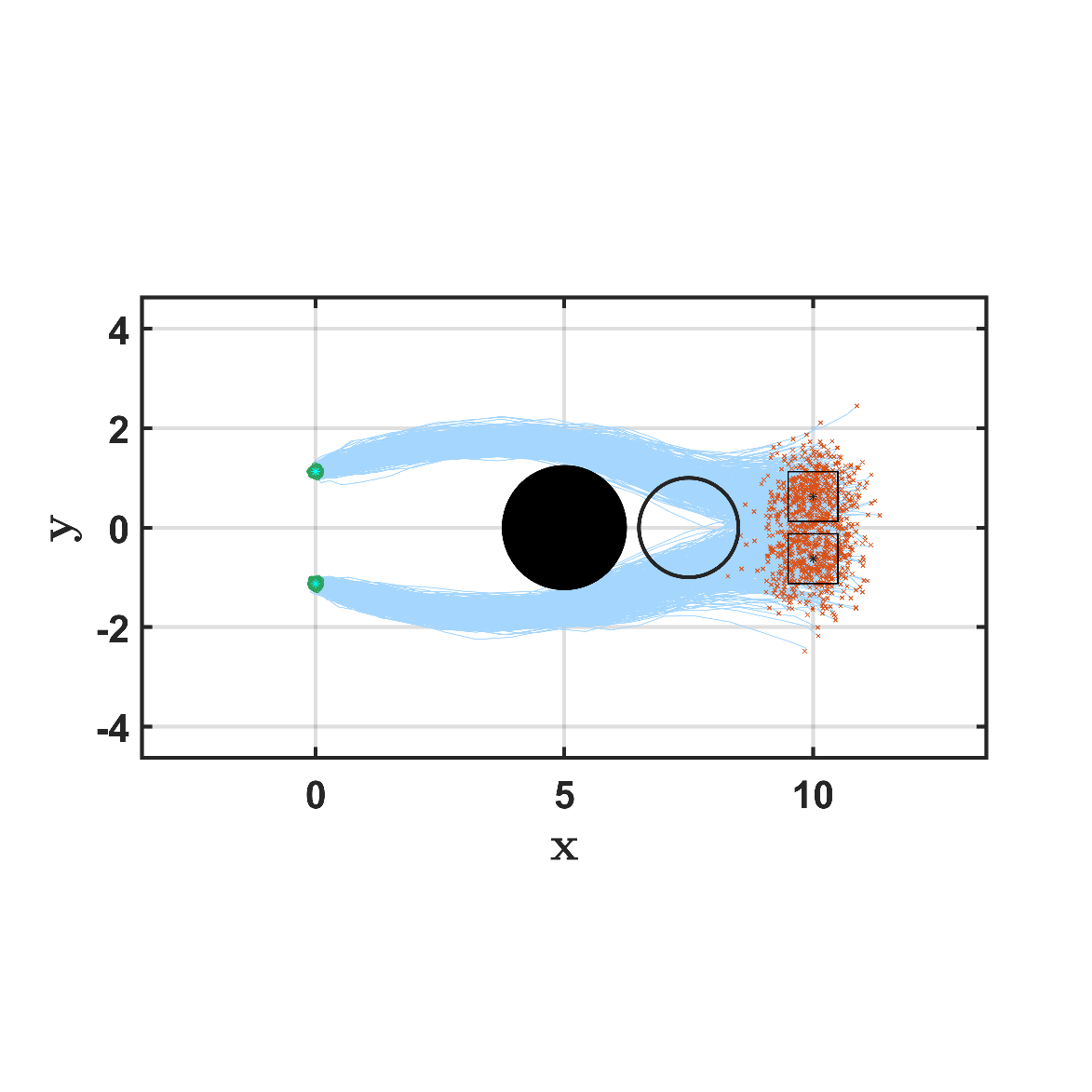}}
\subfloat[ At $k =15$ \label{3b}]{
       \includegraphics[width= 0.45\linewidth, trim={0.95cm 4cm 1.75cm 5cm},clip]{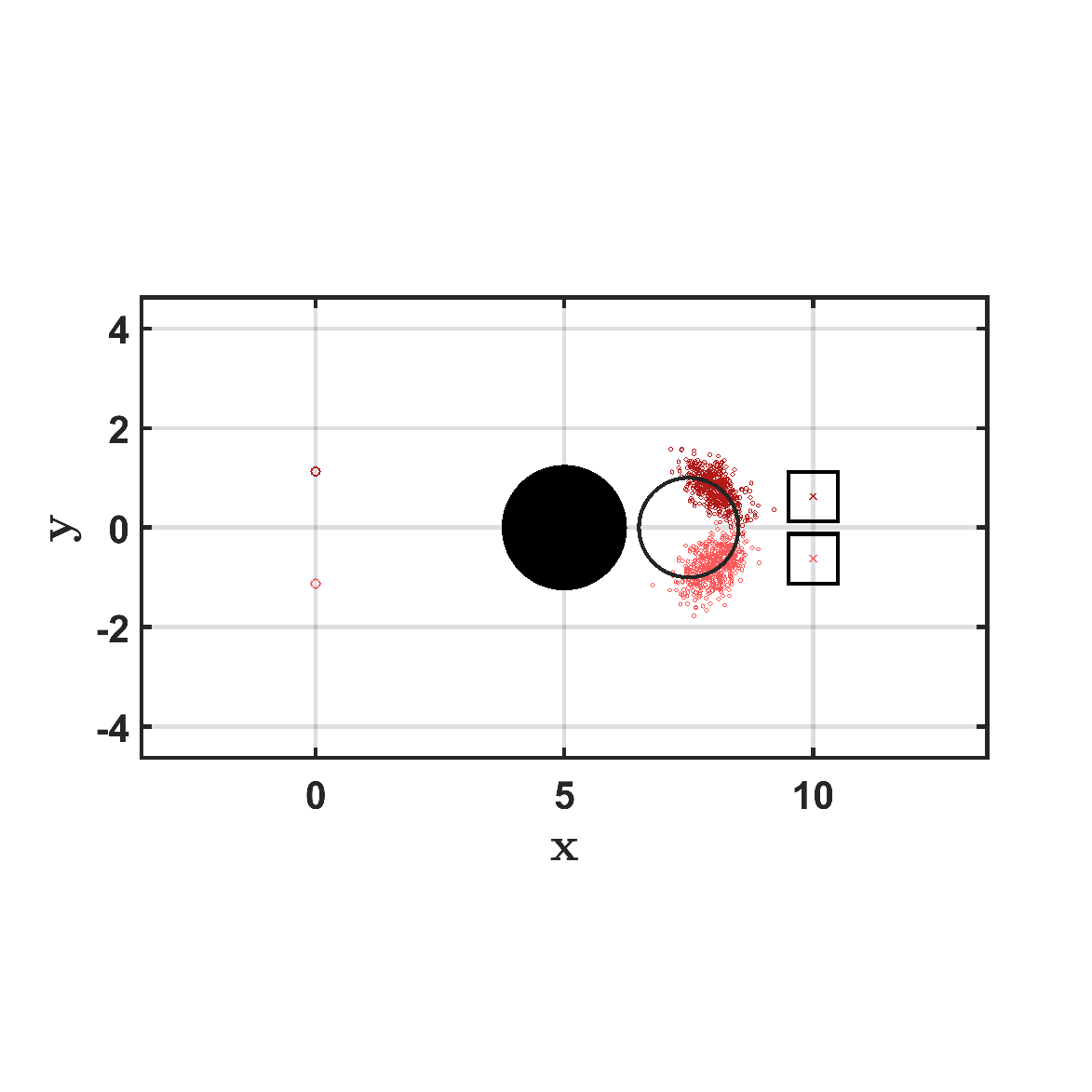}}
       \\
  \subfloat[ With chance constraints  \label{3c}]{
        \includegraphics[width=0.45\linewidth, trim={0.95cm 4cm 1.75cm 5cm},clip]{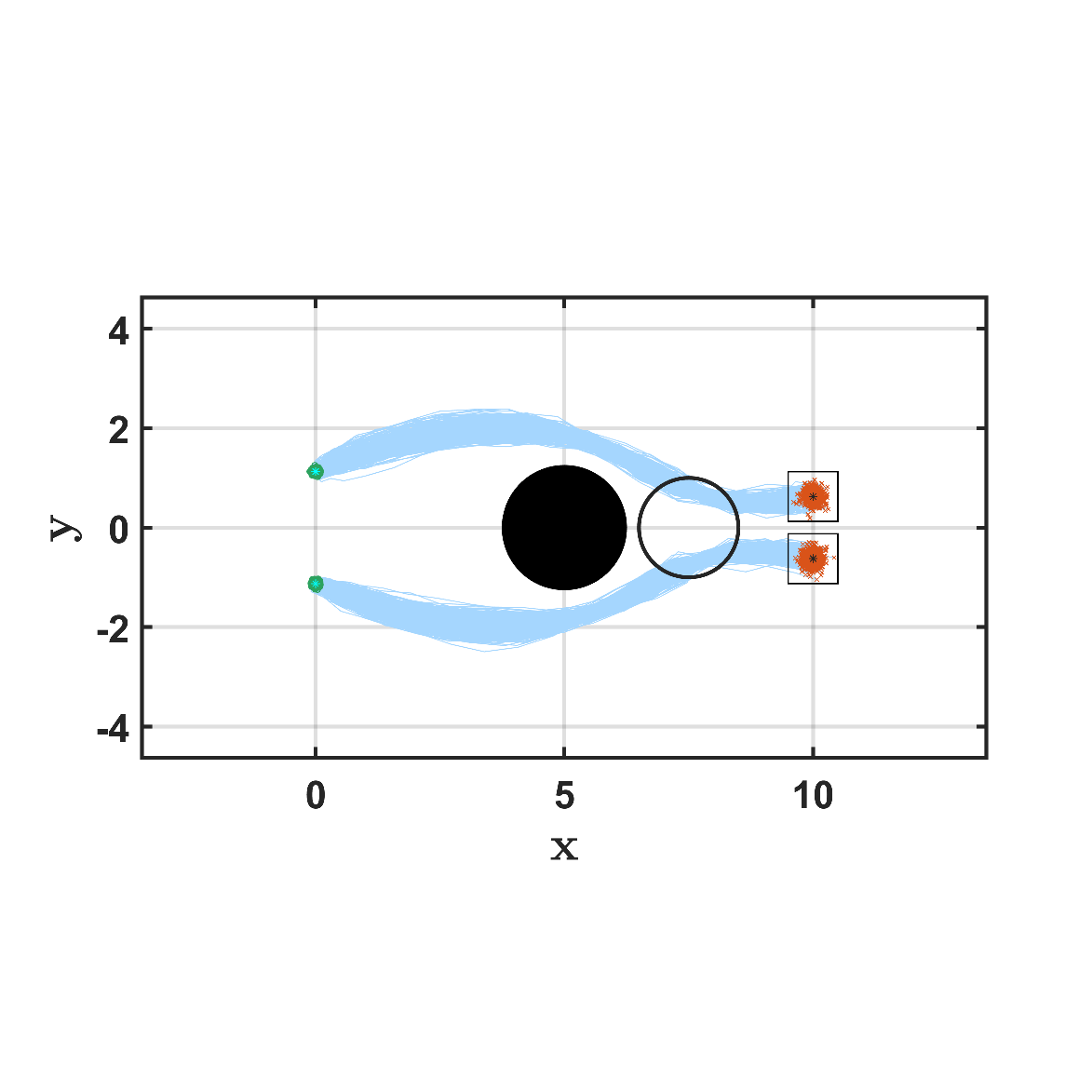}}
 \subfloat[ At $k=15$  \label{3d}]{
        \includegraphics[width=0.45\linewidth, trim={0.95cm 4cm 1.75cm 5cm},clip]{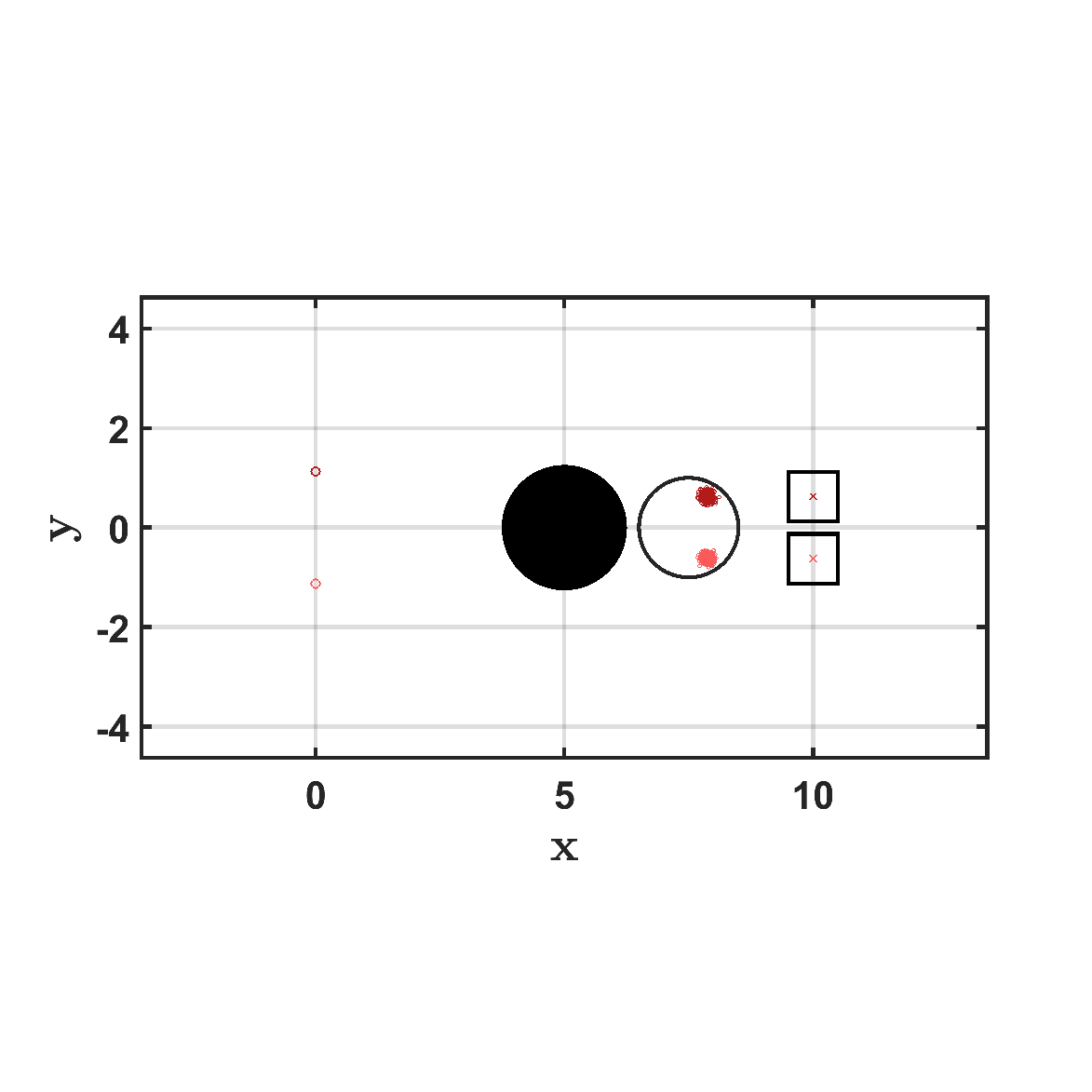}}
        \hfill
  \caption{ \textbf{Mixed Case Scenario-1:}  With terminal linear chance \eqref{Mixed Constraints: Robust linear Chance Constraints}, collision avoidance \eqref{Mixed Constraints: nonconvex norm chance}, \eqref{Mixed Constraints: inter-agent nonconvex chance}  and convex norm chance \eqref{Mixed Constraints: convex norm chance} constraints.}
  \label{fig3} 
  \vspace{-0.25cm}
\end{figure}
% % % % %
\begin{figure}[t!] 
    \centering
  \subfloat[ With constraints on $\Eb(\bx_i)$ \label{4a}]{
       \includegraphics[width= 0.425\linewidth, trim={0.95cm 4.5cm 1.75cm 5cm},clip]{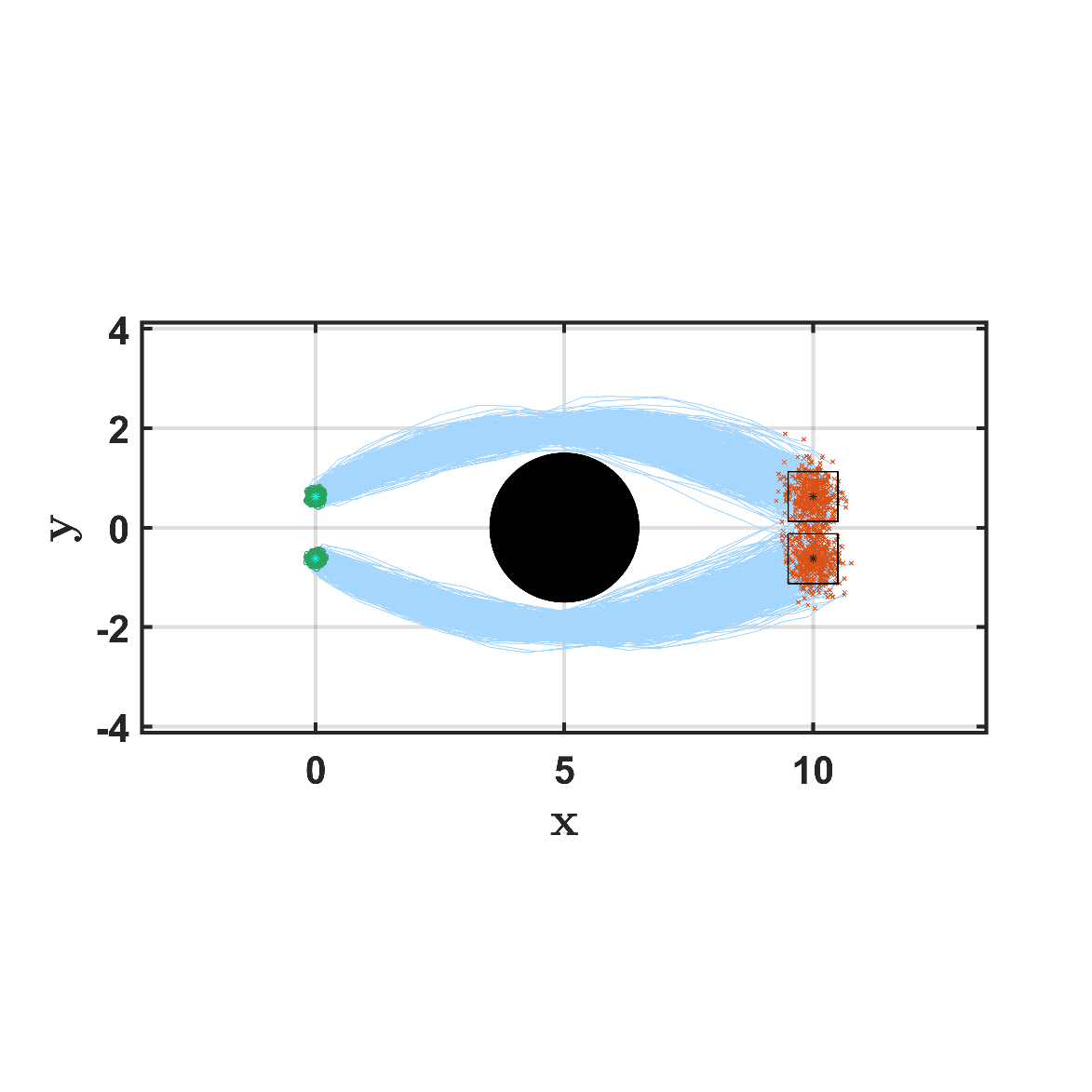}}
  \subfloat[ Adding terminal covariance constraints  \label{4b}]{
        \includegraphics[width=0.425\linewidth, trim={0.95cm 4.5cm 1.75cm 5cm},clip]{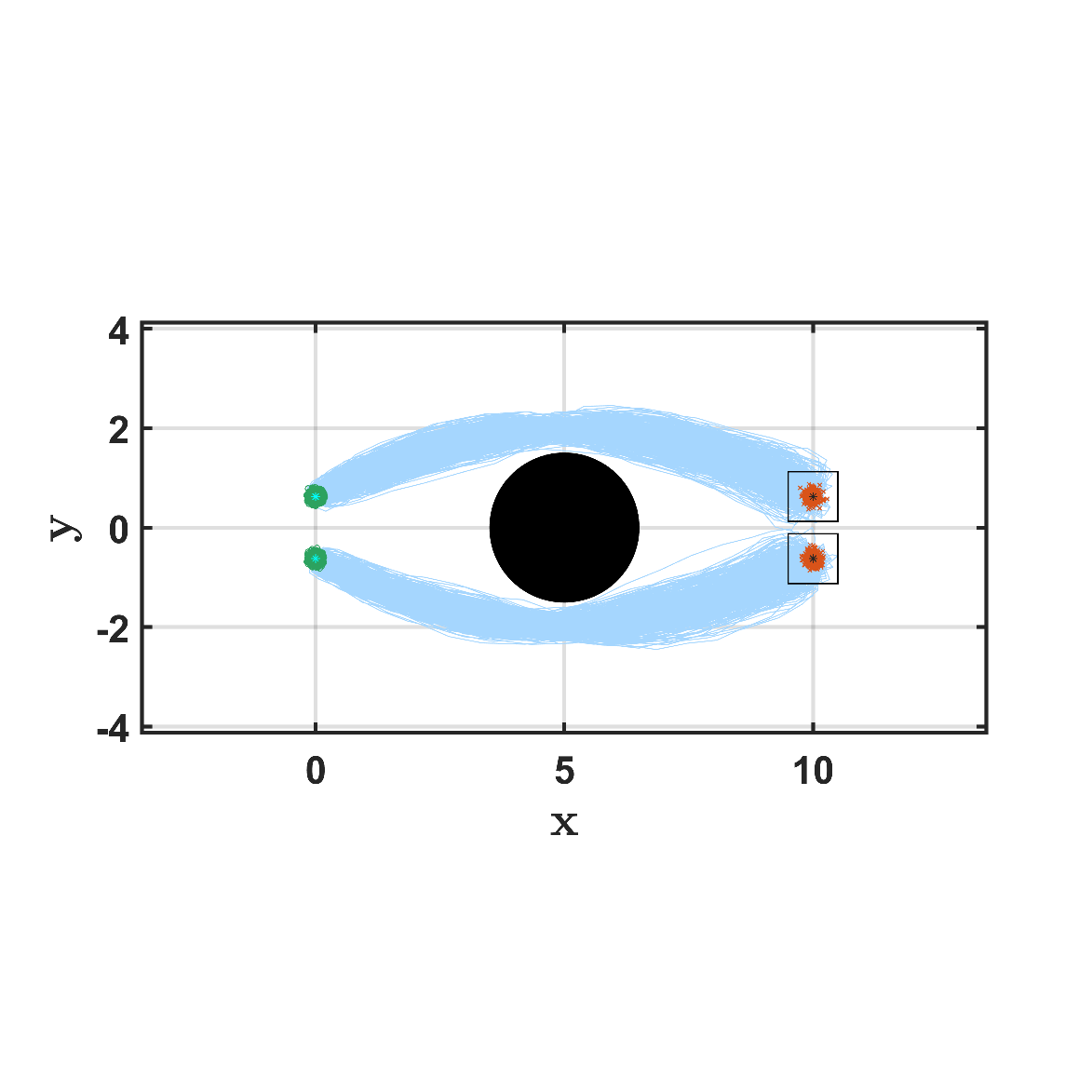}}
        \hfill
  \caption{ \textbf{Mixed Case Scenario-2:} With terminal linear \eqref{Robust Linear Constraints}, obstacle avoidance \eqref{Robust nonconvex norm-of-mean constraints}, \eqref{robust inter-agent collision avoidance}, and terminal covariance \eqref{Mixed Constraints: Covariance} constraints.}
  \label{fig4} 
  \vspace{-0.25cm}
\end{figure}

In Fig. \ref{fig3}, we demonstrate the effectiveness of robust chance constraints for the mixed case (defined in \ref{sec: robust constraints in mixed case}) by considering the same scenario as in Fig. \ref{fig2}. Fig. \ref{3a}- \ref{3b} show the trajectories of the agents when the control policies are obtained without accounting for the stochastic noise (i.e., applying only the deterministic robust constraints). 
% For each agent, the probability of obstacle avoidance constraint satisfaction is {\color{blue} enter}, convex norm constraint satisfaction is {\color{blue} enter}, and the terminal linear constraint satisfaction is {\color{blue} enter}. 
In Fig \ref{3c}-\ref{3d}, we address the effect of stochastic noise by considering the chance constraints. Since the derived constraint reformulations in Section \ref{sec: Extension to Mixed Disturbance case} are conservative, the probability of constraint constraint satisfaction is more than the set value.

In mixed disturbance cases where the deterministic uncertainty is relatively higher than the stochastic noise, the control gains obtained by considering robust constraints only on $\Eb(\bx_i)$ can also effectively control the deviation in the state due to stochastic noise. For instance, in Fig. \ref{fig4}, we showcase a mixed disturbance case wherein in Fig. \ref{4a} presents the trajectories when only the deterministic robust constraints are applied. It can be observed that most of the trajectory realizations for both the agents avoid the obstacle, with only the terminal constraints being violated. In Fig. \ref{4b}, we add the terminal covariance constraint to address the said issue. 
\subsubsection{Handling Complex Scenarios}
\begin{figure}[t!] 
    \centering
  \subfloat[ Robust Trajectory \label{5a}]{
       \includegraphics[width= 0.575\linewidth, trim={0.8cm 3.75cm 1.75cm 5cm},clip]{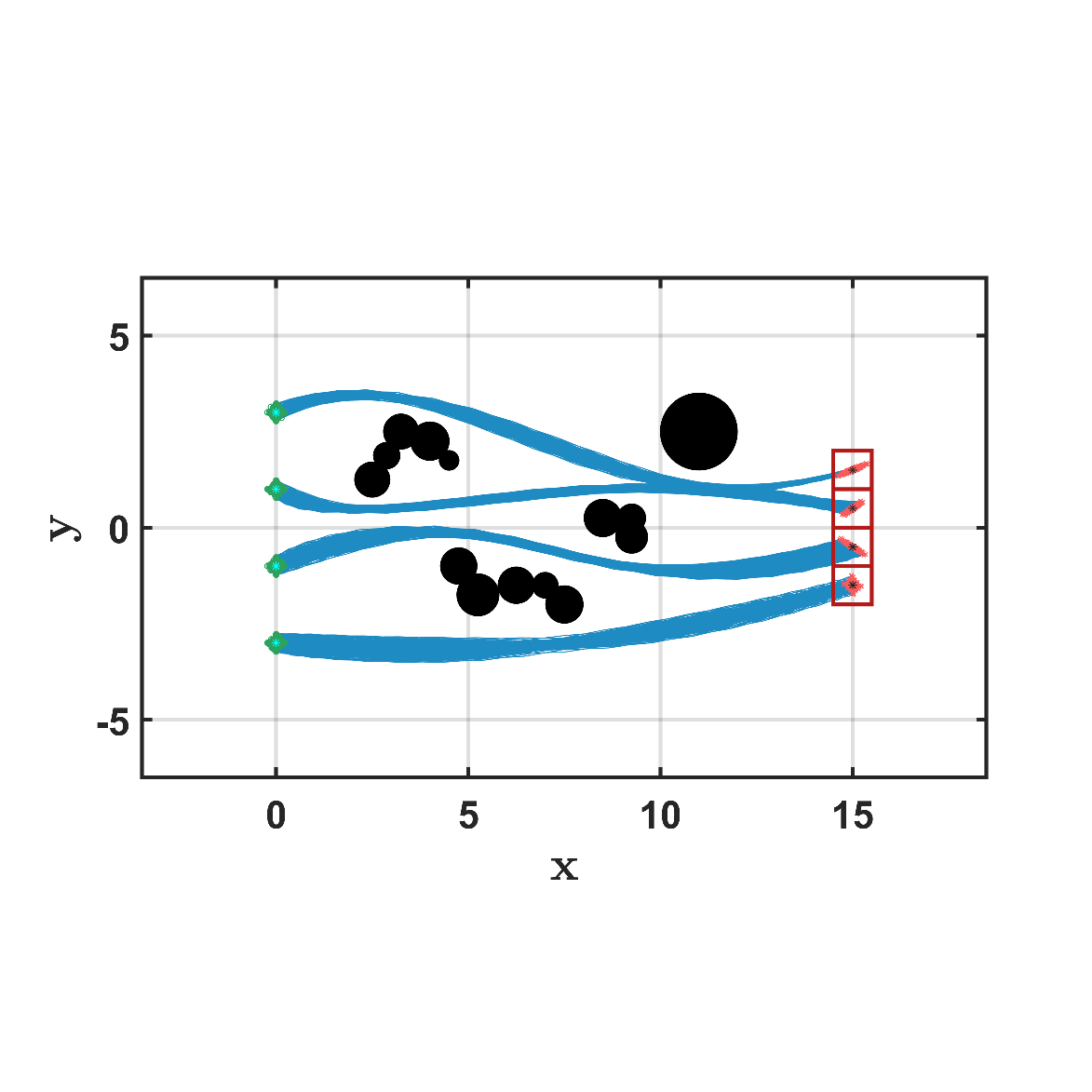}}
  \subfloat[Distance Plot  \label{5b}]{
        \includegraphics[width=0.35\linewidth, trim={0cm -0.35cm 0cm 0cm},clip]{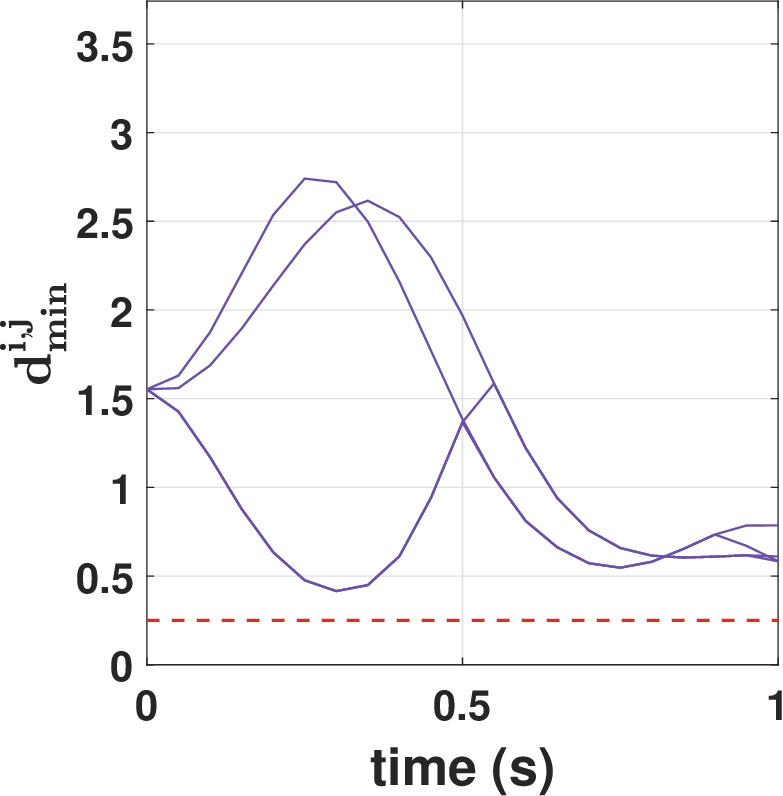}} \\
 \subfloat[ $k = 10$ \label{5c}]{
       \includegraphics[width= 0.45\linewidth, trim={0.8cm 3.75cm 1.75cm 5cm},clip]{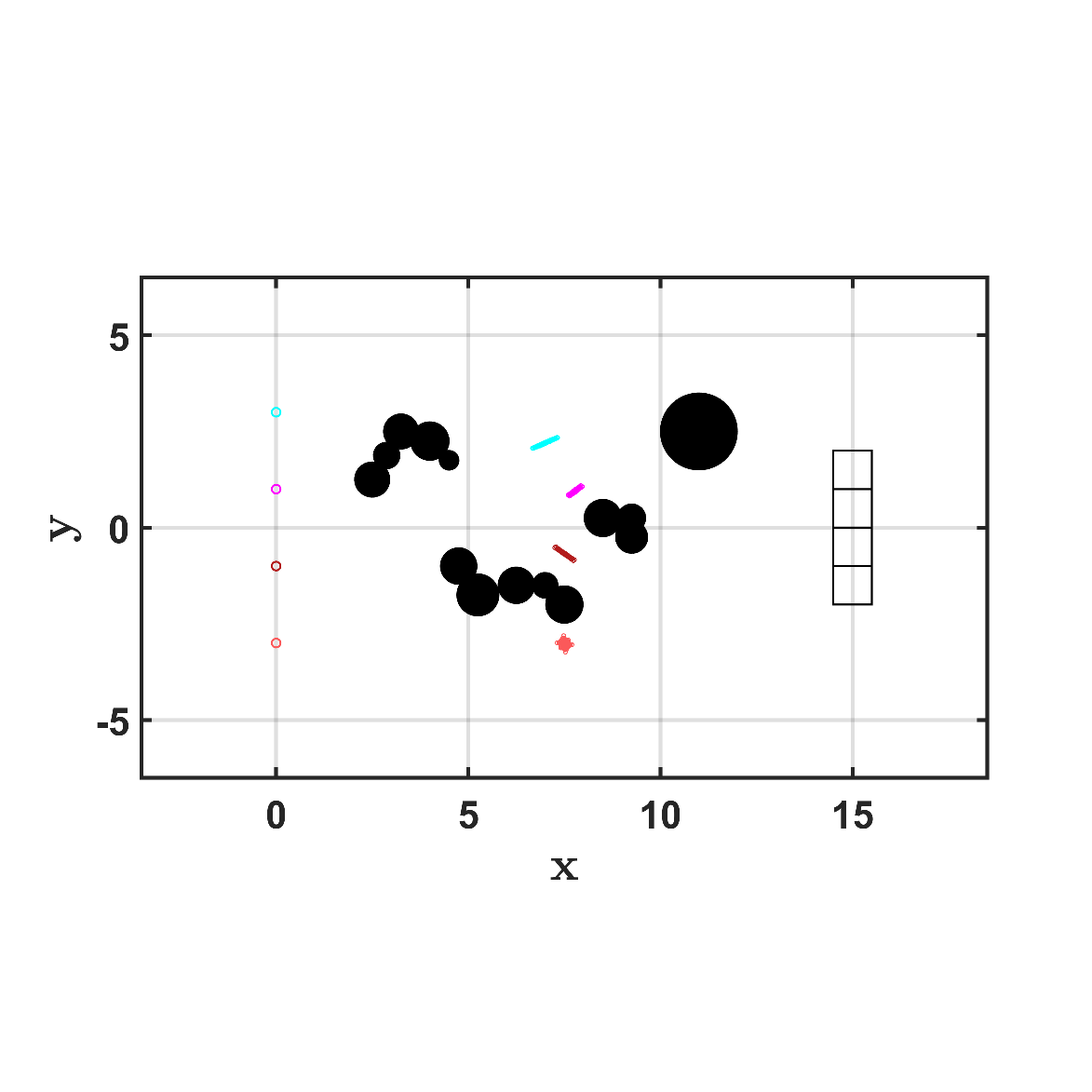}} 
 \subfloat[ $k = 20$ \label{5d}]{
       \includegraphics[width= 0.45\linewidth, trim={0.8cm 3.75cm 1.75cm 5cm},clip]{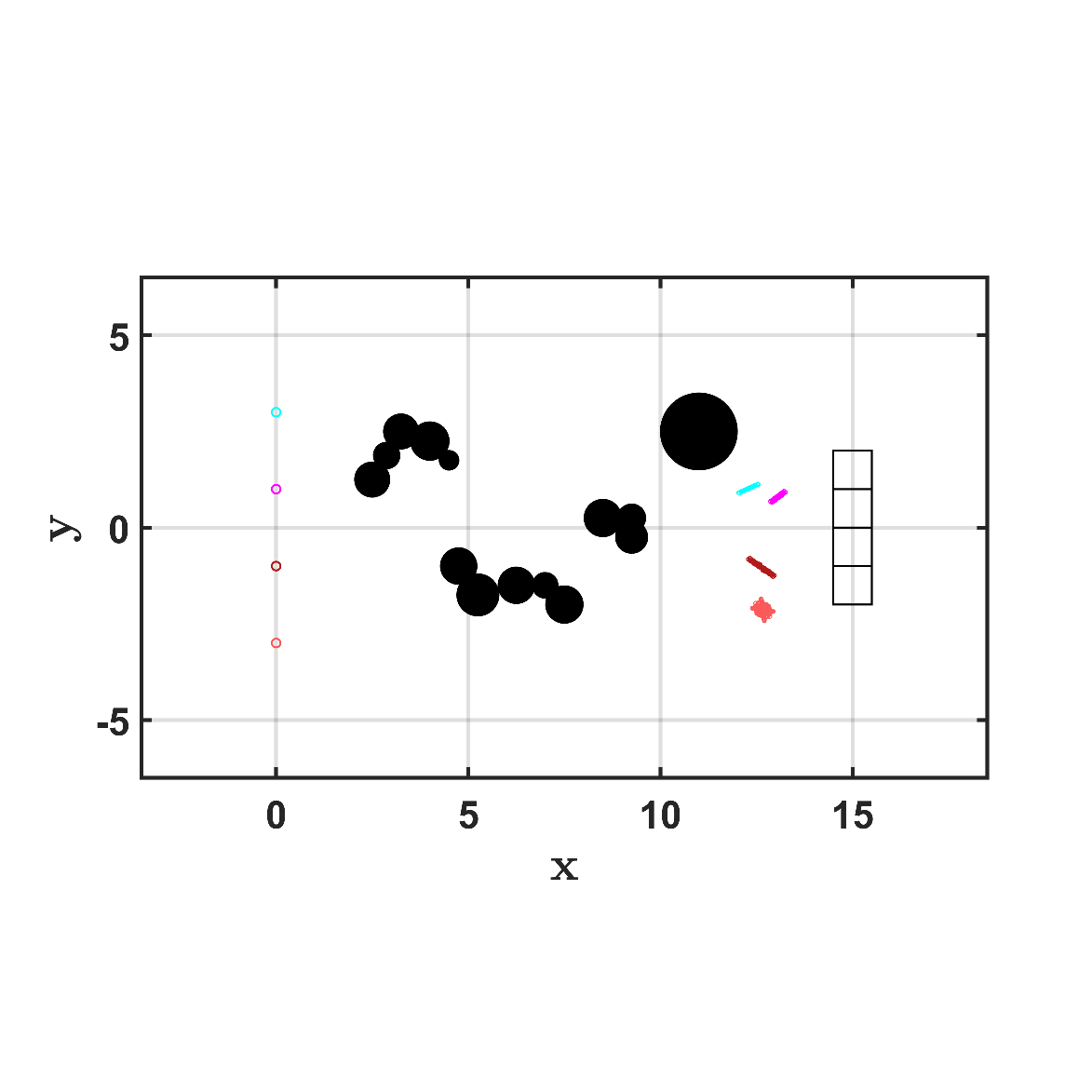}}
        \hfill
  \caption{ \textbf{Deterministic Case with Non-convex Obstacles:} 4-agent system with 14 circular obstacles. Terminal linear \eqref{Robust Linear Constraints} and collision avoidance \eqref{Robust nonconvex norm-of-mean constraints}, \eqref{robust inter-agent collision avoidance} constraints are imposed.}
  \label{fig5} 
  \vspace{-0.25cm}
\end{figure}

\begin{figure}[t!] 
    \centering
    \subfloat[ Robust Trajectory \label{6a}]{
       \includegraphics[width= 0.425\linewidth, trim={1cm 0cm 1.75cm 1cm},clip]{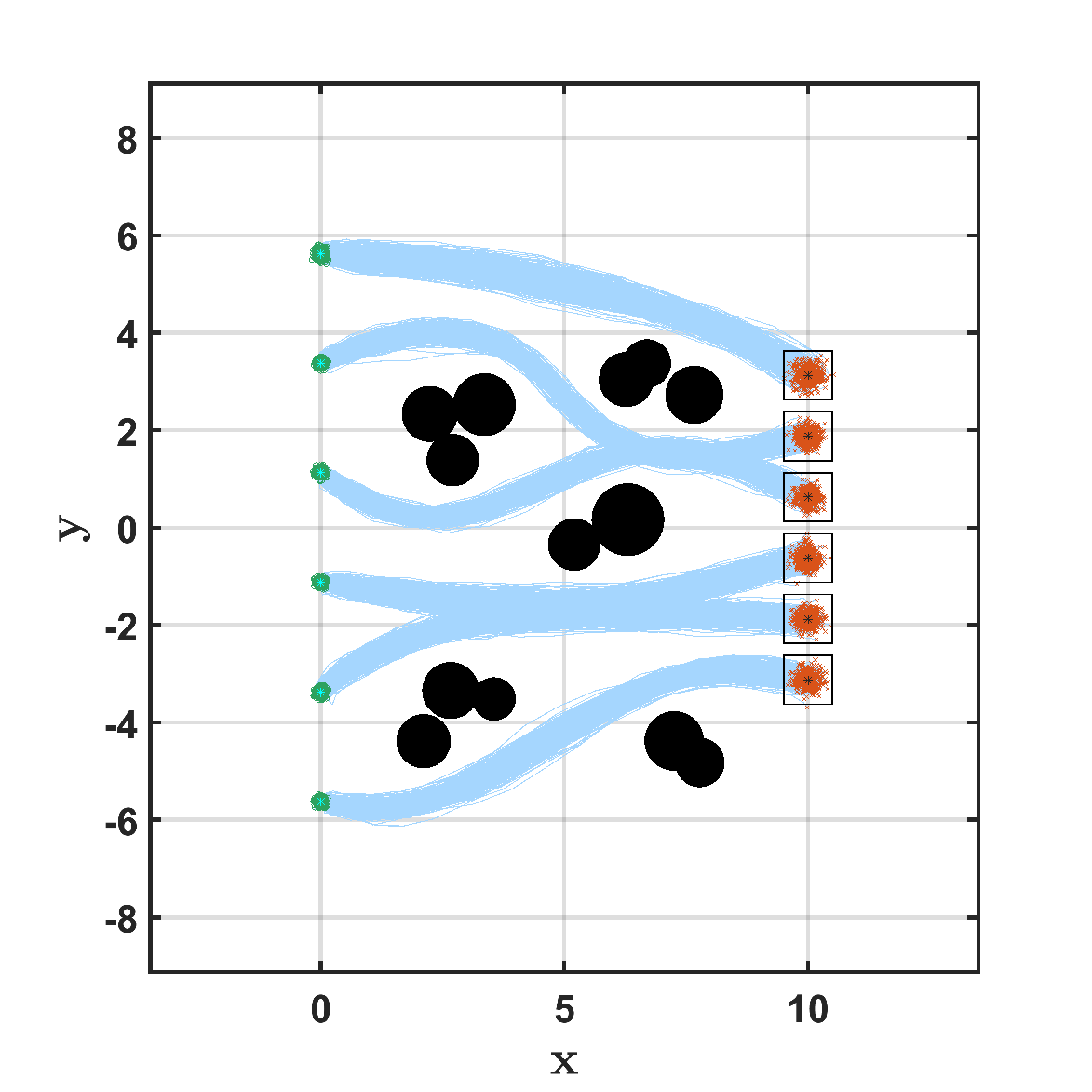}}
    \subfloat[ $k=5$ \label{6b}]{
       \includegraphics[width= 0.425\linewidth, trim={5.5cm 0cm 6.5cm 1cm},clip]{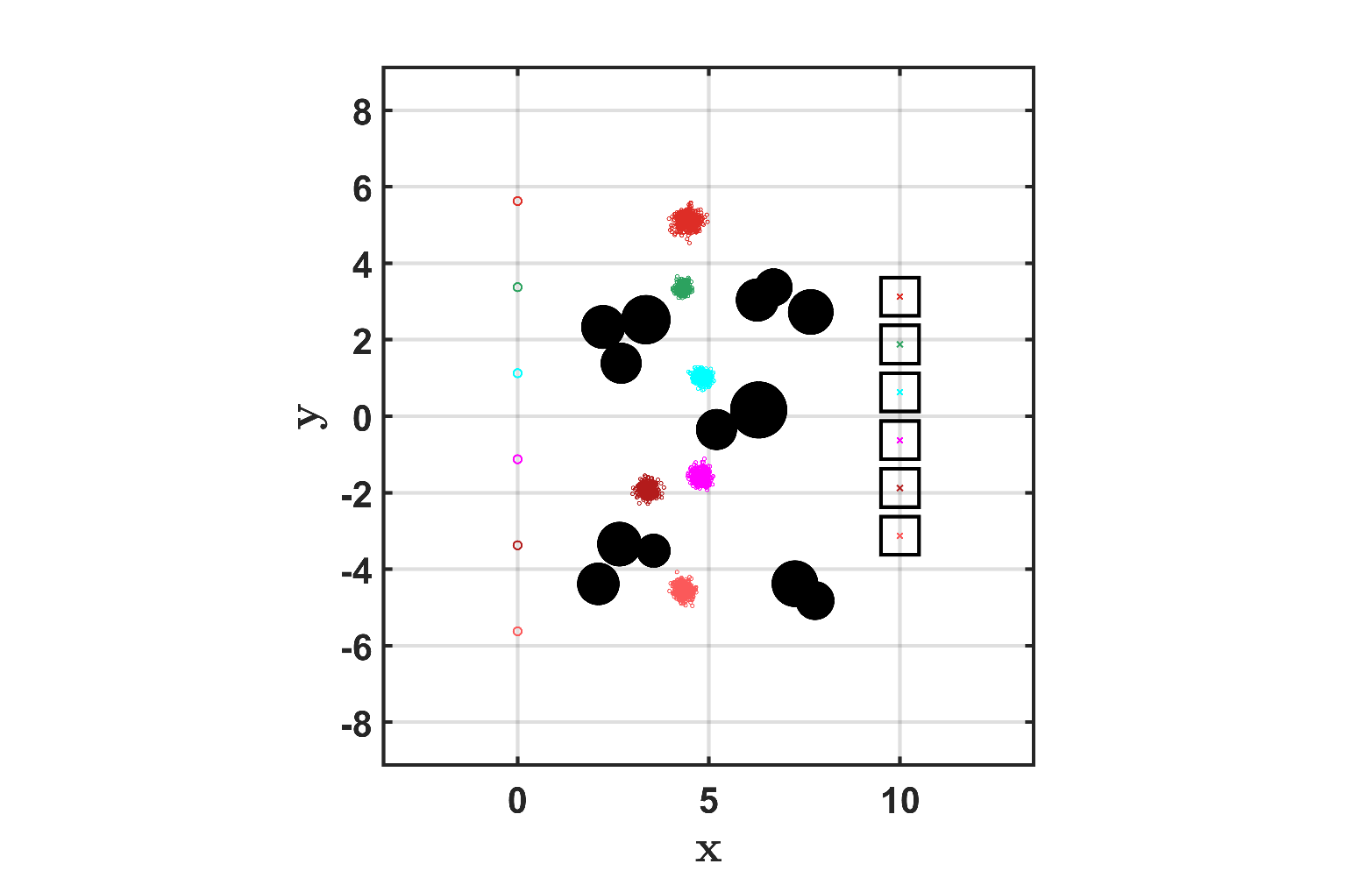}}
       \\
    \subfloat[ $k=10$ \label{6c}]{
       \includegraphics[width= 0.425\linewidth, trim={5.5cm 0cm 6.5cm 1cm},clip]{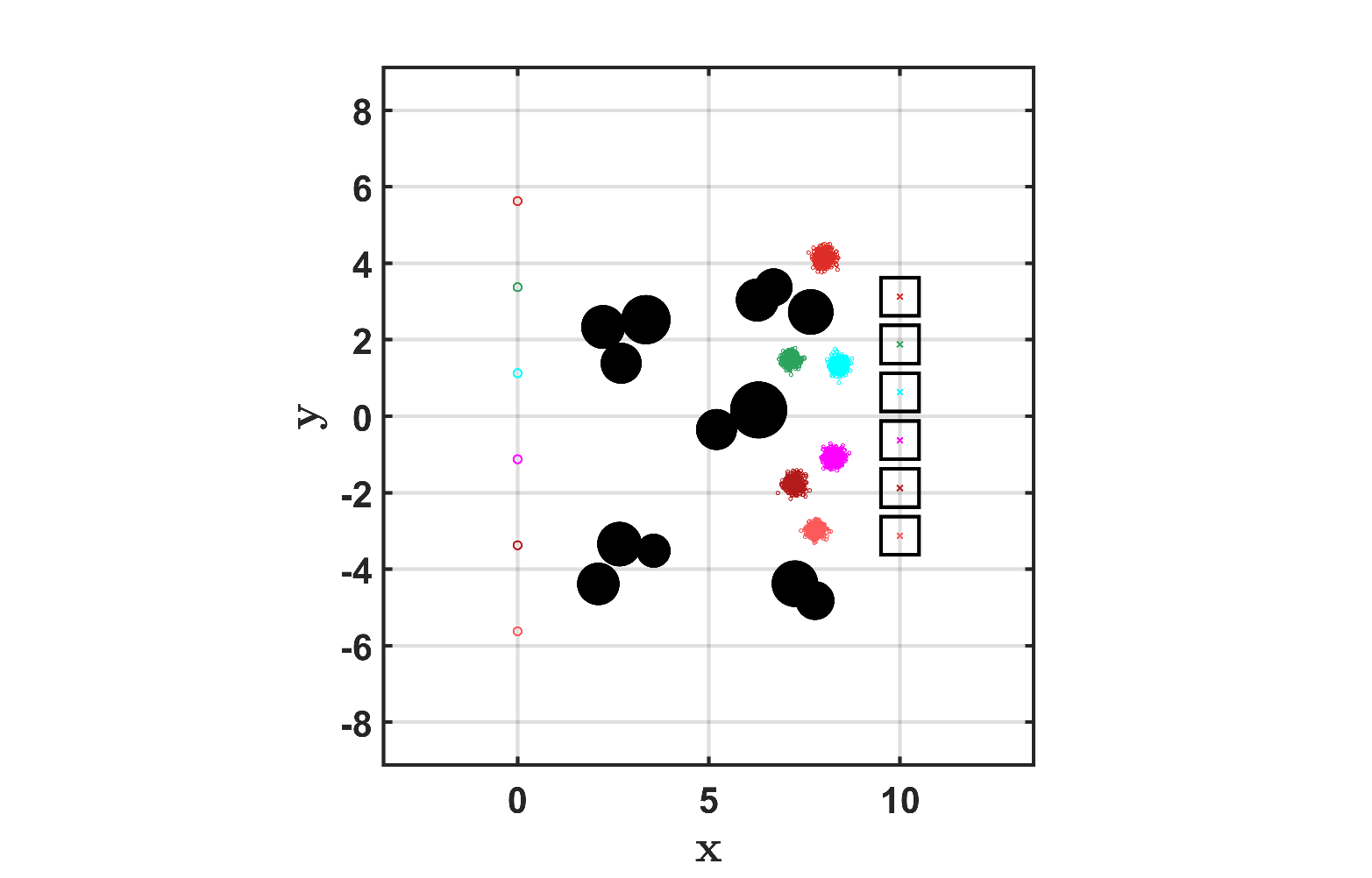}}
  \subfloat[Distance Plot  \label{6d}]{
        \includegraphics[width=0.45\linewidth, trim={0cm -0.28cm 0cm 0cm},clip]{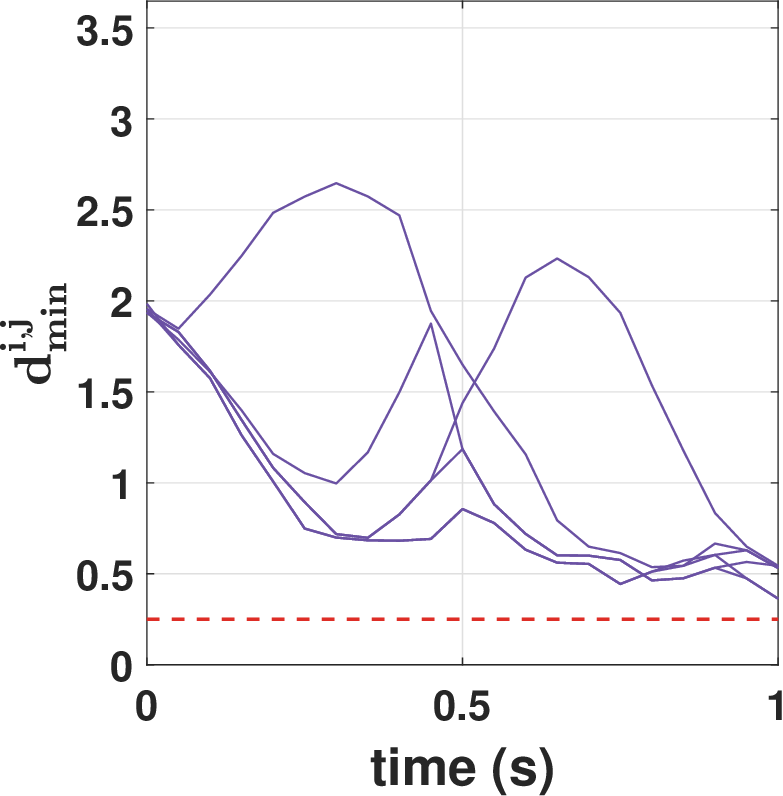}}
        \hfill
  \caption{ \textbf{Mixed Case with Non-convex Obstacles:} 6-agent system with 13 circular obstacles. Terminal linear chance \eqref{Mixed Constraints: Robust linear Chance Constraints} and collision avoidance \eqref{Mixed Constraints: nonconvex norm chance}, \eqref{Mixed Constraints: inter-agent nonconvex chance} constraints are imposed.} 
  \label{fig6} 
\end{figure}

In Fig. \ref{fig5}, we showcase a scenario with nonconvex obstacles in a deterministic case. From Fig. \ref{5a}, it is evident that all the agents successfully avoid collision with the obstacles. It can also be observed that though the terminal target bounds of the agents are intersecting, the terminal position of the agents do not collide. Further, a distance plot shown in Fig. \ref{5b} proves that there are no inter-agent collisions. The term $d^{i,j}_{\text{min}}$ in the plot represents the minimum distance between the trajectory realizations of the neighboring agents and the dashed red line represents the inter-agent collision threshold. Further, we also provide the snapshots of the trajectory realizations at time steps $k=10$ and $k=15$ in Fig. \ref{5c} and \ref{5d} respectively. Likewise, in Fig. \ref{fig6}, we illustrate the performance of robust chance constraints with non-convex obstacles in mixed case. We observe that the agents successfully avoid collision with the obstacles and the other agents. 
% Further, Fig. \ref{6b} and \ref{6c} present the snapshot of the agents' trajectories at timestep $k=5$ and $k=10$ respectively. Fig. \ref{6b}-\ref{6d} illustrate that there are no inter-agent collisions. 
% % %

\begin{figure*}[t!] 
    \centering
  \subfloat[ Robust Trajectory \label{7a}]{
       \includegraphics[width= 0.245\linewidth, trim={0.5cm 0.5cm 1cm 1.75cm},clip]{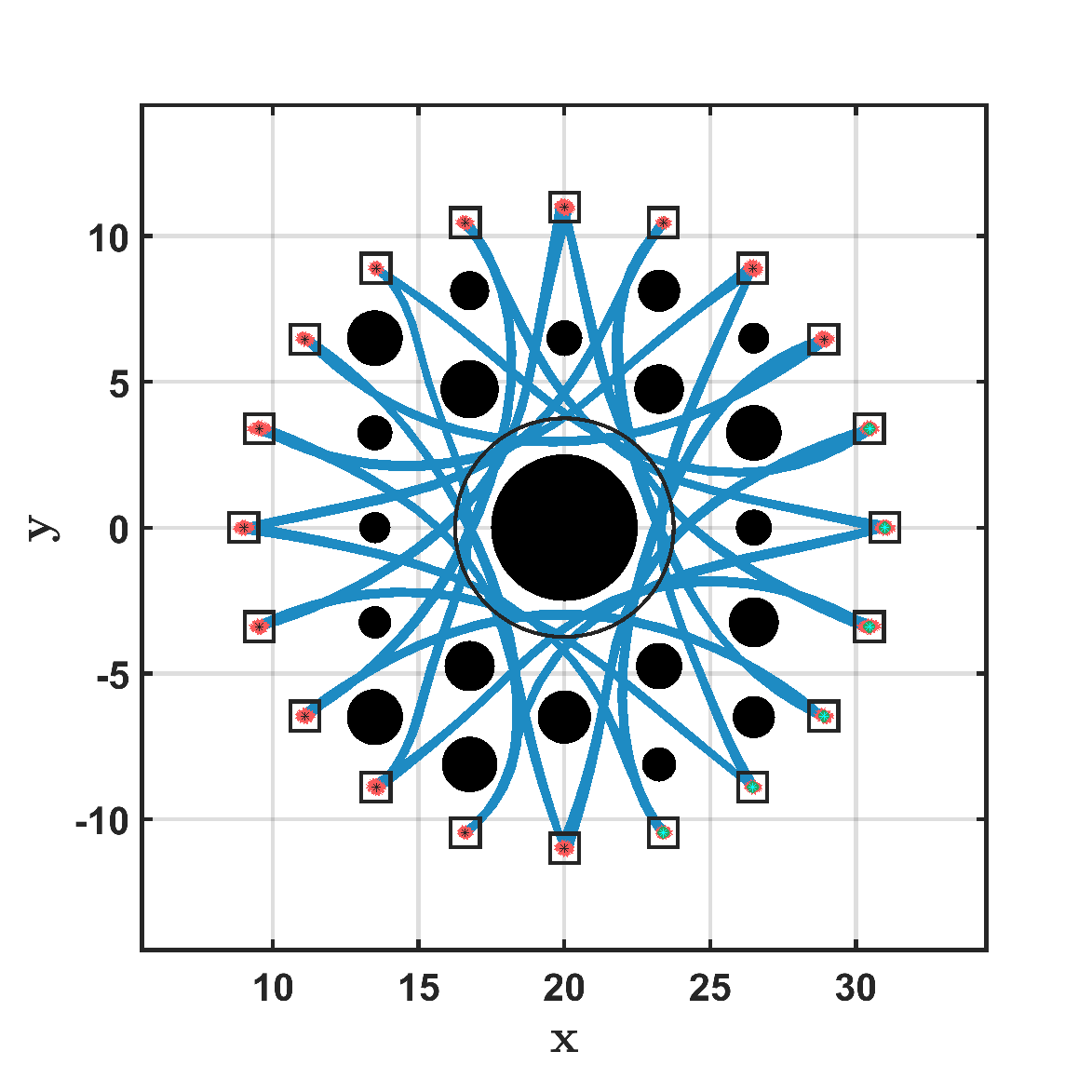}}
  \subfloat[ At $k=10$ \label{7b}]{
       \includegraphics[width= 0.245\linewidth, trim={0.5cm 0.5cm 1cm 1.75cm},clip]{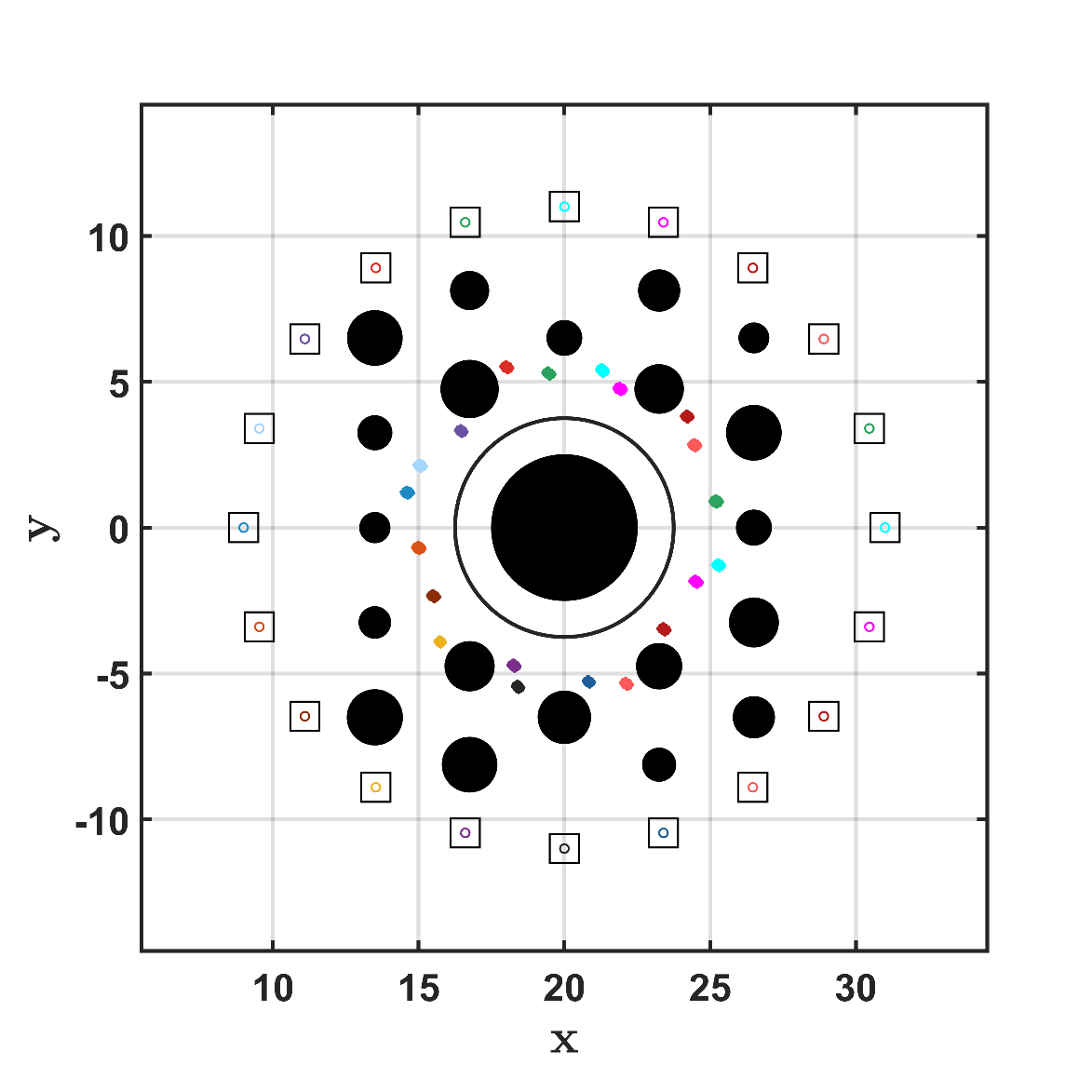}}
  % \\     
  \subfloat[ At $k=15$ \label{7c}]{
       \includegraphics[width= 0.245\linewidth, trim={0.5cm 0.5cm 1cm 1.75cm},clip]{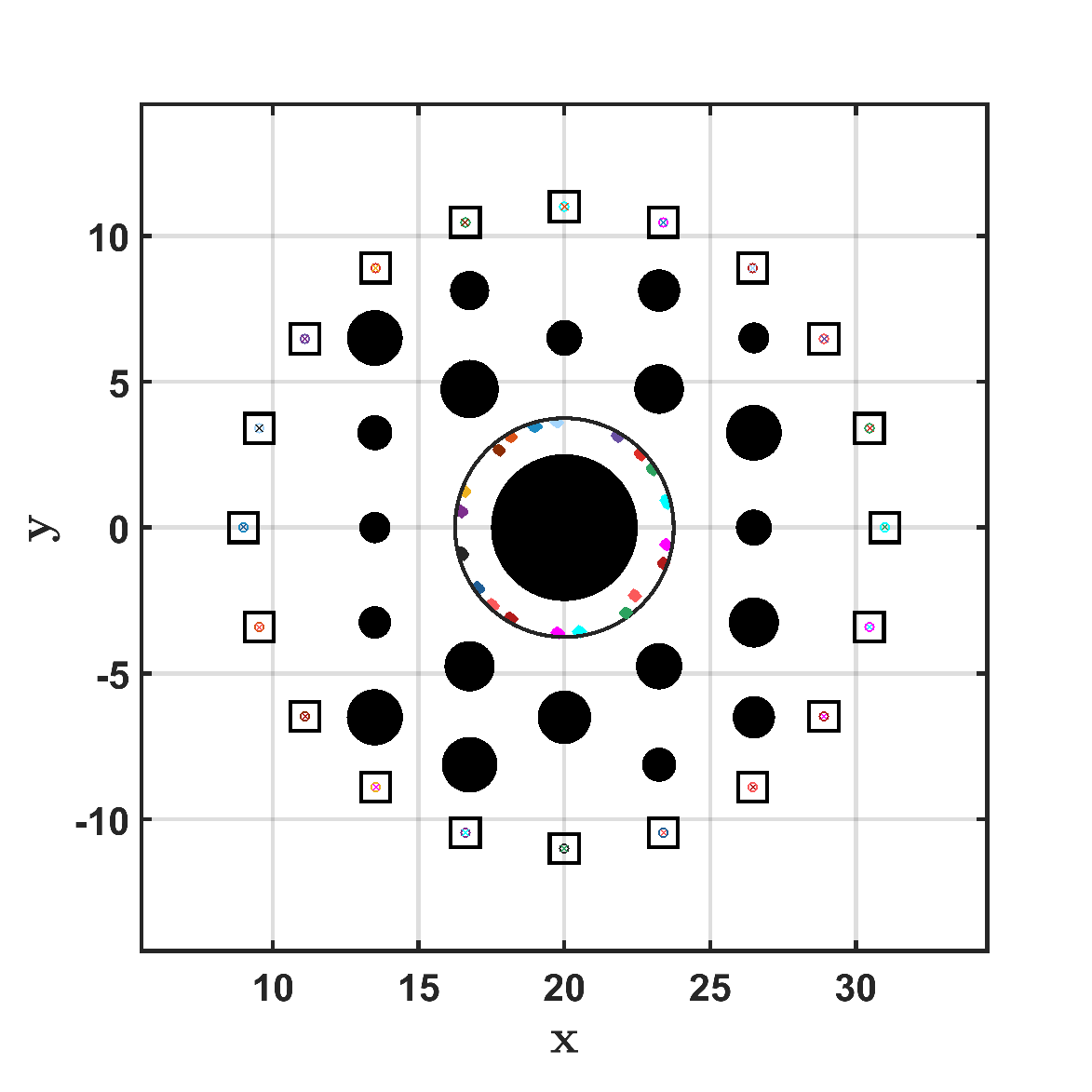}}
  \subfloat[Distance Plot  \label{7d}]{
        \includegraphics[width=0.225\linewidth, trim={0cm 0cm 0cm 0cm},clip]{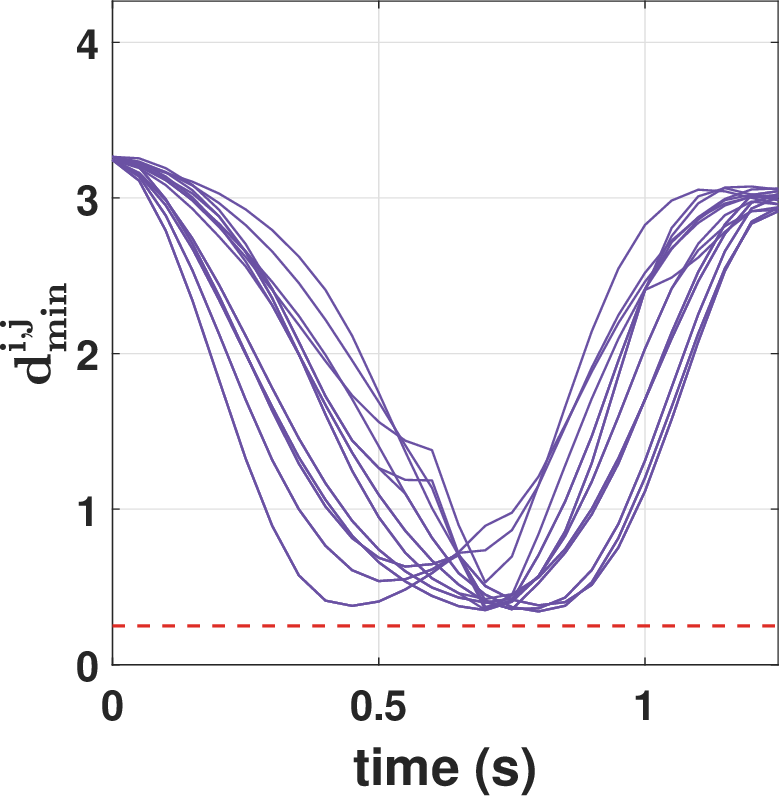}}
        \hfill
  \caption{ \textbf{Deterministic Case with 20 agents and 21 obstacles:} Terminal linear \eqref{Robust Linear Constraints}, collision avoidance \eqref{Robust nonconvex norm-of-mean constraints} \eqref{robust inter-agent collision avoidance} and convex norm \eqref{Robust convex norm-of-mean constraints} constraints are imposed.}
  \label{fig7} 
  \vspace{-0.25cm}
\end{figure*}
\begin{figure*}[t!] 
    \centering
  \subfloat[ Robust Trajectory \label{8a}]{
       \includegraphics[width= 0.245\linewidth, trim={1.5cm 0.25cm 2.5cm 1cm},clip]{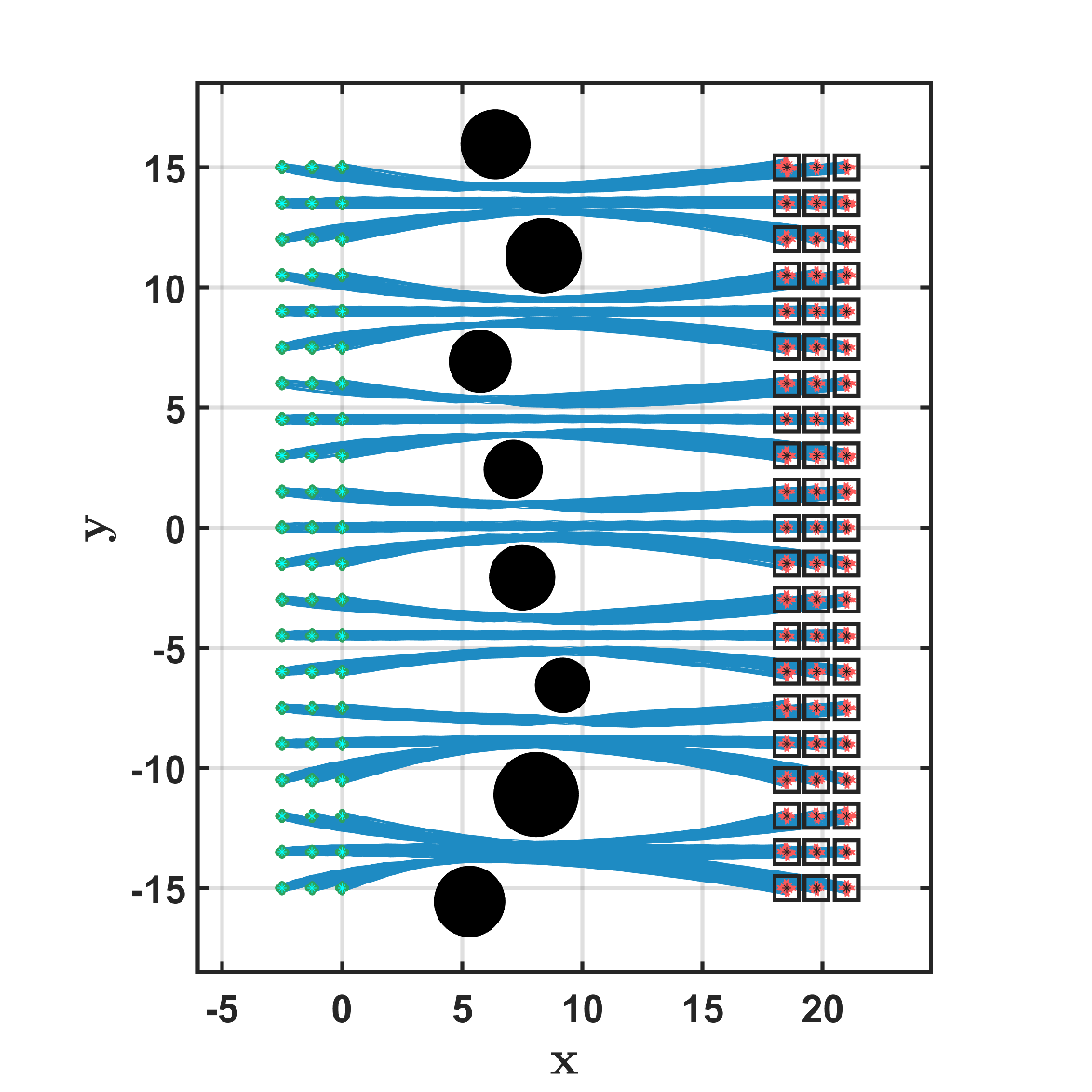}}
  \subfloat[ At $k=10$ \label{8b}]{
       \includegraphics[width= 0.245\linewidth, trim={1.5cm 0.25cm 2.5cm 1cm},clip]{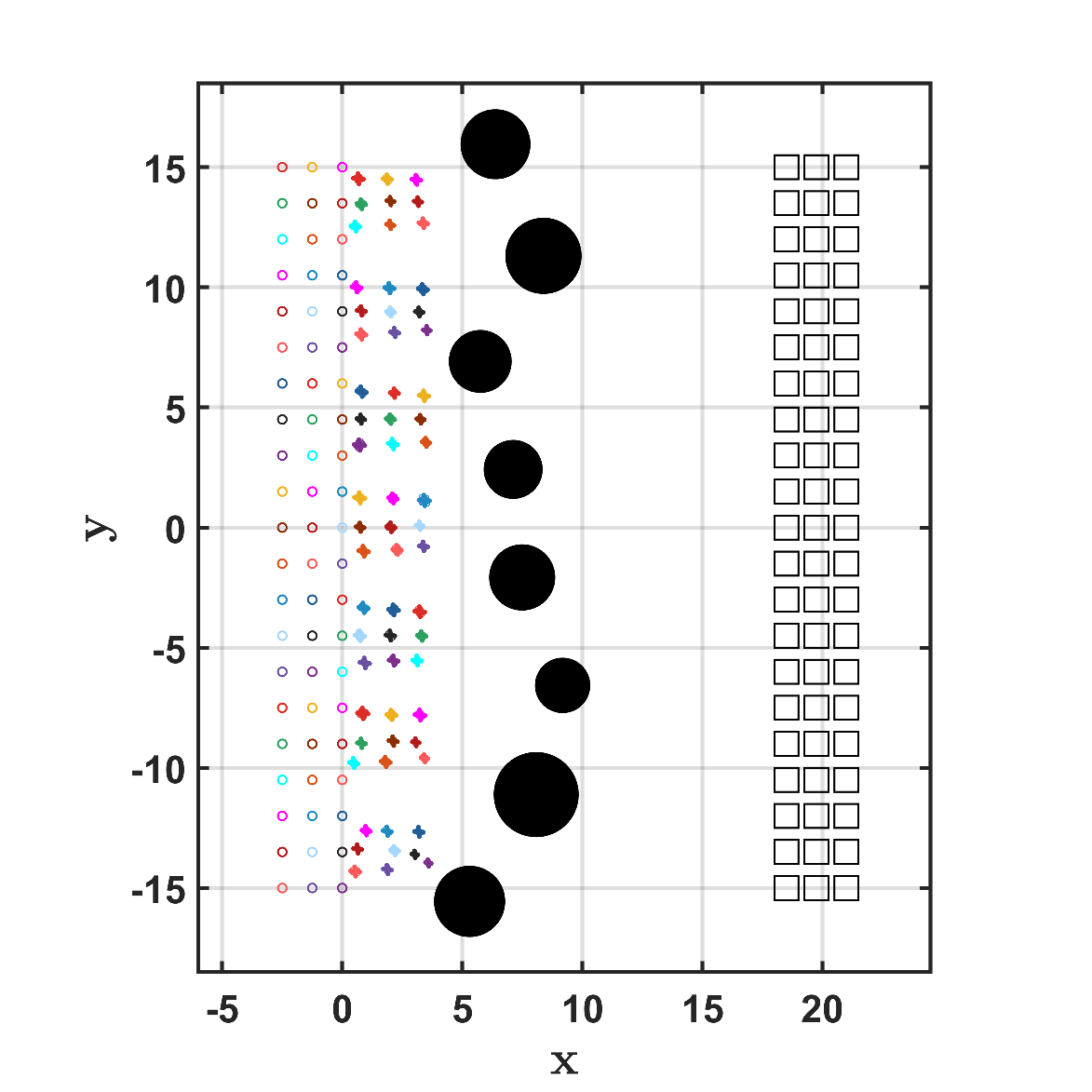}}
  % \\     
  \subfloat[ At $k=15$ \label{8c}]{
       \includegraphics[width= 0.245\linewidth, trim={1.5cm 0.25cm 2.5cm 1cm},clip]{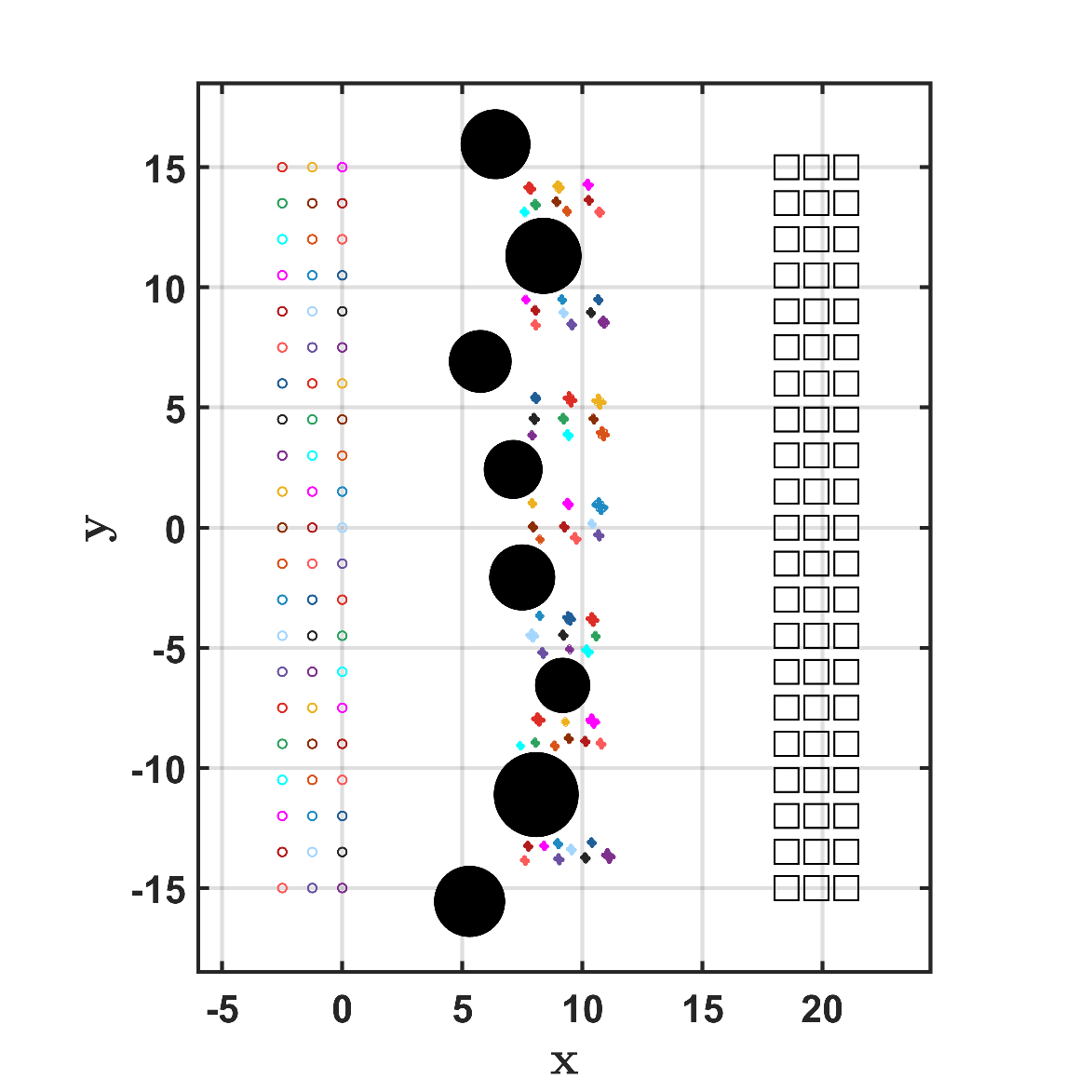}}
  \subfloat[Distance Plot  \label{8d}]{
        \includegraphics[width=0.225\linewidth, trim={0cm -2.5cm 0cm 0cm},clip]{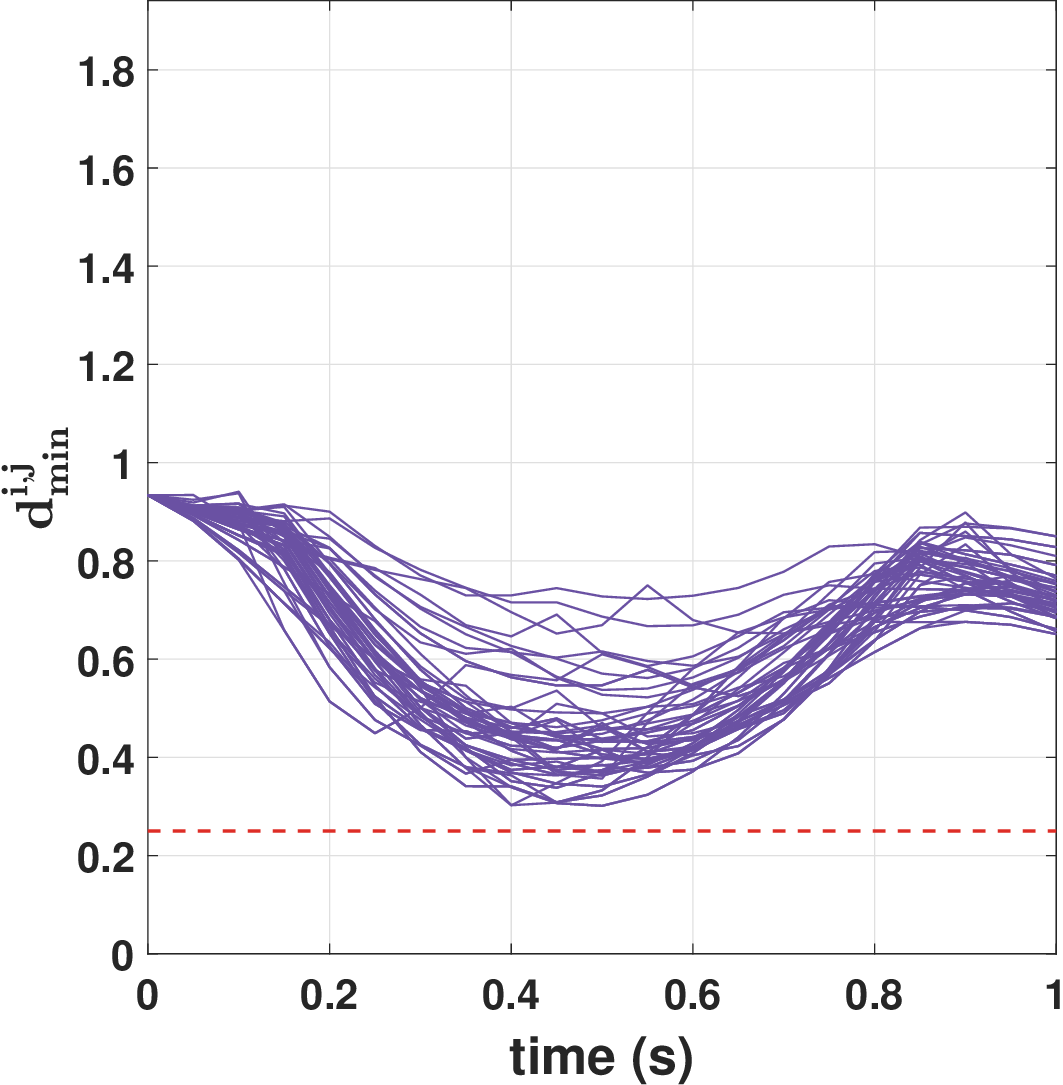}}
        \hfill
  \caption{ \textbf{Deterministic Case with 63 agents and 8 obstacles:} Terminal linear \eqref{Robust Linear Constraints}, collision avoidance \eqref{Robust nonconvex norm-of-mean constraints} \eqref{robust inter-agent collision avoidance} and convex norm \eqref{Robust convex norm-of-mean constraints} constraints are imposed.}
  \label{fig8} 
  \vspace{-0.25cm}
\end{figure*}
\begin{figure*}[t!] 
    \centering
  \subfloat[ Robust Trajectory \label{9a}]{
       \includegraphics[width= 0.25\linewidth, trim={0.5cm 0.5cm 1cm 1.75cm},clip]{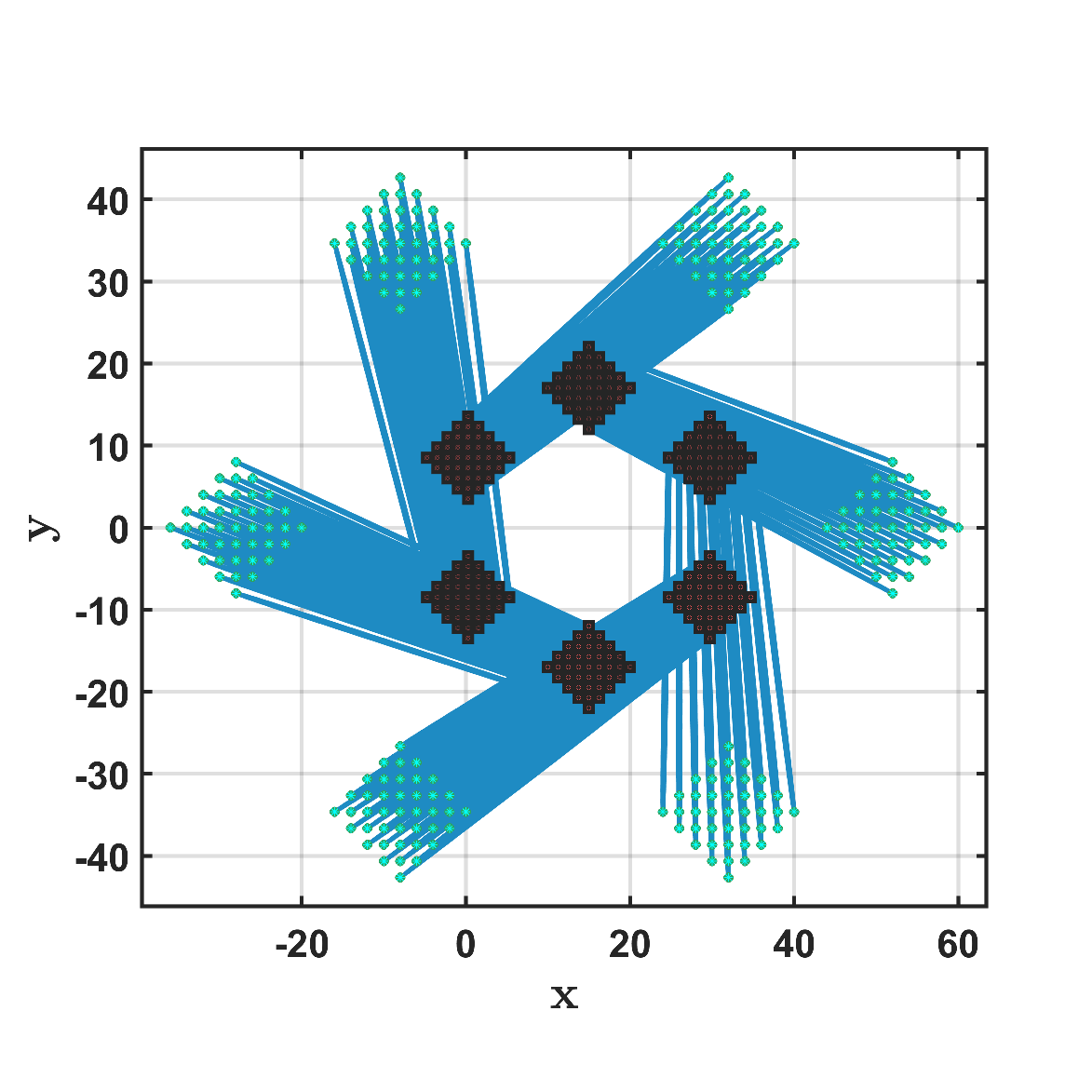}}
    \subfloat[ At $k=15$ \label{9b}]{
       \includegraphics[width= 0.25\linewidth, trim={0.5cm 0.5cm 1cm 1.75cm},clip]{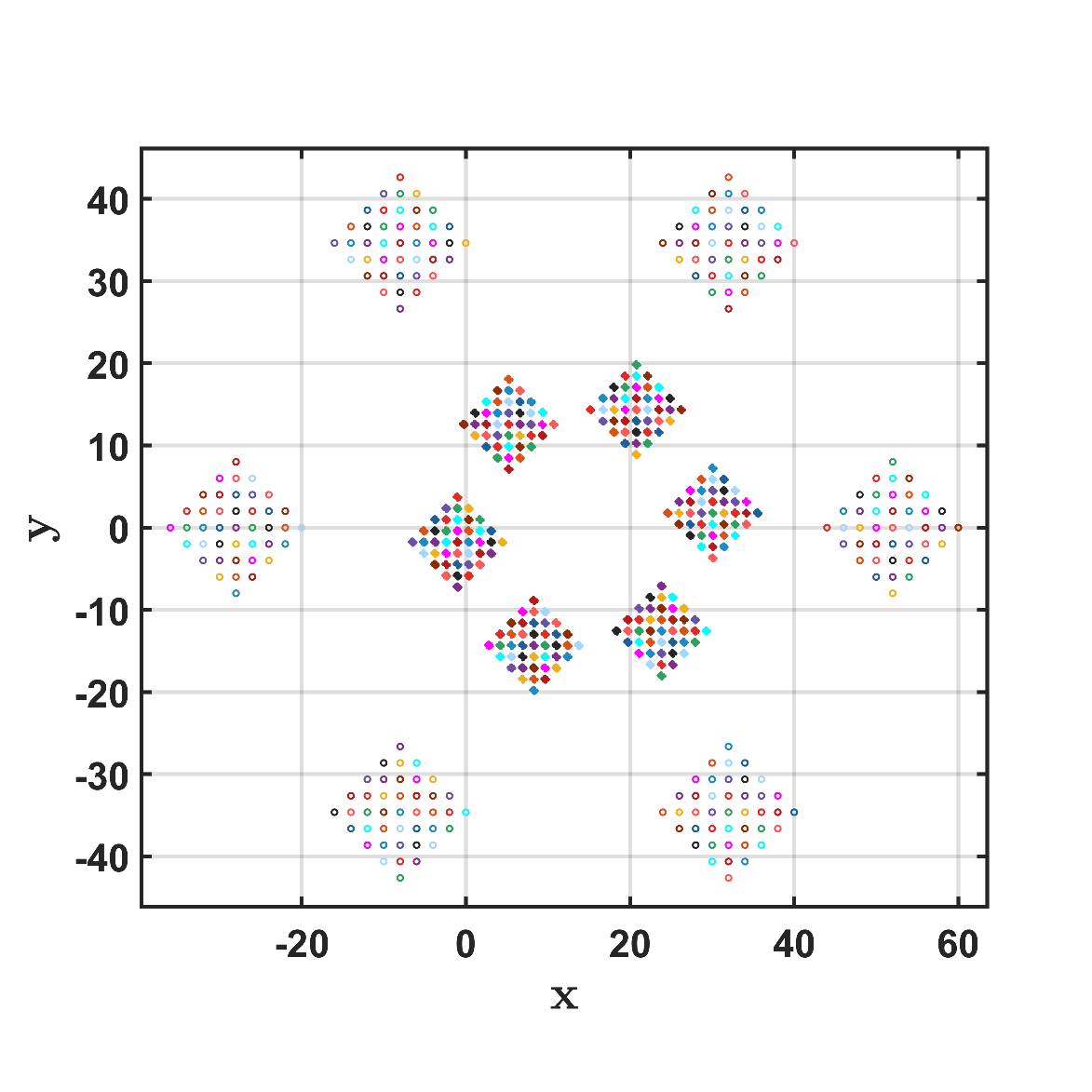}}
    % \\
    \subfloat[ At $k=20$ \label{9c}]{
        \includegraphics[width=0.215\linewidth, trim={0.5cm -0.5cm 1.5cm 0cm},clip]{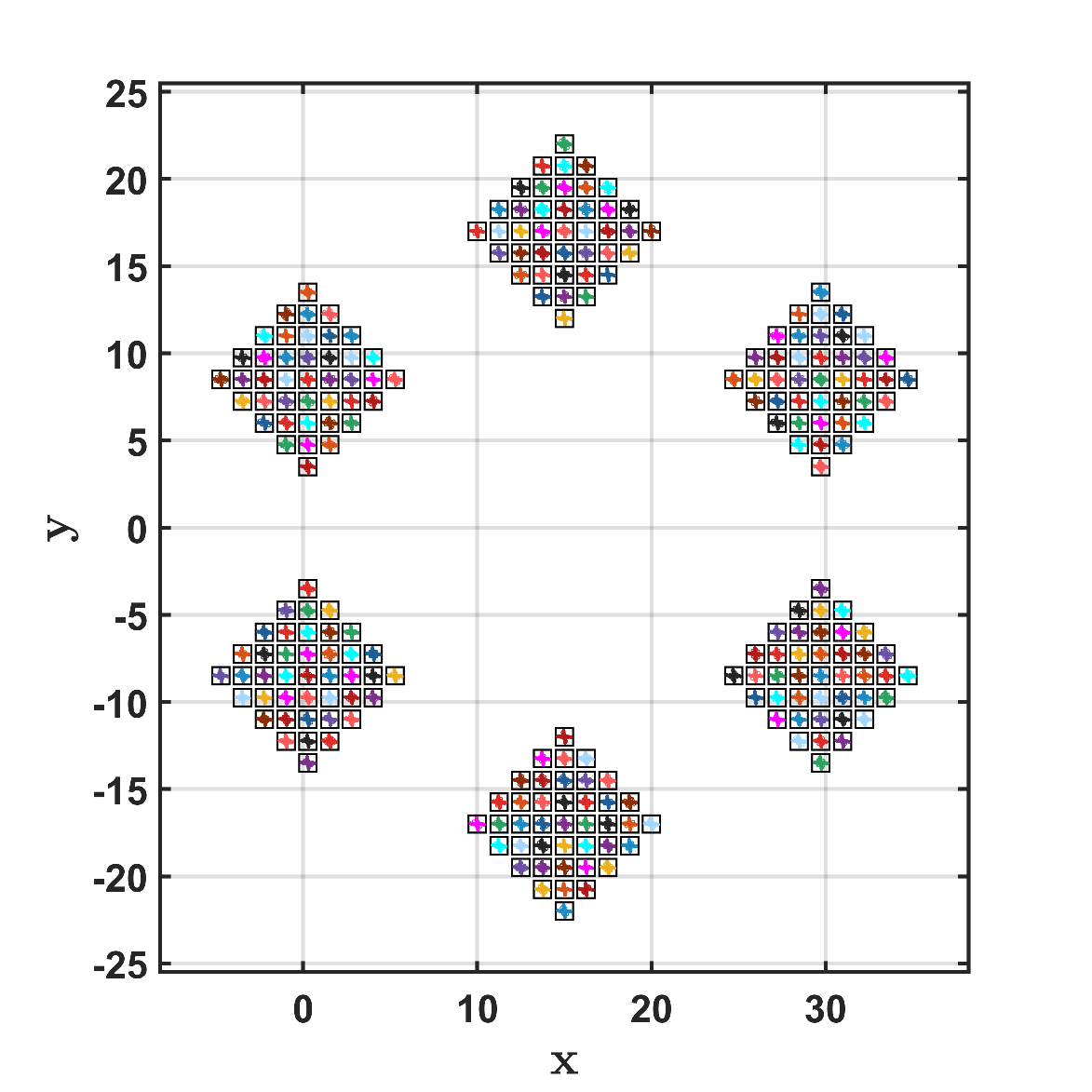}}
    \subfloat[Distance Plot  \label{9d}]{
        \includegraphics[width=0.215\linewidth, trim={0cm -1cm 0cm 0cm},clip]{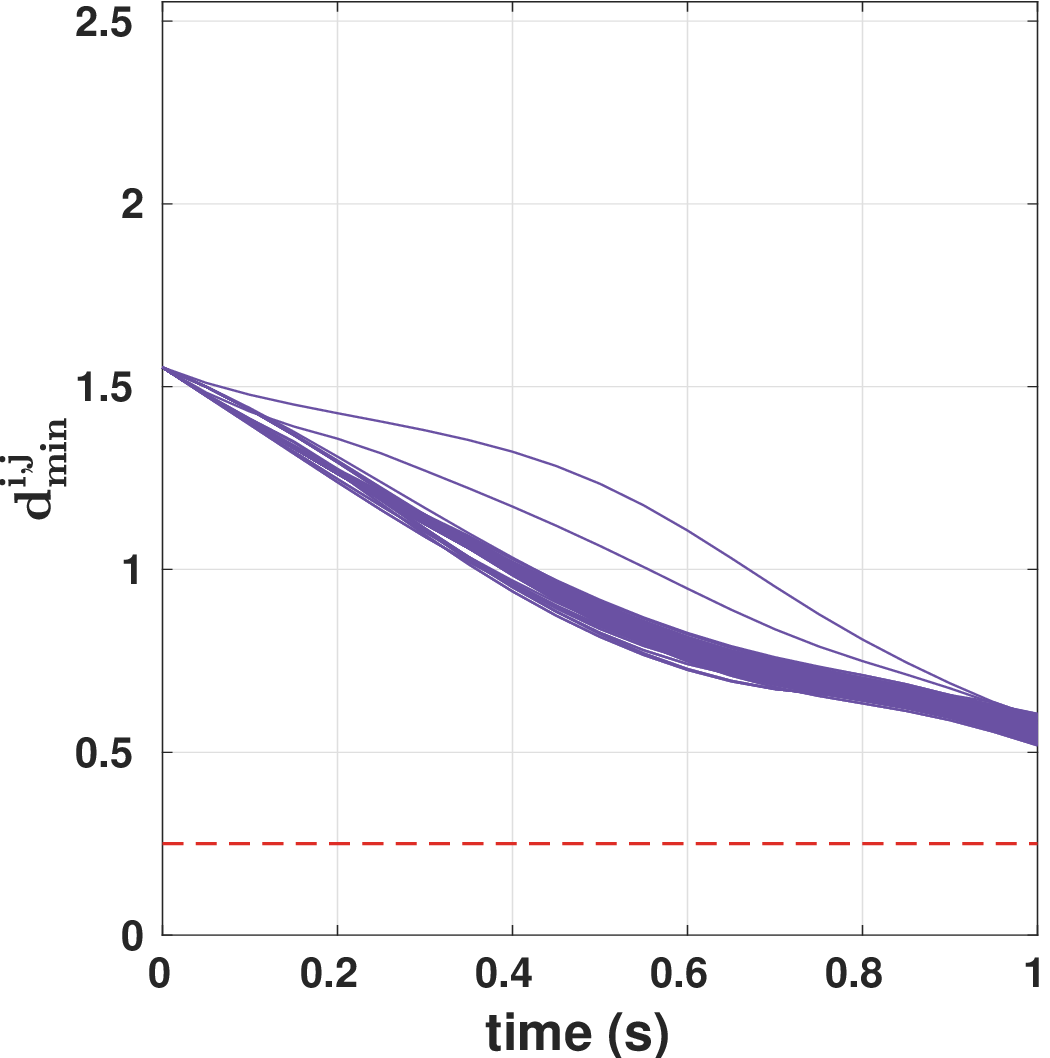}}
        \hfill
  \caption{ \textbf{Deterministic case with 246 agents:} Terminal linear \eqref{Robust Linear Constraints} and collision avoidance \eqref{robust inter-agent collision avoidance} constraints are imposed.}
  \label{fig9} 
\end{figure*}
\subsubsection{Scalability of the Proposed Framework}
We demonstrate the scalability of the proposed framework through Figs. \ref{fig7}, \ref{fig8} and \ref{fig9} in deterministic case. In Fig. \ref{fig7}, we illustrate a complex scenario with 20 agents and 21 obstacles. Each agent is tasked to move six positions to the left in the circular formation while also lying inside the shown circle at a timestep $k=15$. In Fig. \ref{fig8}, we showcase a scenario with 63 agents and 8 obstacles where each agent is tasked to reach target without colliding with obstacles.
In Fig. \ref{fig9}, we demonstrate a large-scale scenario with 246 agents, where each agent is tasked to reach the target bounds. In all three scenarios, all the agents satisfy the imposed constraints while avoiding collision with other obstacles and agents. 
Further, we provide a computational time comparison between the proposed distributed approach and the centralized approach in Table \ref{Table-2}. It can be observed that the there is not much increase in the optimizer time with the increase in agents for the distributed approach. Whereas, the optimizer time increases significantly with the number of agents for the centralized approach, and fails for scenarios involving more than 50 agents. 

All simulations were carried out in Matlab2022b \cite{MATLAB} environment using CVX \cite{cvx} as the modeling software and MOSEK as the solver \cite{mosek} on a system with an Intel Corei9-13900K. 
\begin{table}[t!]
    \centering
    \begin{tabular}{ |c|c|c|c|c|c|c| } 
\hline
 Number of Agents & 15 & 25 & 50 & 75 & 100  \\
 \hline
 Proposed Framework & 6.93 & 6.42 & 7.31 & 9.23 & 9.73  \\
 \hline
 Centralized Approach & 38.2 & 80.45 & - & - & - \\
\hline
\end{tabular}
\caption{Optimizer time comparison - Distributed Vs Centralized}
\label{Table-2}
\vspace{-0.5cm}
\end{table}
\section{Conclusion}
This article presents a decentralized robust optimization framework for multi-agent trajectory optimization under stochastic and deterministic uncertainty by integrating robust optimization, stochastic optimization, and distributed optimization techniques. 
We proposed novel robust constraint reformulations for deterministic and mixed disturbance cases, which are computationally efficient compared to prior approaches.
We introduced a decentralized framework based on CADMM with a discounted dual update step. Additionally, a convergence study is presented grounded on empirical observations. 
Further, we provide worst-case computational complexity bounds for the proposed distributed framework, prior work, and the centralized approach. 
Finally, extensive simulation experiments across varied scenarios validate the performance, robustness, and scalability of the proposed framework.
Future research directions include extending this framework to address nonlinear systems and further improve the convergence speed of the framework by leveraging distributed learning to optimize architectures \cite{saravanos2025deep}.

% \section*{Acknowledgment}

\newpage
\section*{Appendix}
\appendices
%
%
% \subsection{Compact Dynamics Matrices}
% \label{Appendix: Compact Dynamics Matrices}
% Given the matrix $\vPhi^i(k_1,k_2)$ defined as $ \vPhi^i(k_1,k_2) = A_{k_1 -1}^i A_{k_1-2}^i \dots A_{k_2}^i$ for $ k_1 > k_2$, the matrices $\vG_0^i$, $\vG_u^i$, $\vG_{\bw}^i$, and $\vG_{\zeta}^i$ are given as
% \begin{subequations}
% \begin{align}
% &
% \vG_0^i = 
% \begin{bmatrix}
% \vI; & 
% \vPhi^i(1,0) ; & 
% \vPhi^i(2,0) ; &
% \dots ; &
% \vPhi^i(T,0)
% \end{bmatrix}, 
% %\in \Rb^{(T+1) n_x \times n_x }, 
% \nonumber
% \\[0.2cm]
% &
% \bar{\vG}^i (B^i) =  \begin{bmatrix}
% \vzero & \vzero & \dots & \vzero \\ 
% B_0^i & \vzero & \dots & \vzero \\ 
% \vPhi^i(2,1) B_0^i & B_1^i & \dots & \vzero \\ 
% \vdots & \vdots & \vdots & \vdots\\ 
% \vPhi^i(T,1) B_0^i & \vPhi^i(T,2) B_1^i & \dots & B_{T-1}^i
% \end{bmatrix}
% \nonumber
% \end{align}
% \end{subequations}
% such that $\vG_u^i = \bar{\vG}^i ( B^i )$, $\vG_{\zeta}^i = \begin{bmatrix}
%     \vG_0^i, & \bar{\vG}^i (C^i )
% \end{bmatrix}$, and
% $\vG_{\bw}^i = \begin{bmatrix}
%     \vG_0^i, & \bar{\vG}^i (D^i ) \end{bmatrix}$.
% %
%
\subsection{Robust Constraint Reformulation}
%
%
% \subsubsection{Initial Approach} \label{Appendix: Sec 1 initial approach}
% \begin{subequations}
% \begin{align}
%     & \mathcal{Q}^l ( \vK^i) = (\vH^i \tilde{\vM}_k^{i,l}) \T \vH^i \tilde{\vM}_k^i \nonumber
%     % \\
%     % & \qquad \qquad \qquad
%     + (\vH^i (\tilde{\vM}_k^i - \tilde{\vM}_k^{i,l})) \T \vH^i \tilde{\vM}_k^{i,l}
%     \\
%     & \bar{\mathcal{Q}}^l ( \vK^i, \Bar{\bu}^i) = (\vH^i \bmu_{x_k,\bar{u}}^{i,l} - \bm{b}_i) \T (\vH^i \tilde{\vM}_k^i) \nonumber
%     \\
%     & \qquad \qquad \qquad
%     + (\vH^i(\bmu_{x_k,\bar{u}}^i - \bmu_{x_k,\bar{u}}^{i,l})) \T (\vH^i \tilde{\vM}_k^{i,l}) \nonumber
%     \\
%     & q^l(\Bar{\bu}^i) = (\vH^i \bmu_{x_k,\bar{u}}^{i,l} - \bm{b}_i) \T (\vH^i (2\bmu_{x_k,\bar{u}}^{i}-\bmu_{x,\bar{u}}^{i,l})  - \bm{b}_i ) \nonumber
% \end{align}
% \end{subequations}
%
%
\subsubsection{Proof of Proposition 1} \label{Appendix: Sec1 linear constraint proof}
From \eqref{ellipsoid uncertainty set}, we have $\bzeta^i = \boldsymbol{\Gamma}_i \bz_i$. Thus, we can rewrite the problem $\max_{\bzeta^i \in \calU_i} \tilde{\ba}_i \T \bzeta^i $ as follows
\begin{align}
       \max_{\bzeta^i \in \calU_i} \tilde{\ba}_i \T \bzeta^i 
       & = 
       \max_{\bz_i^T \vS_i  \bz_i  \leq \tau^i} \tilde{\ba}_i \T \boldsymbol{\Gamma}_i \bz_i \nonumber \\
       & =
       - \min_{\bz_i^T \vS_i  \bz_i  \leq \tau^i} 
       - \tilde{\ba}_i \T \boldsymbol{\Gamma}_i \bz_i
       \label{prop-1 rel1}
\end{align}
It should be observed that if $\tilde{\ba}_i \T \boldsymbol{\Gamma}_i = 0$, then $\min_{\bz_i^T \vS_i  \bz_i  \leq \tau^i} 
       - \tilde{\ba}_i \T \boldsymbol{\Gamma}_i \bz_i = 0$ since $\bz_i$ is bounded. 

       Let us now consider the case where $\tilde{\ba}_i \T \boldsymbol{\Gamma}_i \neq 0$.
   The Lagrangian for the problem 
   $\min_{\bz_i^T \vS_i  \bz_i  \leq \tau^i} - \tilde{\ba}_i \T \boldsymbol{\Gamma}_i \bz_i$ 
   is given as follows.
   \begin{align}
    \calL(\bz_i, \tilde{\lambda}) =  - \tilde{\ba}_i \T \boldsymbol{\Gamma}_i \bz_i 
    + \tilde{\lambda} (\bz_i^T \vS_i  \bz_i  - \tau^i)
    \nonumber
\end{align}
where $\tilde{\lambda} \in \Rb$ is the Lagrange multiplier. The above convex problem satisfies Slater's conditions. Thus, we can use the Karush–Kuhn–Tucker (KKT) conditions to derive the solution $(\bz^*, \tilde{\lambda}^*)$. By simplifying the condition $\nabla_{\bz_i}  \calL (\bz_i^*, \tilde{\lambda}^*) = 0$, we get
\begin{equation}
2 \tilde{\lambda}^* \bz_i^* = \vS_i^{-1} \boldsymbol{\Gamma}_i \T \tilde{\ba}_i
\nonumber
\end{equation}
Since $\tilde{\ba}_i \T \boldsymbol{\Gamma}_i \neq 0$ and $\vS_i \in \Sb_{++}^{\bar{n}_i}$, the RHS of the above equation is non-zero, which implies that both $\bz_i^*$ and $\tilde{\lambda}^*$ need to be non-zero. Therefore, we obtain
\begin{equation}
\bz_i^* 
= (2 \tilde{\lambda}^*)^{-1} \vS_i^{-1} \boldsymbol{\Gamma}_i^T \tilde{\ba}_i
\label{prop 1 proof z_i star}
\end{equation}
Combining the fact $\tilde{\lambda}^* \neq 0$, and the KKT slackness condition $\tilde{\lambda}^* (\bz_i^*{}\T \vS_i  \bz_i^* - \tau^i) = 0$, we get $\bz_i^*{}\T \vS_i  \bz_i^* - \tau^i = 0$. Now, substituting \eqref{prop 1 proof z_i star} in the aforementioned equation, we get
\begin{equation}
2 \tilde{\lambda}^* = (\tau^i)^{-1/2} \; \| \boldsymbol{\Gamma}_i \T \tilde{\ba}_i \|_{\vS_i^{-1}},
\label{prop-1 rel 2}
\end{equation}
Using \eqref{prop 1 proof z_i star} and \eqref{prop-1 rel 2}, we can obtain $\calL(\bz_i^*, \tilde{\lambda}^*)$ as follows.
\begin{equation}
    \calL (\bz_i^*, \tilde{\lambda}^*) = - \ba_i \T \vM_i \boldsymbol{\Gamma}_i \bz_i^* 
    =
    - (\tau^i)^{1/2} \| \boldsymbol{\Gamma}_i \T \tilde{\ba}_i \|_{\vS_i^{-1}}.
    \label{prop-1 rel3}
\end{equation}
Therefore, combining \eqref{prop-1 rel1} and \eqref{prop-1 rel3}, we get
\begin{align}
    \max_{\bzeta^i \in \calU_i} \tilde{\ba}_i \T \bzeta^i 
    = - \calL (\bz_i^*, \tilde{\lambda}^*)
    = (\tau^i)^{1/2} \| \boldsymbol{\Gamma}_i \T \tilde{\ba}_i \|_{\vS_i^{-1}} \nonumber
\end{align}
Using this, we can further derive the following.
\begin{equation}
     \min_{\bzeta^i \in \calU_i} \tilde{\ba}_i \T \bzeta^i
     = - \max_{\bzeta^i \in \calU_i} - \tilde{\ba}_i \T \bzeta^i
     = -(\tau^i)^{1/2} \| \boldsymbol{\Gamma}_i \T \tilde{\ba}_i \|_{\vS_i^{-1}} \nonumber
\end{equation}
\subsubsection{Proof of Proposition 2} \label{Appendix: Sec 1 nonconvex obstacle mu_d}
    The first step is to construct a tight upper bound for $\| \vH^i \tilde{\vM}_k^i \bzeta^i \|_2^2$. For that, we use the following relation
%\begin{subequations}
\begin{align}
    \max_{\bzeta^i \in \calU_i}  \| \vH^i \tilde{\vM}_k^i \bzeta^i \|_2 
    & = 
    \max_{\bzeta^i \in \calU_i} \bigg(  \sum_{\bar{m} = 1}^{n_{\text{pos}}} 
    (\bm{h}_{k,\bar{m}}^i {}\T \vM_i \bzeta^i)^2 \bigg)^{\frac{1}{2}} \nonumber \\
    & 
    \leq 
    \bigg( \sum_{\bar{m} = 1}^{n_{\text{pos}}} 
    \max_{\bzeta^i \in \calU_i} (\bm{h}_{k,\bar{m}}^i {}\T \vM_i \bzeta^i)^2 \bigg)^{\frac{1}{2}} 
    \label{prop 2 rel1}
    % \\
    % & ~~ =
    % \bigg( \sum_{\bar{m} = 1}^{m} \tau^i \| \boldsymbol{\Gamma}_i \T \vM_i \T \bm{h}_{k,\bar{m}}^i \|_{\vS_i^{-1}}^2 
    %  \bigg)^{1/2} \nonumber
\end{align} 
%\end{subequations}
where $\bm{h}_{k,\bar{m}}^i {}\T $ is $\bar{m}^{th}$ row of $ \vH^i \vP_k^i$.  
From Proposition \ref{prop: Max bound on mean}, we have $\max_{\bzeta^i \in \calU_i} \bm{h}_{k,\bar{m}}^i {}\T \vM_i \bzeta^i 
= 
- \min_{\bzeta^i \in \calU_i} \bm{h}_{k,\bar{m}}^i {}\T \vM_i \bzeta^i$ which implies the following
\begin{equation}
    \max_{\bzeta^i \in \calU_i} (\bm{h}_{k,\bar{m}}^i {}\T \vM_i \bzeta^i)^2
    =
    \bigg( \max_{\bzeta^i \in \calU_i} \bm{h}_{k,\bar{m}}^i {}\T \vM_i \bzeta^i \bigg)^2
    \label{prop-2 rel2}
\end{equation}
Further, from Proposition \ref{prop: Max bound on mean} we also get the following
\begin{align}
\max_{\bzeta^i \in \calU_i} \bm{h}_{k,\bar{m}}^i {}\T \vM_i \bzeta^i 
& = (\tau^i)^{1/2} \| \boldsymbol{\Gamma}_i \T \vM_i \T \bm{h}_{k,\bar{m}}^i \|_{\vS_i^{-1}},
\label{prop 2 proof min}
\end{align}
Next, by introducing a variable $\bmu_{d,k}^i$ such that 
\begin{equation}
(\tau^i)^{1/2} \| \boldsymbol{\Gamma}_i \T \vM_i \T \bm{h}_{k,\bar{m}}^i\|_{\vS_i^{-1}} 
\leq (\bmu_{d,k}^i)_{\bar{m}}, \nonumber
\end{equation}
for every $\bar{m} = 1, \dots, n_{\text{pos}}$, and by combining \eqref{prop 2 rel1}, \eqref{prop-2 rel2}, and \eqref{prop 2 proof min}, we obtain
\begin{align}
    \max_{\bzeta^i \in \calU_i}
    \| \vH^i \tilde{\vM}_k^i \bzeta^i \|_2  \leq \|\bmu_{d,k}^i {} \|_2 . \nonumber
\end{align}
Therefore, we can consider the following tighter approximation of the constraint \eqref{nonconvex obstacle K1}, given by
\begin{align}
    \| \bmu_{d,k}^i {} \|_2  \leq \tilde{c}_i^k. \nonumber
\end{align}
\subsubsection{Proof of Proposition 3} \label{Appendix: sec1 nonconvex interagent}
This proof is similar to the one of Proposition \ref{prop- obs avoidance}. We start with the following relation.
\begin{equation}
\begin{aligned}
    & \max_{\bzeta^i \in \calU_i, \bzeta^j \in \calU_j}
          \| \vH^i \tilde{\vM}_k^i \bzeta^i 
         -  \vH^j \tilde{\vM}_k^j \bzeta^j \|_2 \\
    & \quad \leq
    \bigg( \sum_{\bar{m}= 1}^{n_{\text{pos}}} \max_{\substack{\bzeta^i \in \calU_i, \\ \bzeta^j \in \calU_j}} 
    (\bm{h}_{k,\bar{m}}^i {}\T \tilde{\vM}_k^i \bzeta^i 
         -  \bm{h}_{k,\bar{m}}^j {}\T \tilde{\vM}_k^j \bzeta^j)_{\bar{m}}^2 \bigg)^{\frac{1}{2}} 
        \label{prop-3 rel1}
\end{aligned} 
\end{equation}
where $\bm{h}_{k,\bar{m}}^i {}\T, \bm{h}_{k,\bar{m}}^j {}\T $ are the $\bar{m}^{th}$ rows of the matrices $\vH^i \vP_k^i, \vH^j \vP_k^j$ respectively.
Let us now introduce the variables $\bmu_{d,k}^i, \bmu_{d,k}^j$ such that
\begin{align}
    & 
    (\bmu_{d,k}^i)_{\bar{m}} \geq  
    (\tau^i)^{1/2} \| \boldsymbol{\Gamma}_i \T \vM_i \T \bm{h}_{k,\bar{m}}^i\|_{\vS_i^{-1}} \nonumber \\
    &
    (\bmu_{d,k}^j)_{\bar{m}} \geq  
    (\tau^j)^{1/2} \| \boldsymbol{\Gamma}_j \T \vM_j \T \bm{h}_{k,\bar{m}}^j\|_{\vS_j^{-1}} \nonumber
\end{align}
for every $\bar{m} = 1, \dots, n_{\text{pos}}$. Using proposition \ref{prop: Max bound on mean}, 
we get
\begin{equation}
\begin{aligned}
    & \max_{\bzeta^i \in \calU_i, \bzeta^j \in \calU_j}  \bm{h}_{k,\bar{m}}^i {}\T \tilde{\vM}_k^i \bzeta^i 
         -  \bm{h}_{k,\bar{m}}^j {}\T \tilde{\vM}_k^j \bzeta^j
    \\
    &~~ =
    \max_{\bzeta^i \in \calU_i} \bm{h}_{k,\bar{m}}^i {}\T \tilde{\vM}_k^i \bzeta^i
    - \min_{\bzeta^j \in \calU_j} \bm{h}_{k,\bar{m}}^j {}\T \tilde{\vM}_k^j \bzeta^j
    \\
    &~~ = 
    - \min_{\bzeta^i \in \calU_i, \bzeta^j \in \calU_j}  (\bm{h}_{k,\bar{m}}^i {}\T \tilde{\vM}_k^i \bzeta^i 
         -  \bm{h}_{k,\bar{m}}^j {}\T \tilde{\vM}_k^j \bzeta^j ), \nonumber
\end{aligned}
\end{equation}
and employing a similar approach as in the proof of proposition \ref{prop- obs avoidance}, we obtain the following relation from \eqref{prop-3 rel1} -
\begin{equation}
\begin{aligned}
    \max_{\bzeta^i \in \calU_i, ~ \bzeta^j \in \calU_j}
         \| \vH^i \tilde{\vM}_k^i \bzeta^i 
         & -  \vH^j \tilde{\vM}_k^j \bzeta^j \|_2 \\
    & \qquad \qquad
    \leq
    \| \bmu_{d,k}^i + \bmu_{d,k}^j \|_2. \nonumber
\end{aligned}
\end{equation}
Thus, we consider the following tighter approximation to the constraint \eqref{interagent collision avoidance K1},
\begin{align}
    \| \bmu_{d,k}^i + \bmu_{d,k}^j \|_2 \leq  \tilde{c}_{ij}^k
\end{align}
\subsubsection{Proof of Proposition 4} \label{Appendix: sec1 convex norm s-lemma}
We would start by rewriting the constraint \eqref{convex norm-of-mean constraints: interm1} in terms of $\bz^i$ (from \eqref{ellipsoid uncertainty set}) as follows
\begin{align}
    \max_{\bz_i^T \vS_i  \bz_i \leq \tau^i} \| \vH^i\tilde{\vM}_{\tilde{k}}^i \boldsymbol{\Gamma}_i \bz^i 
    + \hat{\bmu}^i_{\text{pos},k} - \bm{p}_{\text{tar}} \|_2^2 
    \leq c_i^2. 
    \label{Appendix: sec1 prop 5 eq1}
\end{align}
Let us now rewrite the RHS term of the above as following
\begin{equation}
\begin{aligned}
    & \| \vH^i\tilde{\vM}_{\tilde{k}}^i \boldsymbol{\Gamma}_i \bz^i 
    + \hat{\bmu}^i_{\text{pos},k} - \bm{p}_{\text{tar}} \|_2^2 \\
    & \quad 
    = \bz_i^T \vF_{1,k}^i(\vK^i) \bz_i
    + 2 \vF_{2,k}^i (\vK^i, \bar{\bu}^i) \bz_i
    + \vF_{3,k}^i(\bar{\bu}^i), \nonumber
\end{aligned}
\end{equation}
where $\vF_{1,k}^i(\vK^i) = 
    (\vH^i\tilde{\vM}_{\tilde{k}}^i\boldsymbol{\Gamma}_i)^T (\vH^i\tilde{\vM}_{\tilde{k}}^i \boldsymbol{\Gamma}_i)$, 
    $\vF_{2,k}^i (\vK^i, \bar{\bu}^i) =
    (\hat{\bmu}^i_{\text{pos},k} - \bm{p}_{\text{tar}})^T
    \vH^i\tilde{\vM}_{\tilde{k}}^i\boldsymbol{\Gamma}_i$, and $\vF_{3,k}^i(\bar{\bu}^i) =
    \| \hat{\bmu}^i_{\text{pos},k} - \bm{p}_{\text{tar}} \|_2^2$.
%
% Further, the constraint \eqref{Appendix: sec1 prop 5 eq1} can be equivalently given as 
% \begin{equation}
% \begin{aligned}
%      & \bz_i^T \vS_i  \bz_i \leq \tau^i \\
%     & \qquad
%     \implies 
%     \| \vH^i\tilde{\vM}_{\tilde{k}}^i \boldsymbol{\Gamma}_i \bz^i 
%     + \hat{\bmu}^i_{\text{pos},k} - \bm{p}_{\text{tar}} \|_2^2 
%     \leq c_i^2. \nonumber
% \end{aligned}
% \end{equation}
%
Using S-lemma, the constraint \eqref{Appendix: sec1 prop 5 eq1} can be equivalently given by the following set of constraints
\begin{align}
&
\begin{bmatrix}
    -  \vF_{1,k}^i(\vK^i)  + \alpha \vS_i 
    & - \vF_{2,k}^i (\vK^i, \bar{\bu}^i)  \T  \\
    - \vF_{2,k}^i (\vK^i, \bar{\bu}^i)  & -\vF_{3,k}^i (\bar{\bu}^i) + c_i^2 - \alpha \tau^i
\end{bmatrix}
 \succeq 0,
 \label{convex norm constraints - non-convex equivalent approach 1}
 \\[0.2cm]
 & ~~~~~~~~~~~~~~~~~~~~~~~~~~~ \alpha \geq 0.
\end{align}
It can be observed that the constraint \eqref{convex norm constraints - non-convex equivalent approach 1} is nonconvex because the terms $\vF_{1,k}^i(\vK^i)$, $\vF_{2,k}^i (\vK^i, \bar{\bu}^i)$, and $\vF_{3,k}^i (\bar{\bu}^i)$ are non-linear in the decision variables. Thus, we need to further reformulate the above constraint. For that, we start by rewriting the constraint \eqref{convex norm constraints - non-convex equivalent approach 1} as follows
\begin{equation}
\begin{aligned}
& 
\begin{bmatrix}
    \alpha \vS_i 
    & 0 \\
    0 &  c_i^2 - \alpha \tau^i
\end{bmatrix} \succeq 
\begin{bmatrix}
    \vF_{1,k}^i(\vK^i) 
    & \vF_{2,k}^i (\vK^i, \bar{\bu}^i)  \T  \\
    \vF_{2,k}^i (\vK^i, \bar{\bu}^i) 
    & \vF_{3,k}^i (\bar{\bu}^i)
\end{bmatrix} \nonumber
\end{aligned}
\end{equation}
which can be further rewritten as 
\begin{equation}
\begin{aligned}
& 
\begin{bmatrix}
    \alpha \vS_i 
    & 0 \\
    0 &  c_i^2 - \alpha \tau^i
\end{bmatrix} 
- \vF_{4,k}^i (\vK^i, \bar{\bu}^i)^T \vF_{4,k}^i (\vK^i, \bar{\bu}^i) \succeq 0,
\end{aligned}
\label{convex norm S-lemma constraint}
\end{equation}
where $ \vF_{4,k}^i (\vK^i, \bar{\bu}^i) = [
    \vH^i\tilde{\vM}_{\tilde{k}}^i \boldsymbol{\Gamma}_i, 
    \hat{\bmu}^i_{\text{pos},k} - \bm{p}_{\text{tar}}]$.
Using Schur Complement, the constraint \eqref{convex norm S-lemma constraint} can be equivalently given as
\begin{align}
    \begin{bmatrix}
    \alpha \vS_i & 0 
    & (\vH^i \tilde{\vM}_{\tilde{k}}^i \boldsymbol{\Gamma}_i) \T \\
    0 &  c_i^2 - \alpha \tau^i 
    & (\hat{\bmu}^i_{\text{pos},k} - \bm{p}_{\text{tar}} ) \T \\
    \vH^i \tilde{\vM}_{\tilde{k}}^i \boldsymbol{\Gamma}_i
    & \hat{\bmu}^i_{\text{pos},k} - \bm{p}_{\text{tar}} 
    & \vI_{n_{\text{pos}}}
    \end{bmatrix}
    \succeq 0 \nonumber
\end{align}
%
%
%
% \subsection{Mixed Case: Robust Constraint Reformulation}
%
\subsubsection{Proof of Proposition 5} \label{Appendix: sec2 proposition 5}
Let us rewrite the constraint \eqref{mixed: nonconvex norm chance wbar interm2 - appendix} as follows
\begin{equation}
\begin{aligned}
    & \max_{
    \substack{\bzeta^i  \in \calU_i,  \hat{\bw}^i \in \calW_i[\bar{\eta}_i] }} 
    \| \vH^i \tilde{\vM}_k^i \bzeta^i
     +  \bSigma_{\bar{\bw}^i}^{1/2} \hat{\bw}^i \|_2 
    \leq \tilde{c}_i^k .
\end{aligned}
\label{mixed nonconvex obstacle determinisitic interm}
\end{equation}
We use a similar intuition as in Proposition \ref{prop- obs avoidance} to derive a tractable version of the above constraint. First, we write the following relation -
\begin{equation}
\begin{aligned}
    & \max_{
    \substack{\bzeta^i  \in \calU_i, \hat{\bw}^i \in \calW_i[\bar{\eta}_i] }} 
    \| \vH^i \tilde{\vM}_k^i \bzeta^i
     +  \bSigma_{\bar{\bw}^i}^{1/2} \hat{\bw}^i \|_2 
    \\
    & \quad \leq
     \bigg( \sum_{\bar{m} = 1}^{n_{\text{pos}}} 
    \max_{
    \substack{\bzeta^i  \in \calU_i, 
     \hat{\bw}^i \in \calW_i[\bar{\eta}_i] }} 
    \big[ \boldsymbol{e}_{\bar{m}} {}\T ( \vH^i \tilde{\vM}_k^i \bzeta^i
     +  \bSigma_{\bar{\bw}^i}^{1/2} \hat{\bw}^i ) \big]^2 \bigg)^{\frac{1}{2}} \nonumber
\end{aligned}
\end{equation}
where $\boldsymbol{e}_{\bar{m}} {}\T$ is $\bar{m}^{th}$ row of $\vI_{n_{\text{pos}}}$.
Then, we establish the following relation -
\begin{equation}
\begin{aligned}
     & \max_{
    \substack{\bzeta^i  \in \calU_i,  \hat{\bw}^i \in \calW_i[\bar{\eta}_i] }} 
    \boldsymbol{e}_{\bar{m}} {}\T ( \vH^i \tilde{\vM}_k^i \bzeta^i
     +  \bSigma_{\bar{\bw}^i}^{1/2} \hat{\bw}^i )
     \\[0.2cm]
     &~~ =
     \max_{\bzeta^i  \in \calU_i} 
     \bm{h}_{k,\bar{m}}^i {}\T \vM_i \bzeta^i
     + 
    \max_{ \hat{\bw}^i \in \calW_i[\bar{\eta}_i] } 
    \boldsymbol{e}_{\bar{m}} {}\T \bSigma_{\bar{\bw}^i}^{1/2} \hat{\bw}^i
\end{aligned}
\label{prop 5 appendix interm3}
\end{equation}
where $\bm{h}_{k,\bar{m}}^i {}\T = \boldsymbol{e}_{\bar{m}} {}\T \vH^i \vP_k^i$ (which is $\bar{m}^{th}$ row of the matrix $\vH^i \vP_k^i$). 
It should be noted that the disturbance set $\calW_i[\bar{\eta}_i ]$ is an ellipsoidal set of the same form as $ \calU_i$. Thus, using Proposition \ref{prop: Max bound on mean}, we have
\begin{equation}
\begin{aligned}
    & \max_{ \hat{\bw}^i \in \calW_i[\bar{\eta}_i] } 
    \boldsymbol{e}_{\bar{m}} {}\T \bSigma_{\bar{\bw}^i}^{1/2} \hat{\bw}^i
    = - \min_{ \hat{\bw}^i \in \calW_i[\bar{\eta}_i] } 
    \boldsymbol{e}_{\bar{m}} {}\T \bSigma_{\bar{\bw}^i}^{1/2} \hat{\bw}^i 
\end{aligned}
\end{equation}
Now, using the above relation and \eqref{prop 5 appendix interm3}, we get
\begin{equation}
\begin{aligned}
    & \max_{
    \substack{\bzeta^i  \in \calU_i, \hat{\bw}^i \in \calW_i[\bar{\eta}_i] }} 
   \boldsymbol{e}_{\bar{m}} {}\T ( \vH^i \tilde{\vM}_k^i \bzeta^i
     +  \bSigma_{\bar{\bw}^i}^{1/2} \hat{\bw}^i )
     \\
     & \qquad 
    = - \min_{
    \substack{\bzeta^i  \in \calU_i, \hat{\bw}^i \in \calW_i[\bar{\eta}_i] }} 
    \boldsymbol{e}_{\bar{m}} {}\T ( \vH^i \tilde{\vM}_k^i \bzeta^i
     +  \bSigma_{\bar{\bw}^i}^{1/2} \hat{\bw}^i )
\end{aligned}
\end{equation}
By establishing the above relation, we can use the same steps as in Proposition \ref{prop- obs avoidance} to derive approximate tighter version of the constraint \eqref{mixed nonconvex obstacle determinisitic interm} as follows -
\begin{align}
    & \| \bmu_{w,k}^i {} \|_2  \leq \tilde{c}_i^k, \\
    & 
    \begin{aligned}
    & (\bmu_{w,k}^i)_{\bar{m}} \geq 
    \max_{
    \substack{\bzeta^i  \in \calU_i, 
    \\ \hat{\bw}^i \in \calW_i[\bar{\eta}_i] }} 
    \boldsymbol{e}_{\bar{m}} {}\T ( \vH^i \tilde{\vM}_k^i \bzeta^i
     +  \bSigma_{\bar{\bw}^i}^{1/2} \hat{\bw}^i )
    \\[-0.4cm]
    & \qquad \qquad \qquad \qquad \qquad \qquad \qquad \quad
    \forall \; \bar{m} \in \llbracket 1, n_{\text{pos}} \rrbracket. 
    \end{aligned}
    \label{prop 5 nonconvex obs interm4}
\end{align}
% We now simplify the constraint \eqref{prop 5 nonconvex obs interm4}. 
Now, using Proposition \ref{prop: Max bound on mean}, we have
\begin{align}
    & \max_{\bzeta^i \in \calU_i} \bm{h}_{k,\bar{m}}^i {}\T \vM_i \bzeta^i 
    = \sqrt{\tau^i} \| \boldsymbol{\Gamma}_i \T \vM_i \T \bm{h}_{k,\bar{m}}^i \|_{\vS_i^{-1}}
    \label{max limit deter part} \\
    & \max_{\hat{\bw}^i \in \calW_i[\bar{\eta}_i] }
    \boldsymbol{e}_{\bar{m}} {}\T \bSigma_{\bar{\bw}^i}^{1/2} \hat{\bw}^i 
    =
    \sqrt{\bar{\eta}_i} \| \bSigma_{\bar{\bw}^i}^{1/2} \boldsymbol{e}_{\bar{m}} \|_2
    \label{max limit sto part}
\end{align}
Substituting the fact that $\bSigma_{\bar{\bw}^i} = \vH^i \tilde{\vM}^{w,i}_k \bSigma_{\boldsymbol{\bw}^i}  \tilde{\vM}^{w,i}_k{}^T \vH^i{}^T$, we can rewrite the constraint \eqref{max limit sto part} as follows 
\begin{equation}
\begin{aligned}
    \max_{\hat{\bw}^i \in \calW_i[\bar{\eta}_i] }
    \boldsymbol{e}_{\bar{m}} {}\T \bSigma_{\bar{\bw}^i}^{1/2} \hat{\bw}^i 
    =
    \sqrt{\bar{\eta}_i} \| \boldsymbol{\varphi}^i \vM^{w}_i {}\T \bm{h}_{k,\bar{m}}^i \|_2
\end{aligned}
 \label{max limit sto part - final}
\end{equation}
By using \eqref{max limit deter part} and \eqref{max limit sto part - final}, we can rewrite the constraint \eqref{prop 5 nonconvex obs interm4} for all $\bar{m} \in \llbracket 1, n_{\text{pos}} \rrbracket$ as follows -
\begin{equation}
\begin{aligned}
    (\bmu_{w,k}^i)_{\bar{m}} 
    \geq 
     \sqrt{\tau^i} \| \boldsymbol{\Gamma}_i \T \vM_i \T \bm{h}_{k,\bar{m}}^i \|_{\vS_i^{-1}} 
    % \\
    % & \qquad \qquad
    +
    \sqrt{\bar{\eta}_i} \| \boldsymbol{\varphi}^i \vM^{w}_i {}\T \bm{h}_{k,\bar{m}}^i \|_2 \nonumber
\end{aligned}
\end{equation}
\subsubsection{Proof of Proposition 6} \label{Appendix: Proof of Proposition 6 nonconvex interagent chance}
% We combine the intuition from Propositions \ref{Proposition: nonconvex interagent deterministic} and \ref{proposition: nonconvex chance} for this derivation. 
% First, let us rewrite the constraint \eqref{mixed: interagent collision avoidance eq2} as follows
% \begin{equation}
% \begin{aligned}
%     & \forall ~ \bzeta^i \in \calU_i, ~ \bzeta^j \in \calU_j, \\
%     & \quad 
%     \Pb \big( \| \vH^i \tilde{\vM}_k^i \bzeta^i
%     -  \vH^j \tilde{\vM}_k^j \bzeta^j
%     \\
%     & \qquad \quad
%     + \bSigma_{\bar{\bw}^i}^{1/2} \hat{\bw}^i 
%     - \bSigma_{\bar{\bw}^j}^{1/2} \hat{\bw}^j \|_2  
%     \leq \tilde{c}_{ij}^k \big) \geq 1-p_{ij}
% \end{aligned}
% \label{Appendix: proposition 6 interm1}
% \end{equation}
%
Similar to Proposition \ref{proposition: nonconvex chance}, we can assert that if the constraint within the probability bracket (of \eqref{mixed: interagent collision avoidance eq2}) is satisfied for all stochastic disturbances $\hat{\bw}^i$ and $\hat{\bw}^j$ lying inside their respective $100 (\sqrt{1- p_{ij}})$ confidence ellipsoids, then the constraint \eqref{mixed: interagent collision avoidance eq2} is satisfied. Thus, we consider the following approximation of the constraint 
\begin{equation}
\begin{aligned}
    &
    \| \vH^i \tilde{\vM}_k^i \bzeta^i
    -  \vH^j \tilde{\vM}_k^j \bzeta^j
    + \bSigma_{\bar{\bw}^i}^{1/2} \hat{\bw}^i 
    - \bSigma_{\bar{\bw}^j}^{1/2} \hat{\bw}^j \|_2  
    \leq \tilde{c}_{ij}^k, 
    \\
    &~~~ 
    \forall \bzeta^i \in \calU_i, ~ \bzeta^j \in \calU_j, 
    ~ \hat{\bw}^{i} \in \calW_{i} [\bar{\eta}_{ij}], 
    ~ \hat{\bw}^{j} \in \calW_{j} [\bar{\eta}_{ij}] \nonumber
\end{aligned}
\end{equation}
where $\calW_{i} [\bar{\eta}_{ij}]$ and $\calW_{j} [\bar{\eta}_{ij}]$ are $100 (\sqrt{1- p_{ij}})$ confidence ellipsoids of $\hat{\bw}^i$ and $\hat{\bw}^j$ respectively, and are defined as in \eqref{Stochastic noise uncertainty set - confidence ellipsoid}. 
Now, using a similar intuition from the Propositions \ref{Proposition: nonconvex interagent deterministic} and \ref{proposition: nonconvex chance}, we approximate the above constraint by the following
 \begin{align}
        & 
       \| \bar{\bmu}_{w,k}^i +  \bar{\bmu}_{w,k}^j \|_2
        \leq 
        \tilde{c}_{ij}^k, \nonumber
        \\
        &
        (\bar{\bmu}_{w,k}^i)_{\bar{m}} \geq  
        \sqrt{\tau^i} || \boldsymbol{\Gamma}_i \T \vM_i \T \bm{h}_{k,\bar{m}}^i||_{\vS_i^{-1}}, \nonumber
        % \nonumber \\
        % & \qquad \qquad \qquad
        + \sqrt{\bar{\eta}_{ij}} \|  \boldsymbol{\varphi}^i \vM_w^i {}\T h_{k,\bar{m}}^i \|_2
        \label{mixed: nonconvex interagent mu_di lowerbound} \\
        &
        (\bar{\bmu}_{w,k}^j)_{\bar{m}} \geq  
        \sqrt{\tau^j} || \boldsymbol{\Gamma}_j \T \vM_j \T \bm{h}_{k,\bar{m}}^j ||_{\vS_j^{-1}}, 
        % \nonumber \\
        % & \qquad \qquad \qquad
        + \sqrt{\bar{\eta}_{ij}} \|  \boldsymbol{\varphi}^j \vM_w^j {}\T h_{k,\bar{m}}^j \|_2 \nonumber
    \end{align}
where $\bm{h}_{k,\bar{m}}^i$ and $\bm{h}_{k,\bar{m}}^j$ are $\bar{m}^{th}$ rows of the matrices $\vH^i \vP_k^i$ and $\vH^j \vP_k^j$ respectively.
\subsection{Distributed Robust Optimization Framework}
\subsubsection{Derivation of Global Updates} \label{Appendix: sec D Derivation of Global Updates}
The global update $\bnu^{i,l+1}$can be obtained by solving the following problem
\begin{align}
    \min_{\{\bnu^i\}_{i \in \calV}} 
    \sum_{i \in \calV}
    \blambda^{i,l} {}\T \mathrm{\vC} ( \Tilde{\bmu}^{i,l+1} - \Tilde{\bnu}^i)
    + \frac{\rho}{2} || \mathrm{\vC}( \Tilde{\bmu}^{i,l+1} - \Tilde{\bnu}^i) ||_2^2 \nonumber
\end{align}
which can be solved in a decentralized manner where each agent $i$ solves the following problems to obtain $\bnu^{i, l+1}_{\text{pos}}$ and $\bnu^{i, l+1}_{d}$
\begin{equation}
\begin{aligned}
    & \min_{\bnu^{i}_{\text{pos}}}
    \blambda^{i,l}_{\text{p}}{}^T ( \hat{\bmu}^{i, l+1}_{\text{pos}} 
    - \bnu^{i}_{\text{pos}})
    + \frac{\rho}{2} ||   \hat{\bmu}^{i, l+1}_{\text{pos}} 
    - \bnu^{i}_{\text{pos}} ||_2^2
    \\
    &~~~~~ 
    + \sum_{j \in \calP_i} \big[ \blambda^{j,l}_{i, \text{p}}{}^T (  \Bar{\bmu}^{j, l+1}_{i, \text{pos}} 
    - \bnu^{i}_{\text{pos}} )
    + \frac{\rho}{2} ||    \Bar{\bmu}^{\hat{j},l+1}_{i, \text{pos}} 
    - \bnu^{i}_{\text{pos}} ||_2^2  \big] \nonumber
    \\
    & \min_{\bnu^{i}_{d}}
    \blambda^{i,l}_{d}{}^T ( \hat{\bmu}^{i, l+1}_{d} 
    - \bnu^{i}_{\text{pos}})
    + \frac{\rho \rho_r}{2} ||   \hat{\bmu}^{i, l+1}_{d} 
    - \bnu^{i}_{d} ||_2^2
    \\
    &~~~~~ 
    + \sum_{j \in \calP_i} \big[ \blambda^{j,l}_{i, d}{}^T (  \Bar{\bmu}^{j, l+1}_{i, d} 
    - \bnu^{i}_{d} )
    + \frac{\rho \rho_r}{2} ||    \Bar{\bmu}^{\hat{j},l+1}_{i, d} 
    - \bnu^{i}_{d} ||_2^2  \big] \nonumber
\end{aligned}
\end{equation}
The above problems have closed form solutions. 
For simplicity of analysis, we only derive closed form for  $\Tilde{\bnu}_{ \text{pos}}^{i,l+1}$ update. 
Since $\bnu^{i, l+1}_{\text{pos}}$ is the minimizer of the problem considered above, we have
\begin{equation}
\begin{aligned}
    \blambda^{i,l}_{\text{p}}
    + & \rho ( \hat{\bmu}^{i, l+1}_{\text{pos}} - \bnu^{i, l+1}_{\text{pos}} )
    \\
    & + \sum_{j \in \calP_i} \big[ \blambda^{j,l}_{i, \text{p}}
    + \rho ( \Bar{\bmu}^{\hat{j},l+1}_{i, \text{pos}} 
    - \bnu^{i, l+1}_{\text{pos}}
    )\big] = 0 
\end{aligned}
\label{Appendix: Global update eq1}
\end{equation}
From the dual update step in the algorithm \ref{alg1}, we can rewrite the above as follows -
\begin{equation}
\begin{aligned}
    \blambda^{i,l+1}_{\text{p}} + \delta_{\text{pos}} \blambda^{i,l}_{\text{p}}
    + \sum_{j \in \calP_i} \big[ \blambda^{j,l+1}_{i, \text{p}}
    + \delta_{\text{pos}} \blambda^{j,l}_{i, \text{p}} \big] = 0
    \nonumber
\end{aligned}
\end{equation}
If we initialize the dual variables $\{ \blambda^{i,0} \}_{i \in \calV}$ such that $\blambda^{i,0}_{\text{p}} 
    + \sum_{j \in \calP_i} \blambda^{j,0}_{i, \text{p}} = 0$ for all $i \in \calV$, then we have for all ADMM iterations $l$, $\blambda^{i,l}_{\text{p}} 
    + \sum_{j \in \calP_i} \blambda^{j,l}_{i, \text{p}} = 0$.
Substituting this in \eqref{Appendix: Global update eq1} would give us the $\bnu^{i, l+1}_{\text{pos}}$ update as follows 
\begin{equation}
    \bnu^{i, l+1}_{\text{pos}}
    =
    \frac{1}{1 + n(\calP_i)} \bigg( \hat{\bmu}^{i, l+1}_{\text{pos}} 
    + \sum_{j \in \calP_i} \Bar{\bmu}^{\hat{j},l+1}_{i, \text{pos}} \bigg)
    \nonumber
\end{equation}

\subsection{Computational Complexity Proposition} 
\label{Appendix: sec computational complexity}
We derive the complexity of solving each local subproblem of an agent in the distributed SDP framework versus the proposed distributed framework. 
Since the SDP and SOCP constraints are significantly more computationally expensive than the linear constraints, only those are considered for the analysis.
% The computational burden of the SDP and SOCP constraints is significant compared to the linear constraints. Therefore, we consider only these constraints to determine the worst-case complexity of the algorithms. 
Further, the complexity of solving an optimization problem with SDP or SOCP constraints \cite{ben2001lectures} depends on the number of variables ($n_{\text{var}}$), the number of SDP or SOCP constraints ($n_{\text{SDP}}$ or $n_{\text{SOCP}}$), and the size of each SDP or SOCP constraint involved in the optimization problem. We define $n_c$  to represent the size of $c^{th}$ constraint. Further, we also define $n_{obs}$ and $n_{inter}$ as number of obstacles, and the maximum number of neighboring agents.
\\
\textit{Assumptions for the Analysis:} For all the agents, per time step, the state, control, and deterministic disturbance dimensions are equal and are represented by $n_{u_i}$, $n_{x_i}$, and $n_{d_i}$ respectively in the analysis. Further, we consider $\boldsymbol{\Gamma}_i = I$  for all $i \in \calV$, which gives us $\bar{n}_i = T n_{d_i} + n_{x_i}$. 
\subsubsection{Complexity of the Distributed SDP Framework}
For this framework. we assume that the shared variables among the agents are the control parameters since it is not straightforward to have the mean of the states as the shared variables in the distributed SDP approach. 
For local sub-problem of each agent $i \in \calV$, we have
\begin{enumerate}
    \item $n_{\text{var}} =
    T \big[ ( 1+ n_{\text{inter}} ) ( n_{u_i} + n_{u_i} \gamma_h n_{x_i}) + n_{obs}+ n_{\text{inter}} \big]$
    \item Number of obstacle avoidance constraints = $Tn_{obs}$, and size of each constraint $n_c^{\text{obs}} = \bar{n}_i + 1$.
    \item Number of inter-agent collision avoidance constraints = $Tn_{\text{inter}}$, and size of each constraint $n_c^{\text{inter}} = \bar{n}_i + \bar{n}_j + 1$.
\end{enumerate}
The complexity of solving the local sub-problem to $\epsilon$-solution  using the Interior point method with the SDP constraints is given as \cite{ben2001lectures}
\begin{equation}
\begin{aligned}
    O (1) \bigg( 1 + \sum_{c = 1}^{n_{\text{SDP}}} n_c  \bigg)^{1/2}
    n_{\text{var}} 
    \bigg(  n_{\text{var}}^2 
    + n_{\text{var}} \sum_{c = 1}^{n_{\text{SDP}}} n_c^2
    + \sum_{c = 1}^{n_{\text{SDP}}} n_c^3
    \bigg) \text{D}(\epsilon)
\end{aligned}
\label{Appendix: Computational SDP bound original}
\end{equation}
where D$(\epsilon)$ are the accuracy digits in the $\epsilon-$solution.

Let us now simplify each term in the above bound. Based on the constraint details mentioned earlier, we can write $\sum_{c = 1}^{n_{\text{SDP}}} n_c$ as follows 
\begin{align}
    \sum_{c= 1}^{n_{\text{SDP}}} n_{c} & =
    Tn_{\text{obs}} (\bar{n}_i + 1)
    + T n_{\text{inter}} (\bar{n}_i + \bar{n}_j + 1), \nonumber
\end{align}
which can be rewritten using the fact $\bar{n}_i = T n_{d_i} + n_{x_i}$, as 
\begin{equation}
\begin{aligned}
    \sum_{c= 1}^{n_{\text{SDP}}} n_{c} & =
    T \big[ (T n_{d_i} + n_{x_i}) (n_{\text{obs}} + 2 n_{\text{inter}})  
    + n_{\text{obs}} + n_{\text{inter}} \big] \nonumber
\end{aligned}
\end{equation}
We will simplify the contribution of each of the terms involved since we are interested in finding the complexity in terms of $O(1)$. Thus, the contribution of the term $\sum_{c= 1}^{n_{\text{SDP}}} n_{c}$ can be simplified to $T^2  n_{d_i} ( n_{\text{obs}} + n_{\text{inter}} )$.
% \begin{equation}
%     T^2  n_{d_i} ( n_{\text{obs}} + n_{\text{inter}} ).
%     \label{Appendix: Computational SDP derivation n_c eq1}
% \end{equation}
%
%
Similarly, we can simplify the contribution of the terms 
$\sum_{c= 1}^{n_{\text{SDP}}} n_{c}^2$, 
$\sum_{c= 1}^{n_{\text{SDP}}} n_{c}^3$, 
and $\sqrt{ 1 + \sum_{c = 1}^{n_{\text{SDP}}} n_c  }$ 
to 
$T^3  n_{d_i}^2 ( n_{\text{obs}} + n_{\text{inter}} )$, 
$T^4  n_{d_i}^3 ( n_{\text{obs}} + n_{\text{inter}} )$, and 
$T  \sqrt{n_{d_i} ( n_{\text{obs}} + n_{\text{inter}} )}$ respectively.

Further, using a similar intuition, we will simplify the contribution of the terms $n_{\text{var}}$ and $n_{\text{var}}^2$ to 
$T \big[ n_{u_i} \gamma_h n_{x_i} n_{\text{inter}} + n_{obs} \big]$ and $T^2 \big[ n_{u_i} \gamma_h n_{x_i} n_{\text{inter}} + n_{obs} \big]^2$ respectively.

Now, combining the above simplifications, we simplify the contribution of the term $n_{\text{var}} \big(  n_{\text{var}}^2 
    + n_{\text{var}} \sum_{c = 1}^{n_{\text{SDP}}} n_c^2
    + \sum_{c = 1}^{n_{\text{SDP}}} n_c^3
    \big)$ to the following 
\begin{equation}
\begin{aligned}
    & T^3 \big[ n_{u_i} \gamma_h n_{x_i} n_{\text{inter}} + n_{obs} \big]^3 \\
    & \quad
    +
    T^5  n_{d_i}^2 \big[ n_{u_i} \gamma_h n_{x_i} n_{\text{inter}} + n_{obs} \big]^2 ( n_{\text{obs}} + n_{\text{inter}} )
    \\
    & \qquad 
    + 
    T^5  n_{d_i}^3 \big[ n_{u_i} \gamma_h n_{x_i} n_{\text{inter}} + n_{obs} \big] ( n_{\text{obs}} + n_{\text{inter}} ) \nonumber
\end{aligned}
\end{equation}
which can be further simplified to
\begin{equation}
\begin{aligned}
    T^5 n_{d_i}^3 \big[ n_{u_i} \gamma_h n_{x_i} n_{\text{inter}} + n_{obs} \big]^3.
\end{aligned}
\label{Appendix: Computational SDP derivation semif}
\end{equation}
Consequently, we can rewrite the bound in \eqref{Appendix: Computational SDP bound original} as follows
\begin{equation}
\begin{aligned}
    O \bigg( 
    T^6 n_{d_i}^{7/2} \big[ n_{u_i} \gamma_h n_{x_i} n_{\text{inter}} + n_{obs} \big]^3 ( n_{\text{obs}} + n_{\text{inter}} )^{1/2}
    \bigg) \nonumber
\end{aligned}
\end{equation}
\subsubsection{Complexity of the Proposed Distributed Framework}
Similar to the previous derivation, we start by writing the details of the local sub-problem of each agent $i \in \calV$ as follows
\begin{enumerate}
    \item $n_{\text{var}} = T \big[ n_{u_i} + n_{u_i} \gamma_h n_{x_i} + n_{\text{pos}} 
    + n_{\text{obs}} + n_{\text{inter}} (1+ 2 n_{\text{pos}}) \big]$
    \item Number of obstacle avoidance constraints = $Tn_{obs}$, and size of each constraint $n_c^{\text{obs}} = n_{\text{pos}} + 1$.
    \item Number of inter-agent collision avoidance constraints = $Tn_{\text{inter}}$, and size of each constraint $n_c^{\text{inter}} = n_{\text{pos}} + 1$.
    \item Number of $\bmu_d^i$ bound constraints (constraint (35))= $Tn_{\text{pos}}$, and size of each constraint = $\bar{n}_i + 1$.
\end{enumerate}
The complexity of solving the local sub-problem at each agent to an $\epsilon$-solution using the Interior point method with the SOCP constraints is given as \cite{ben2001lectures}
\begin{equation}
\begin{aligned}
    O (1) ( 1 + n_{\text{SOCP}}  )^{1/2}
    n_{\text{var}} 
    \bigg(  n_{\text{var}}^2
    + n_{\text{SOCP}}
    + \sum_{c = 1}^{n_{\text{SOCP}}} n_c^2
    \bigg) \text{D}(\epsilon)
\end{aligned}
\label{Appendix: Computational SOCP bound original}
\end{equation}
where D$(\epsilon)$ are the accuracy digits in the $\epsilon-$solution.

Let us now simplify each term in the above bound. Based on the constraint details mentioned earlier, we have $n_{\text{SOCP}} = Tn_{\text{obs}} + Tn_{\text{inter}} + T n_{\text{pos}}$, whose contribution can be simplified to $T(n_{\text{obs}} + n_{\text{inter}})$.
Similarly, we can simplify the contribution of the term $\sqrt{( 1 + n_{\text{SOCP}} } $ to $\sqrt{T (n_{\text{obs}} + n_{\text{inter}})}$. 
Subsequently, we can write $\sum_{c = 1}^{n_{\text{SOCP}}} n_c^2$ as follows 
\begin{equation}
\begin{aligned}
    \sum_{c= 1}^{n_{\text{SOCP}}} n_{c}^2 &=
    (Tn_{\text{obs}} + Tn_{\text{inter}}) (n_{\text{pos}} + 1)^2 
    + T n_{\text{pos}} (\bar{n}_i + 1)^2, \nonumber
\end{aligned}
\end{equation}
which can be rewritten using the fact $\bar{n}_i = T n_{d_i} + n_{x_i}$, as 
\begin{equation}
\begin{aligned}
    \sum_{c= 1}^{n_{\text{SOCP}}} n_{c}^2 &=
    (Tn_{\text{obs}} + Tn_{\text{inter}}) (n_{\text{pos}} + 1)^2 
    \\
    &~~~~~~~~~ 
    + T n_{\text{pos}} (T n_{d_i} + n_{x_i} + 1)^2. \nonumber
\end{aligned}
\end{equation}
Thus, we can simplify the contribution of the term $\sum_{c= 1}^{n_{\text{SOCP}}} n_{c}^2$ to $T^3 n_{d_i}^2$.

Further, the contribution of the terms $n_{\text{var}}$ and $n_{\text{var}}^2$ to $T (n_{u_i} \gamma_h n_{x_i} + n_{\text{obs}} + n_{\text{inter}})$ and $T^2 (n_{u_i} \gamma_h n_{x_i} + n_{\text{obs}} + n_{\text{inter}})^2$ respectively.

For convenience, let us define $\hat{n}_{\text{com}} = n_{u_i} \gamma_h n_{x_i} + n_{\text{obs}} + n_{\text{inter}}$. 
Thereupon, we combine the above results to provide the contribution of the term $n_{\text{var}} 
    \big(  n_{\text{var}}^2
    + n_{\text{SOCP}}
    + \sum_{c = 1}^{n_{\text{SOCP}}} n_c^2
    \big)$ 
    as $T^3 \hat{n}_{\text{com}}^3
    + T^2 \hat{n}_{\text{com}} (n_{\text{obs}} + n_{\text{inter}}) 
    +T^4 \hat{n}_{\text{com}} n_{d_i}^2$, which can be further simplified to $T^3 \hat{n}_{\text{com}} [ \hat{n}_{\text{com}}^2 + n_{d_i}^2 T ]$.

Consequently, we can rewrite the bound \eqref{Appendix: Computational SOCP bound original} as follows 
\begin{equation}
\begin{aligned}
    & O \bigg( T^{7/2} ( n_{\text{obs}} + n_{\text{inter}})^{1/2}
    \hat{n}_{\text{com}}
    \big[ \hat{n}_{\text{com}}^2 + n_{d_i}^2 T \big] \bigg) \nonumber
\end{aligned}
\end{equation}
\subsubsection{Complexity of Centralized approach with proposed constraint reformulation}
We use the same analysis as in the previous case. In centralized approach, the problem is solved iteratively for $L_c$ iterations, with the non-convex constraints linearized in each iteration. The details of the problem in each iteration of this approach are as follows
\begin{itemize}
    \item $n_{\text{var}} = \sum_{i \in \calV} T \bigg[ n_{u_i} + n_{u_i} \gamma_h n_{x_i} 
    + n_{\text{pos}} 
    + n_{\text{obs}} 
    + \frac{n_{\text{inter}}}{2} \bigg]$
    \item Number of obstacle avoidance constraints = $NTn_{obs}$, and size of each constraint $n_c^{\text{obs}} $ $ = n_{\text{pos}} + 1$.
    \item Number of inter-agent collision avoidance constraints = $NTn_{\text{inter}}/2$, and size of each constraint $n_c^{\text{inter}}$ $ = n_{\text{pos}} + 1$.
    \item Number of $\bmu_d^i$ bound constraints (constraint (35))= $NTn_{\text{pos}}$, and size of each constraint = $\bar{n}_i + 1$.
\end{itemize}
We use \eqref{Appendix: Computational SOCP bound original} to derive the complexity bound. We have $n_{\text{SOCP}} = 
    NTn_{\text{obs}} + NTn_{\text{inter}}/2 + N T n_{\text{pos}}$, whose contribution can be simplified to $NT(n_{\text{obs}} + n_{\text{inter}})$. 
Similarly, we can simplify the contribution of the term $ \sqrt{ 1 + n_{\text{SOCP}} }$ to $\sqrt{NT (n_{\text{obs}} + n_{\text{inter}} )}$.
Subsequently, using the fact $\bar{n}_i = T n_{d_i} + n_{x_i}$, we can write $\sum_{c = 1}^{n_{\text{SOCP}}} n_c^2$ as follows 
\begin{equation}
\begin{aligned}
    \sum_{c= 1}^{n_{\text{SOCP}}} n_{c}^2 &=
    NTn_{\text{obs}} (n_{\text{pos}} + 1)^2
    + N Tn_{\text{inter}} (n_{\text{pos}} + 1)^2 \\
    & \qquad \qquad \qquad 
    + N T n_{\text{pos}} (T n_{d_i} + n_{x_i} + 1)^2, \nonumber
\end{aligned}
\end{equation}
and the contribution of $\sum_{c= 1}^{n_{\text{SOCP}}} n_{c}^2$ can be reduced to $N T^3 n_{d_i}^2$.

Further, the contribution of the terms $n_{\text{var}}$ and $n_{\text{var}}^2$ can be simplified to $N T [ n_{u_i} \gamma_h n_{x_i} 
    + n_{\text{obs}} 
    + n_{\text{inter}} ]$
and $N^2 T^2 [ n_{u_i} \gamma_h n_{x_i} + n_{\text{obs}} + n_{\text{inter}} ]^2$ respectively.

Let us define $\hat{n}_{\text{com}} = n_{u_i} \gamma_h n_{x_i} + n_{\text{obs}} + n_{\text{inter}}$.
Thereupon, we combine the above results to provide the contribution of the term $n_{\text{var}} 
    \big(  n_{\text{var}}^2
    + n_{\text{SOCP}}
    + \sum_{c = 1}^{n_{\text{SOCP}}} n_c^2
    \big)$ as
    $ N^2 T^2 \hat{n}_{\text{com}} [ NT \hat{n}_{\text{com}}^2 
    + n_{\text{obs}} + n_{\text{inter}}
    + T^2 n_{d_i}^2 ]$.
Since we are interested in large-scale problems, we consider that the term $N \hat{n}_{\text{com}}^2 $ is significantly larger than $T n_{d_i}^2$, thus we can further simplify the above contribution to $N^3 T^3 \hat{n}_{\text{com}}^3$. 
Consequently, the computational complexity of the centralized approach using \eqref{Appendix: Computational SOCP bound original} is given as
\begin{equation}
\begin{aligned}
    O \big( L_c (NT)^{7/2} ( n_{\text{obs}} + n_{\text{inter}})^{1/2} \hat{n}_{\text{com}}^3 \big) \nonumber
\end{aligned}
\end{equation}
%
%
%
% *********** Convergence Analysis *****************
\subsection{Convergence Analysis} \label{Appendix sec: Convergence Analysis}
For the simplicity of analysis, we will first rewrite Problem  \ref{Problem 1 distributed formulation} in a more simplified form. For that, 
We define the variable $\by^i = [\Tilde{\bmu}^i; \bu_i; \{ \tilde{c}_i^k; \tilde{c}_{ij}^k \}_{k = 1}^{T} ; vec(\vK_i)]$. 
We also write $\by = [\by^1; \by^2; \dots ; \by^N]$, 
$\bnu = [\bnu^1; \bnu^2; \dots ; \bnu^N]$, 
$\blambda = [\blambda^1; \blambda^2; \dots; \blambda^N]$.
Further, we define matrices $\vC_y^i, \vC_g^i$ such that $\Tilde{\bmu}^i = \vC_y^i \by^i$, $\tilde{\bnu}^i = \vC_g^i \bnu^i$.
Using the above notations, we can rewrite Problem \ref{Problem 1 distributed formulation} in a simplified form as follows - 
\begin{problem}[Distributed Robust Optimization Problem - Convergence Analysis] \label{Problem: DRO - compact form}
Solve the following optimization problem \\[-0.6cm]
\begin{subequations}
\begin{align}
& \min f(\by) \nonumber
\\
\mathrm{s.t.} \quad & \bg (\by) \leq 0,
\label{DRO problem - convex constraints}
\\
& h_c (\by) - \tilde{h}_c (\by) \leq 0 \quad \forall c = \llbracket 1, n_{\text{nc}} \rrbracket, 
\label{DRO problem - nonconvex constraint}
\\ &  \mathrm{\vC}_{y} \by = \mathrm{\vC}_g \bnu.
\label{DRO problem - consensus constraints}
\end{align} 
\end{subequations}
where $\mathrm{\vC}_{y} = \bdiag( \{ \vC \vC_y^i \}_{i \in \calV} )$, $\vC_g = \bdiag( \{ \vC \vC_g^i \}_{i \in \calV} )$. 
The functions $\bg \in \Rb^{n_{\text{cc}}}$, $f$, $h_c$, and $\tilde{h}_c$ are convex. Recall that the nonconvex constraints \eqref{inter-agent collision avoidance ubar} and \eqref{inter-agent collision avoidance ubar} are in the form of difference of convex functions and are represented by \eqref{DRO problem - nonconvex constraint}. In addition, note that the matrix $\mathrm{\vC}_g$ has full column rank.
\end{problem}
The points $(\bar{\by}, \bar{\bnu}, \bar{\blambda}, \bar{\boldsymbol{\vartheta}}_g, \bar{\boldsymbol{\vartheta}}_h)$ are stationary points of the above problem iff they satisfy the following Karush–Kuhn–Tucker (KKT) conditions-
\begin{align}
    & \nabla f(\bar{\by}) 
    + \mathrm{\vC}_{y} \T \bar{\blambda}
    + \rho \mathrm{\vC}_{y}\T ( \mathrm{\vC}_{y} \bar{\by} - \mathrm{\vC}_g \bar{\bnu} )
     \label{convergence: KKT gradient condition y - main} \\
    & 
    = - \nabla \bg (\bar{\by}) \bar{\boldsymbol{\vartheta}}_g 
    - \sum_{c=1}^{n_{\text{nc}}} \bar{\vartheta}_{hc} \big( \nabla h_c (\bar{\by}) 
    - \nabla \tilde{h}_c (\bar{\by}) \big) \nonumber
    \\
    &
    \mathrm{\vC}_g \T \bar{\blambda} = 0,
    \\
    & 
    \bg (\bar{\by}) \leq 0,
    ~~ h_c (\bar{\by}) - \tilde{h}_c (\bar{\by}) \leq 0 ~ \forall c \in \llbracket 1, n_{nc} \rrbracket
    \label{convergence: KKT constraint satisfaction - main}
    \\
    &
    \bar{\boldsymbol{\vartheta}}_{g} \geq 0, 
     ~~ \bar{\vartheta}_{g,c} g_c (\bar{\by}) = 0 ~\forall c \in \llbracket 1, n_{cc} \rrbracket
     \label{convergence: KKT g constraint slackness condition - main}
     \\
     & \bar{\vartheta}_{hc} \geq 0, 
     ~~ \bar{\vartheta}_{hc} ( h_c (\bar{\by}) - \tilde{h}_c (\bar{\by}) ) = 0
     ~ \forall c \in \llbracket 1, n_{nc} \rrbracket
     \label{convergence: KKT h constraint slackness condition - main}
     \\
     & \mathrm{\vC}_{y} \bar{\by} = \mathrm{\vC}_g \bar{\bnu}
      \label{convergence: KKT consensus constraint satisfaction - main}
\end{align}
We will now write the Lagrangian and the CADMM update steps as per Algorithm \ref{alg1} to solve Problem \ref{Problem: DRO - compact form} as follows
\begin{equation}
\begin{aligned}
    \calL_{\rho}^{\text{conv}} ( \by, \bnu; \blambda)
    & = 
    f(\by) + \calI_{Y_g} + \calI_{Y_h} 
    + \blambda^T ( \mathrm{\vC}_{y} \by - \mathrm{\vC}_g \bnu ) 
    \\
    &~~~~~
    + \frac{\rho}{2} \| \mathrm{\vC}_{y} \by - \mathrm{\vC}_g \bnu \|^2_2 
    \nonumber
\end{aligned}
\end{equation}
where $Y_g$ and $Y_h$ are the feasible sets of the constraints \eqref{DRO problem - convex constraints} and \eqref{DRO problem - nonconvex constraint} respectively.

\noindent
In $l+1$ CADMM iteration of Algorithm \ref{alg1}, we have

\noindent
\underline{Local update:} $\by^{l+1}$ is obtained by solving the following
\begin{align}
    & \min_{\by} ~
    f(\by) 
    + \blambda^l {}\T ( \mathrm{\vC}_{y} \by - \mathrm{\vC}_g \bnu^l ) 
    + \frac{\rho}{2} \| \mathrm{\vC}_{y} \by - \mathrm{\vC}_g \bnu^l \|^2_2
    \nonumber
    \\
    &~~~~~
    \text{s.t.  } 
    \bg (\by) \leq 0, 
    \quad 
    h_c^{\text{lin}}(\by; \by^l) \leq 0 \quad \forall c = \llbracket 1, n_{\text{nc}} \rrbracket 
    \label{convergence analysis: local sub-problem}
\end{align}
where $h_c^{\text{lin}}(\by; \by^l) = h_c (\by) - \tilde{h}_c (\by^l) - \nabla \tilde{h}_c (\by^l)\T (\by - \by^l) $.
Using the KKT conditions, we obtain the following relations
\begin{align}
    & \nabla f(\by^{l+1}) 
    + \mathrm{\vC}_{y} \T \blambda^l
    + \rho \mathrm{\vC}_{y}\T ( \mathrm{\vC}_{y} \by^{l+1} - \mathrm{\vC}_g \bnu^l )
     \label{convergence: KKT gradient condition} \\
    & 
    = - \nabla \bg (\by^{l+1}) \boldsymbol{\vartheta}^{l+1}_g 
    - \sum_{c=1}^{n_{\text{nc}}} \vartheta^{l+1}_{hc} \big( \nabla h_c (\by^{l+1}) 
    - \nabla \tilde{h}_c (\by^l) \big) \nonumber
    \\
    & 
    \bg (\by^{l+1}) \leq 0,
    ~~ h_c^{\text{lin}}(\by^{l+1}; \by^l) \leq 0 ~ \forall c \in \llbracket 1, n_{nc} \rrbracket
    \label{convergence: KKT constraint satisfaction}
    \\
    &
    \boldsymbol{\vartheta}_{g}^{l+1} \geq 0, 
     ~~ \vartheta_{g,c}^{l+1} g_c (\by^{l+1}) = 0 ~\forall c \in \llbracket 1, n_{cc} \rrbracket
     \label{convergence: KKT g constraint slackness condition}
     \\
     & \vartheta^{l+1}_{hc} \geq 0, 
     ~~ \vartheta^{l+1}_{hc} h_c^{\text{lin}}(\by^{l+1}; \by^l) = 0
     ~ \forall c \in \llbracket 1, n_{nc} \rrbracket
     \label{convergence: KKT h constraint slackness condition}
\end{align}
\underline{Global update:} $\bnu^{l+1}$ is obtained by solving the following
\begin{align}
    & \min_{\bnu} ~~ \blambda^l {}\T ( \mathrm{\vC}_{y} \by^{l+1} - \mathrm{\vC}_g \bnu ) 
    + \frac{\rho}{2} \| \mathrm{\vC}_{y} \by^{l+1} - \mathrm{\vC}_g \bnu \|^2_2 \nonumber
\end{align}
A closed form solution can be derived for the above. Using first order optimality condition, we get
\begin{equation}
    -\mathrm{\vC}_g\T \big( \blambda^l
    + \rho ( \mathrm{\vC}_{y} \by^{l+1} - \mathrm{\vC}_g \bnu^{l+1} ) \big) = 0
    \label{conv analysis: global update derivative eq1}
\end{equation}
which implies $\mathrm{\vC}_g\T ( \blambda^{l+1} + \delta \blambda^l ) = 0$. Thus, if $\mathrm{\vC}_g\T \blambda^0 = 0$ then $\mathrm{\vC}_g\T \blambda^l = 0$ for all CADMM iterations $l$. Substituting this in \eqref{conv analysis: global update derivative eq1}, we get the global update as follows
\begin{equation}
    \bnu^{l+1} = ( \mathrm{\vC}_g\T\mathrm{\vC}_g  )^{-1}
    \mathrm{\vC}_g\T \mathrm{\vC}_{y} \by^{l+1} 
    ~~ (= \text{Prj}(\mathrm{\vC}_g) \mathrm{\vC}_{y} \by^{l+1} ),
    \label{convergence analysis: global update}
\end{equation}
which is nothing but projection of $\mathrm{\vC}_{y} \by^{l+1}$ onto $Im(\mathrm{\vC}_g\T)$.

\noindent
\underline{Discounted dual update:}
\begin{equation}
    \blambda^{l+1} = (1-\delta) \blambda^l 
    + \rho ( \mathrm{\vC}_{y} \by^{l+1} - \mathrm{\vC}_g \bnu^{l+1} ) 
    \label{convergence analysis: dual update}
\end{equation}
\subsubsection{Sketch of the Analysis}
The Analysis involves the following steps -
\begin{enumerate}
    \item We consider a regularized Lagrangian function \cite{yang2022proximal}, and derive an upper bound of the change in this function value over each iteration. We show that this bound is dependent on the change in the dual variables. 
    \item The next step is to upper bound the change in the dual variables in terms of the local and global variables. This is a crucial step that provides a detailed understanding of the algorithm as well as highlights the importance of using a discounted dual update step.
    \item By combining the results from the above steps, we construct our merit function and prove convergence under empirical assumption.
\end{enumerate}
Note: Throughout this section, we use for any vector $\bp$, the notation $[\Delta \bp]_l = \bp^{l+1} - \bp^{l}$. 
\subsubsection{Boundedness of Dual Variable} 
Let us consider a closed set $Y = \{ (\by, \bnu) \; | \; \by \in Y_g, \by \in Y_h, \bnu = ( \mathrm{\vC}_g\T\mathrm{\vC}_g  )^{-1}
\mathrm{\vC}_g\T \mathrm{\vC}_{y} \by \}$. Since $h_c$ and $\tilde{h}_c$ are convex, $h_c^{\text{lin}}(\by; \by^l) \leq 0$ implies that $\by \in Y_h$. Thus, we have $(\by^{l}, \bnu^{l}) \in Y$ for all $l$. Using Proposition 4 of \cite{yang2022proximal}, we can prove that $\| \blambda^l \|_2$ is bounded for all $l$.
\subsubsection{Regularized Lagrangian function}
A regularized Lagrangian function  \cite{yang2022proximal} for the Problem \ref{Problem: DRO - compact form} is defined as follows
\begin{equation}
    \calL_{\rho}^{\text{reg},l} 
    := \calL_{\rho}^{\text{reg}} ( \by^l, \bnu^l, \blambda^l)
    = 
    \calL_{\rho}^{\text{conv}} ( \by^l, \bnu^l, \blambda^l)
    - \frac{\delta}{2 \rho} \| \blambda^l \|^2_2 \nonumber
\end{equation}
%
% \paragraph{Lower boundedness of $\calL_{\rho}^{\text{reg},l}$}
% Let us consider a closed set $Y = \{ (\by, \bnu) \; | \; \by \in Y_g, \by \in Y_h, \bnu = ( \mathrm{\vC}_g\T\mathrm{\vC}_g  )^{-1} \mathrm{\vC}_g\T \mathrm{\vC}_{y} \by \}$. Since $h_c$ and $\tilde{h}_c$ are convex, $h_c^{\text{lin}}(\by; \by^l) \leq 0$ implies that $\by \in Y_h$. Thus, we have $(\by^{l}, \bnu^{l}) \in Y$ for all $l$. From Proposition 4 of \cite{yang2022proximal}, we can prove that $\| \blambda^l \|_2$ is bounded for all $l$. Further, 
Using Proposition 5 of \cite{yang2022proximal} and the fact that $f$ is quadratic (lower-bounded), we have that $\calL_{\rho}^{\text{reg}, l}$ is bounded for all $l$. We derive an upper bound of $[\Delta\calL_{\rho}^{\text{reg}}]_l$ in the following three steps 

\noindent 
\underline{Change due to local update:} We have
\begin{align}
    \Delta_y^{l} \calL_{\rho}^{\text{reg}} & :=  \calL_{\rho}^{\text{reg}} ( \by^{l+1}, \bnu^{l}; \blambda^{l}) 
    - \calL_{\rho}^{\text{reg}} ( \by^l, \bnu^l; \blambda^l) \nonumber
    \\
    & =
    f( \by^{l+1} ) - f( \by^l )
    + \blambda^l {}\T \mathrm{\vC}_{y} (\by^{l+1} - \by^{l})
    \label{Appendix: DRO Lreg change LU eq1} \\
    & \quad 
    + \frac{\rho}{2} 
    \big( \| \mathrm{\vC}_{y} \by^{l+1} - \mathrm{\vC}_g \bnu^l \|^2_2 
    - \| \mathrm{\vC}_{y} \by^{l} - \mathrm{\vC}_g \bnu^l \|^2_2 \big) 
    \nonumber
\end{align}
Let us now simplify the terms on the RHS of the above equation. Since $f$ is convex, we have 
\begin{equation}
\begin{aligned}
    & f(\by^{l+1}) - f(\by^l) \leq \nabla f(\by^{l+1})^T ( \by^{l+1} - \by^l)
\end{aligned}
\label{Appendix: DRO Lreg change LU eq2}
\end{equation}
Next, we can simplify the following 
\begin{equation}
\begin{aligned}
    & \| \mathrm{\vC}_{y} \by^{l+1} - \mathrm{\vC}_g \bnu^l \|^2_2 
    - \| \mathrm{\vC}_{y} \by^{l} - \mathrm{\vC}_g \bnu^l \|^2_2
    \\
    % &~ =
    % \big( \mathrm{\vC}_{y} \by^{l+1} + \mathrm{\vC}_{y} \by^{l} 
    % - 2 \mathrm{\vC}_g \bnu^l \big)\T 
    % \mathrm{\vC}_{y} ( \by^{l+1} - \by^{l} )
    % \\
    &~ = 
    % 2( \mathrm{\vC}_{y} \by^{l+1} - \mathrm{\vC}_g \bnu^l )\T 
    % \mathrm{\vC}_{y} ( \by^{l+1} - \by^{l} )
    2( \mathrm{\vC}_{y} \by^{l+1} - \mathrm{\vC}_g \bnu^l )\T 
    \mathrm{\vC}_{y} ( [\Delta \by]_{l} )
    - \| \mathrm{\vC}_{y} [\Delta \by]_{l} \|_2^2
    % \\
    % &~~~~~~~~~~~~~~~~~~~
    % - \| \mathrm{\vC}_{y} ( \by^{l+1} - \by^{l} ) \|_2^2
\end{aligned}   
\label{Appendix: DRO Lreg change LU eq3}
\end{equation}
Substituting \eqref{Appendix: DRO Lreg change LU eq2} and \eqref{Appendix: DRO Lreg change LU eq3} in \eqref{Appendix: DRO Lreg change LU eq1}, we get
\begin{equation}
\begin{aligned}
    \Delta_y^{l} \calL_{\rho}^{\text{reg}}
    & \leq 
    % - \frac{\rho}{2} \| \mathrm{\vC}_{y} ( \by^{l+1} - \by^{l} ) \|_2^2
    - \frac{\rho}{2} \| \mathrm{\vC}_{y}  [\Delta \by]_{l}  \|_2^2
    + \big( \nabla f(\by^{l+1}) + \vC_y \T \blambda^l 
    \\
    &~~~
    + \rho \mathrm{\vC}_{y} \T ( \mathrm{\vC}_{y} \by^{l+1} - \mathrm{\vC}_g \bnu^l ) \big)\T (  [\Delta \by]_{l} ),
    \nonumber
\end{aligned}
\end{equation}
which can be rewritten using \eqref{convergence: KKT gradient condition} as 
\begin{align}
    \Delta_y^{l} \calL_{\rho}^{\text{reg}}
    & \leq 
    % - \frac{\rho}{2} \| \mathrm{\vC}_{y} ( \by^{l+1} - \by^{l} ) \|_2^2
    - \frac{\rho}{2} \| \mathrm{\vC}_{y} [\Delta \by]_{l}  \|_2^2
    - \bigg( \nabla \bg (\by^{l+1}) \boldsymbol{\vartheta}^{l+1}_g  
    \label{Appendix: DRO Lreg change LU eq4}
    \\
    &~~  
    + \sum_{c=1}^{n_{\text{nc}}} \vartheta^{l+1}_{hc} \big( \nabla h_c (\by^{l+1}) 
    - \nabla \tilde{h}_c (\by^l) \big) \bigg)\T
    ( [\Delta \by]_{l} ) \nonumber
\end{align}
Since $\bg$ is convex and $\boldsymbol{\vartheta}^{l+1}_g \geq 0$, we have 
\begin{equation}
    - \boldsymbol{\vartheta}^{l+1}_g {}\T \nabla \bg (\by^{l+1})\T (  \by^{l+1} - \by^{l} )
     \leq 
     \boldsymbol{\vartheta}^g{}\T \bg(y^{l})
     - \boldsymbol{\vartheta}^g{}\T \bg(y^{l+1}),
    \nonumber
\end{equation}
which, using \eqref{convergence: KKT g constraint slackness condition}, would give
\begin{align}
    - \boldsymbol{\vartheta}^{l+1}_g  {}\T \nabla \bg (\by^{l+1})\T (  \by^{l+1} - \by^{l} )
     & \leq 0.
     \label{Appendix: DRO Lreg change LU eq5}
\end{align}
Next, using \eqref{convergence: KKT h constraint slackness condition}, we have
\begin{align}
    & \vartheta^{l+1}_{hc} \big( \nabla h_c (\by^{l+1}) 
    - \nabla \tilde{h}_c (\by^l) \big)\T (  \by^{l+1} - \by^{l} )
    \nonumber \\
    & =
    \vartheta^{l+1}_{hc} \big( \nabla h_c (\by^{l+1})\T (  \by^{l+1} - \by^{l} ) 
    - h_c (\by^{l+1}) 
    + \tilde{h}_c ( \by^l ) \big) \nonumber
    \\
    & 
    \geq - \vartheta^{l+1}_{hc} \big( h_c (\by^l)
    - \tilde{h}_c ( \by^l ) \big)
    \label{Appendix: DRO Lreg change LU eq6}
\end{align}
For the last inequality above, we used the fact that $h_c$ is convex and $\vartheta^{l+1}_{hc} \geq 0$. Substituting \eqref{Appendix: DRO Lreg change LU eq5} and \eqref{Appendix: DRO Lreg change LU eq6} into \eqref{Appendix: DRO Lreg change LU eq4}, we get
\begin{equation}
    \Delta_y^{l} \calL_{\rho}^{\text{reg}}
    \leq r_h^l
    - \frac{\rho}{2} \| \mathrm{\vC}_{y} [\Delta \by]_{l} \|_2^2
    \label{Convergence DRO lreg change local}
\end{equation}
where $r_h^l = \sum_{c=1}^{n_{\text{nc}}} \vartheta^{l+1}_{hc} ( h_c (\by^l) - \tilde{h}_c ( \by^l ) )$.

\noindent
\underline{Change due to global update:} We have
\begin{align}
    \Delta_{\nu}^l \calL_{\rho}^{\text{reg}} & :=  \calL_{\rho}^{\text{reg}} ( \by^{l+1}, \bnu^{l+1}; \blambda^{l}) 
    - \calL_{\rho}^{\text{reg}} ( \by^{l+1}, \bnu^l; \blambda^l) \nonumber
    \\
    & =
    - \blambda^l {}\T \mathrm{\vC}_g ( \bnu^{l+1} - \bnu^l )
    + \frac{\rho}{2} 
    \big( \| \mathrm{\vC}_{y} \by^{l+1} - \mathrm{\vC}_g \bnu^{l+1} \|^2_2 
    \nonumber \\
    & \qquad \quad
    - \| \mathrm{\vC}_{y} \by^{l+1} - \mathrm{\vC}_g \bnu^l \|^2_2 \big) 
    \label{Appendix: DRO Lreg change LU global eq1}
\end{align}
Similar to \eqref{Appendix: DRO Lreg change LU eq3}, we can write
\begin{align}
    & \| \mathrm{\vC}_{y} \by^{l+1} - \mathrm{\vC}_g \bnu^{l+1} \|^2_2 
    - \| \mathrm{\vC}_{y} \by^{l+1} - \mathrm{\vC}_g \bnu^l \|^2_2
    \\
    &~~ = 
    % - \| \vC_g ( \bnu^{l+1} - \bnu^{l} ) \|_2^2 
    %  \\
    % &~~~~~~~~ - 2(  \mathrm{\vC}_{y} \by^{l+1} - \mathrm{\vC}_g \bnu^{l+1}  )\T
    % \vC_g ( \bnu^{l+1} - \bnu^{l} ) \nonumber
    - \| \vC_g [\Delta \bnu]_{l} ) \|_2^2 
    - 2(  \mathrm{\vC}_{y} \by^{l+1} - \mathrm{\vC}_g \bnu^{l+1}  )\T
    \vC_g [\Delta \bnu]_{l} \nonumber
\end{align}
Further, using \eqref{convergence analysis: dual update} and the fact that $\vC_g\T \blambda^l = 0$, we get $2(  \mathrm{\vC}_{y} \by^{l+1} - \mathrm{\vC}_g \bnu^{l+1}  )\T
    \vC_g [\Delta \bnu]_{l} = 0$. 
Thus, using the above result, we can rewrite \eqref{Appendix: DRO Lreg change LU global eq1} as follows
\begin{equation}
    \Delta_{\nu}^l \calL_{\rho}^{\text{reg}}
    = - \frac{\rho}{2} \| \vC_g [\Delta \bnu]_{l} \|_2^2
    \label{Convergence DRO lreg change global}
\end{equation}
\underline{Change due to dual update:} From the Proposition 2 of \cite{yang2022proximal}, we can derive
\begin{align}
    \Delta_{\lambda}^l \calL_{\rho}^{\text{reg}} & :=
    \calL_{\rho}^{\text{reg}} ( \by^{l+1}, \bnu^{l+1}; \blambda^{l+1}) 
    - \calL_{\rho}^{\text{reg}} ( \by^{l+1}, \bnu^{l+1}; \blambda^l) \nonumber
    \\
    & =
    \frac{2- \delta}{2 \rho} \| [\Delta\blambda]_{l} \|_2^2 
    \label{Convergence DRO lreg change dual}
\end{align}

Thus, combining \eqref{Convergence DRO lreg change local}, \eqref{Convergence DRO lreg change global} and \ref{Convergence DRO lreg change dual}, we obtain an upper bound on change in $\calL_{\rho}^{\text{reg}}$ in $l+1$ iteration as follows
\begin{equation}
\begin{aligned}
    \Delta^l \calL_{\rho}^{\text{reg}} 
   & 
   \leq
    - \frac{\rho}{2} \| \mathrm{\vC}_{y} [\Delta \by]_{l} \|_2^2
    - \frac{\rho}{2} \| \vC_g [ \Delta \bnu]_{l} \|_2^2
    \\
    &~~~~~~~~
    + r^l_h
    + \frac{2- \delta}{2 \rho} \| [\Delta \blambda]_{l} \|_2^2 
    \label{convergence change in Lreg bound}
\end{aligned}
\end{equation}
%
% We will prove that the term $r^l_h$ is non-positive in the subsequent sections. 
% Thus, the only nonnegative term in the RHS of the above bound involves $\| [\Delta \blambda]_{l} \|_2^2$. 
% Therefore, we will now derive a bound of the term in terms of the local and global update variables.
 The last term on the RHS involving $\| [\Delta \blambda]_{l} \|_2^2$ is non-positive, thus we will now derive its bound in terms of the local and global update variables.
\subsubsection{Bound on dual variables update}
We will derive a bound on $\| [\Delta \blambda]_{l} \|_2^2$ by leveraging the problem structure. We derive two intermediate relations based on the KKT conditions of the local sub-problems in $(l+1)^{th}$ and $l^{th}$ iterations, and combine them to obtain the bound.

% In this subsection, we will derive a bound on the dual variables by leveraging the structure of the problem and the updates. We derive two intermediate relations based on the KKT conditions of the local sub-problems in $(l+1)^{th}$ and $l^{th}$ iterations. We then bound the dual variables by combining these relations.
% ~~~~~~~~~ Relation-1~~~~~~~~~~~~~
\noindent
\underline{Relation-1:} Multiplying KKT condition \eqref{convergence: KKT gradient condition} of $'l+1'$ subproblem with $([\Delta \by]_{l})\T$, and using \eqref{Appendix: DRO Lreg change LU eq5} and \eqref{Appendix: DRO Lreg change LU eq6}, we get
\begin{align}
    & ([\Delta \by]_{l})\T \big[ \nabla f(\by^{l+1}) 
    + \mathrm{\vC}_{y} \T \blambda^l
    + \rho \mathrm{\vC}_{y}\T ( \mathrm{\vC}_{y} \by^{l+1} - \mathrm{\vC}_g \bnu^l ) \big] \nonumber
    \\
    &~~~~ 
    \leq \sum_{c=1}^{n_{\text{nc}}} \vartheta^{l+1}_{hc} \big( h_c (\by^l) - \tilde{h}_c ( \by^l ) \big) \; \; ( = r_h^l )
    \label{Dual update bound: R1 eq1}
\end{align}
%
% Let us now further simplify the LHS of the above.
From \eqref{convergence analysis: global update} and \eqref{convergence analysis: dual update}, we can obtain
\begin{align}
    & \begin{aligned}
    & \blambda^l
    + \rho ( \mathrm{\vC}_{y} \by^{l+1} - \mathrm{\vC}_g \bnu^l )
    % \\
    % & = 
    % \blambda^l
    % + \rho ( \mathrm{\vC}_{y} \by^{l+1} - \mathrm{\vC}_g \bnu^{l+1} )
    % + \rho \mathrm{\vC}_g ( [\Delta \bnu]_{l} ) 
    % \nonumber
    % \\
    % &~~~~~~~~~~~~~~~~~~~~ 
    =
    \blambda^{l+1} + \delta \blambda^l
    + \rho \mathrm{\vC}_g ( [\Delta \bnu]_{l} ), \nonumber
    % \qquad (\because \eqref{convergence analysis: dual update})
    \end{aligned}
    % \label{Dual update bound: R1 eq2}  
    \\
    & 
    \mathrm{\vC}_g\T \mathrm{\vC}_{y} ( [ \Delta \by]_{l} )
    = ( \mathrm{\vC}_g\T \mathrm{\vC}_g  ) ( [\Delta \bnu]_{l}  ) \nonumber
    % \label{Dual update bound: R1 eq3}
\end{align}
% Further, from \eqref{convergence analysis: global update}, we get 
% \begin{equation}
%     \mathrm{\vC}_g\T \mathrm{\vC}_{y} ( [ \Delta \by]_{l} )
%     = ( \mathrm{\vC}_g\T \mathrm{\vC}_g  ) ( [\Delta \bnu]_{l}  )
%     \label{Dual update bound: R1 eq3}
% \end{equation}
%
Substituting the above into \eqref{Dual update bound: R1 eq1}, we get
\begin{align}
    & ( [ \Delta \by]_{l} )\T \big[ \nabla f(\by^{l+1}) 
    + \mathrm{\vC}_{y} \T ( \blambda^{l+1} + \delta \blambda^l ) \big] 
    \label{Dual update bound: R1 final} 
    \\
    &~~~
    \leq 
    - \rho \| \mathrm{\vC}_g [\Delta \bnu]_{l} \|^2_2
    + r^l_h \nonumber
\end{align}
%
% ~~~~~~~~~~~ Relation-2~~~~~~~~~~~~
\underline{Relation-2:} Multiplying the KKT gradient condition \eqref{convergence: KKT gradient condition} for $l$ subproblem with $-([\Delta \by]_{l})\T$, we get
\begin{align}
    & - ([\Delta \by]_{l})\T \big[ \nabla f(\by^{l}) 
    + \mathrm{\vC}_{y} \T \blambda^{l-1}
    + \rho \mathrm{\vC}_{y}\T ( \mathrm{\vC}_{y} \by^{l} - \mathrm{\vC}_g \bnu^{l-1} ) \big] \nonumber
    \\
    &~~
    = 
    \sum_{c=1}^{n_{\text{nc}}} \vartheta^{l}_{hc} ([\Delta \by]_{l})\T \big( \nabla h_c (\by^{l}) 
    - \nabla \tilde{h}_c (\by^{l-1}) \big) 
    \label{Dual update bound: R2 eq1}
    \\[-0.1cm]
    &~~~~~~~~~~~~~
     + ([\Delta \by]_{l})\T \nabla \bg (\by^{l}) \boldsymbol{\vartheta}^{l}_g 
     \nonumber
\end{align}
Let us now simplify the terms in the above equation. Using a similar intuition as in \eqref{Appendix: DRO Lreg change LU eq5}, we can derive the following
\begin{equation}
    ([\Delta \by]_{l})\T \nabla \bg (\by^{l}) \boldsymbol{\vartheta}^{l}_g  \leq 0.
    \label{Dual update bound: R2 eq2}
\end{equation}
Since $h_c$ is convex, and using \eqref{convergence: KKT h constraint slackness condition}, we can write
\begin{equation}
    \vartheta^{l}_{hc} ([\Delta \by]_{l})\T \nabla h_c (\by^{l}) 
    \leq \vartheta^{l}_{hc} ( h_c (\by^{l+1}) - h_c (\by^{l})).
    \label{Dual update bound: R2 eq4}
\end{equation}
Next, we can write
\begin{align}
    & - \vartheta^{l}_{hc} ([\Delta \by]_{l})\T \nabla \tilde{h}_c (\by^{l-1})  \nonumber
    \\
    &
    =
   \vartheta^{l}_{hc} (\by^{l} - \by^{l-1} + \by^{l-1} - \by^{l+1})\T \nabla \tilde{h}_c (\by^{l-1})  
   \label{Dual update bound: R2 eq3}
   \\
   &
   = 
   \vartheta^{l}_{hc} ( h_c( \by^l ) - \tilde{h}_c(\by^{l-1})
   - \nabla \tilde{h}_c (\by^{l-1})\T (\by^{l+1} - \by^{l-1}) )
   \nonumber
\end{align}
Above, the last equality is obtained using \eqref{convergence: KKT h constraint slackness condition}. 
Subsequently, using similar intuition as in 'Relation-1', we can obtain
\begin{align}
    & -([\Delta \by]_{l})\T \big[ \mathrm{\vC}_{y} \T \blambda^{l-1}
    + \rho \mathrm{\vC}_{y}\T ( \mathrm{\vC}_{y} \by^{l} - \mathrm{\vC}_g \bnu^{l-1} ) \big]
    \label{Dual update bound: R2 eq5} \\
    & = - ( [ \Delta \by]_{l} )\T  \mathrm{\vC}_{y} \T ( \blambda^{l-1}
    + \delta \blambda^l )
    % \label{Dual update bound: R2 eq5}
    % \\
    % &~~~~~~~~~~~
    + \rho ( [\Delta \bnu]_{l} ) \mathrm{\vC}_g\T \mathrm{\vC}_g ( [\Delta \bnu]_{l-1} ) \nonumber
\end{align}
Using the fact that $a \T b \leq ( \|a\|^2_2 + \|b\|_2^2 )/2 $, we can write
\begin{equation}
\begin{aligned}
    & -( [\Delta \bnu]_{l} ) \mathrm{\vC}_g\T \mathrm{\vC}_g ( [\Delta \bnu]_{l-1} )
    \\
    &~~~~~~~~~~~ 
    \leq
    \frac{1}{2} \big( 
    \| \mathrm{\vC}_g  [\Delta \bnu]_{l-1} \|^2_2
    + \| \mathrm{\vC}_g  [\Delta \bnu]_{l} \|^2_2 \big)
\end{aligned}
\label{Dual update bound: R2 eq6}
\end{equation}
Substituting \eqref{Dual update bound: R2 eq2}-\eqref{Dual update bound: R2 eq6} into \eqref{Dual update bound: R2 eq1}, we get
\begin{align}
    & ( [\Delta \by]_{l} )\T \big[ \nabla f(\by^{l}) 
    + \mathrm{\vC}_{y} \T ( \blambda^{l-1}
    + \delta \blambda^l ) \big] 
    \label{Dual update bound: R2 final}
    \\
    &~~~~~~
    \leq
    \frac{\rho}{2} \big( 
    \| \mathrm{\vC}_g [\Delta \bnu]_{l-1} \|^2_2
    + \| \mathrm{\vC}_g [\Delta \bnu]_{l} \|^2_2 \big)
    + \hat{r}^l_h \nonumber
\end{align}
where $\hat{r}^l_h = \sum_{c=1}^{n_{\text{nc}}} \vartheta^{l}_{hc} h_c^{\text{lin}}(\by^{l+1}; \by^{l-1}) $.

%%% ~~~~~~~~~~~~~~ final relation ~~~~~~~~~~~~~~~
\noindent
\underline{Bound on $\| [\Delta \blambda]_{l} \|_2^2$:}
% We derive this by combining \eqref{Dual update bound: R1 final} and \eqref{Dual update bound: R2 final}. Let us now simplify some of the terms that would arise from this combination. 
% Sinnce $f$ is convex, we have
% \begin{equation}
%     (\by^{l} - \by^{l+1})\T ( \nabla f(\by^{l}) - \nabla f(\by^{l+1})) \geq 0.
%     \label{Dual update bound: R3 eq1}
% \end{equation}
%
By combining \eqref{Dual update bound: R1 final} and \eqref{Dual update bound: R2 final}, and using \eqref{Appendix: DRO Lreg change LU eq2}, we get
\begin{equation}
    ([\Delta \by]_{l})\T \mathrm{\vC}_{y} \T 
    \big( [\Delta\blambda]_{l} 
    + \delta  [\Delta\blambda]_{l-1} \big)
    \leq r_{c}^l - r_{\nu}^l
    \label{Dual update bound: R3 eq2} 
\end{equation}
% \begin{align}
%      & (\by^{l+1} - \by^{l})\T \mathrm{\vC}_{y} \T 
%     \big[ \blambda^{l+1} - \blambda^{l} 
%     + \delta (  \blambda^{l} - \blambda^{l-1}   ) \big]
%     \leq r_{c} - r_{\nu}
%     \label{Dual update bound: R3 eq2} 
%     % &~~~
%     % \leq 
%     % \frac{\rho}{2} \big( 
%     % \| \mathrm{\vC}_g ( \bnu^{l} - \bnu^{l-1} ) \|^2_2
%     % - \| \mathrm{\vC}_g ( \bnu^{l+1} - \bnu^{l} ) \|^2_2 \big)
%     % + r_{conv} \nonumber
% \end{align}
where $r_{\nu}^l = \rho ( \| \mathrm{\vC}_g [\Delta \bnu]_{l} \|^2_2
- \| \mathrm{\vC}_g [\Delta \bnu]_{l-1} ) \|^2_2 )/2 $, 
$r_{c}^l = r_h^l + \hat{r}_h^l$.
% Let us now further simplify terms on the LHS of the above equation. 
Using \eqref{convergence analysis: dual update} and the fact that $\mathrm{\vC}_g\T \blambda^l = 0$ for all $l$, we derive 
% \begin{align}
%     & (\by^{l+1} - \by^{l})\T \mathrm{\vC}_{y} \T 
%     \big[ \blambda^{l+1} - \blambda^{l} 
%     + \delta (  \blambda^{l} - \blambda^{l-1}   ) \big]
%     \nonumber
%     \\
%     &~ =
%     \big( \mathrm{\vC}_{y} \by^{l+1} - \mathrm{\vC}_g \bnu^{l+1} 
%     - ( \mathrm{\vC}_{y} \by^{l} - \mathrm{\vC}_g \bnu^{l+1} )   \big)\T
%    \big[ \blambda^{l+1} \nonumber
%    \\
%    &~~~~~~~~~~~~~~ 
%    - \blambda^{l} 
%     + \delta (  \blambda^{l} - \blambda^{l-1}   ) \big] \nonumber
%     \\
%     & =
%     \frac{1}{\rho} \big[ \| \blambda^{l+1} - \blambda^{l} \|_2^2
%     - (1- \delta) (\blambda^l - \blambda^{l-1})\T ( \blambda^{l+1} - \blambda^{l} ) \nonumber
%     \\
%     &~~~~
%     + \delta ( \blambda^{l+1} - \blambda^{l} )\T (  \blambda^{l} - \blambda^{l-1}  )
%     - \delta (1-\delta) \| \blambda^{l} - \blambda^{l-1} \|_2^2
%     \big] \nonumber
%     \\
%     & 
%     = \frac{1}{\rho} \big[ \| \blambda^{l+1} - \blambda^{l} \|_2^2
%     - (1- 2 \delta) (\blambda^l - \blambda^{l-1})\T ( \blambda^{l+1} - \blambda^{l} ) \nonumber
%     \\
%     &~~~~~~~~~~~~~
%     - \delta (1-\delta) \| \blambda^{l} - \blambda^{l-1} \|_2^2
%     \big] \nonumber
% \end{align}
%
\begin{align}
    & (\by^{l+1} - \by^{l})\T \mathrm{\vC}_{y} \T 
    \big( [\Delta\blambda]_{l} 
    + \delta  [\Delta\blambda]_{l-1} \big)
    \nonumber
    \\
    &~ =
    \big( \mathrm{\vC}_{y} \by^{l+1} - \mathrm{\vC}_g \bnu^{l+1} 
    - ( \mathrm{\vC}_{y} \by^{l} - \mathrm{\vC}_g \bnu^{l+1} )   \big)\T
   \big( [\Delta\blambda]_{l} 
   \nonumber
   \\ &~~~~~~~~~~~~~~~~~~~~~~~~
    + \delta  [\Delta\blambda]_{l-1} \big) \nonumber
    % \\
    % & =
    % \frac{1}{\rho} \big( \| [\Delta\blambda]_{l+1} \|_2^2
    % - (1- \delta) [\Delta\blambda]_{l}\T [\Delta\blambda]_{l+1}  \nonumber
    % \\
    % &~~~~~~~~~
    % + \delta [\Delta\blambda]_{l+1}\T [\Delta\blambda]_{l}
    % - \delta (1-\delta) \| [\Delta\blambda]_{l} \|_2^2
    % \big) \nonumber
    \\
    & 
    = \frac{1}{\rho} \big[ \| [\Delta\blambda]_{l+1} \|_2^2
    - (1- 2 \delta) [\Delta\blambda]_{l}\T [\Delta\blambda]_{l+1} 
    \nonumber \\[-0.1cm]
    &~~~~~~~~~~~~~~~~~~~~~~~~~~~~~~~~~~~~~
    - \delta (1-\delta) \| [\Delta\blambda]_{l} \|_2^2
    \big] \nonumber
\end{align}
Given that $\delta \in (0,0.5]$, and using the fact that $-a^Tb \geq - (\|a\|_2^2 + \|b\|^2_2)/2 $, we can rewrite the above as
\begin{align}
    & ([\Delta \by]_{l})\T \mathrm{\vC}_{y} \T 
    \big( [\Delta\blambda]_{l} 
    + \delta  [\Delta\blambda]_{l-1} \big)
    \nonumber
    \\
    &~ \geq
    \frac{1}{\rho} \bigg[ \| [\Delta\blambda]_{l} \|_2^2
    - \frac{(1- 2 \delta)}{2} \| [\Delta\blambda]_{l} \|_2^2 
    \nonumber \\[-0.1cm]
    &~~~~~~~
    - \frac{(1- 2 \delta)}{2}\| [\Delta\blambda]_{l-1} \|_2^2 
    - \delta (1-\delta) \| [\Delta\blambda]_{l-1} \|_2^2
    \bigg] \nonumber
    \\
    &~~ =
    \frac{1+ 2 \delta}{2 \rho} \| [\Delta\blambda]_{l} \|_2^2
    - \frac{1- 2 \delta^2}{2 \rho} \| [\Delta\blambda]_{l-1} \|_2^2
    \nonumber
    \\
    &~~ =
    % \frac{1- 2 \delta^2}{2 \rho} \big(\| [\Delta\blambda]_{l} \|_2^2
    %  - \| [\Delta\blambda]_{l-1} \|_2^2 \big)
    %  \nonumber \\[-0.1cm]
    %  &~~~~~~~~~~
    r_{\lambda}^l
     + \frac{\delta(1+ \delta)}{\rho} \| [\Delta\blambda]_{l} \|_2^2 \nonumber
\end{align}
where $r_{\lambda}^l = (1- 2 \delta^2) ( \| [\Delta\blambda]_{l} \|_2^2
     - \| [\Delta\blambda]_{l-1} \|_2^2 )/( 2 \rho )$.
Substituting the above result in \eqref{Dual update bound: R3 eq2} would give us 
\begin{align}
    \frac{\delta(1+ \delta)}{\rho} \| [\Delta\blambda]_{l} \|_2^2
    +  r_{\lambda}^l + r_{\nu}^l \leq r_{c}^l 
    \label{convergence: dual update upper bound}
\end{align}

\subsubsection{Construction of Merit Function}
The upper bound obtained in \eqref{convergence: dual update upper bound} can then be used to upper bound the term $\Delta^l \calL_{\rho}^{\text{reg}}$ from \eqref{convergence change in Lreg bound} as follows
\begin{align}
    & \Delta^l \calL_{\rho}^{\text{reg}} 
    + \alpha_{\delta} ( r_{\lambda}^l + r_{\nu}^l )
    + \hat{\alpha}_{\delta} \| [\Delta \blambda]_{l} \|_2^2
    \label{Convergence analysis descent relation}
    \\
    &~~~
    \leq \alpha_{\delta} r_{c}^l + r_h^l
    - \frac{\rho}{2} \| \mathrm{\vC}_{y} [\Delta \by]_{l} \|_2^2
    - \frac{\rho}{2} \| \vC_g [ \Delta \bnu]_{l} \|_2^2 
    \nonumber
\end{align}
where $\hat{\alpha}_{\delta} = ( 2\delta(1+ \delta)  \alpha_{\delta} - (2- \delta) )/ (2 \rho) $ and $\alpha_{\delta}$ such that $\hat{\alpha}_{\delta} >0$. We consider the following merit function
\begin{align}
    \tilde{\calL}^{l+1} &:= \tilde{\calL} (\by^{l+1}, \bnu^{l+1}, \blambda^{l+1}; \bnu^{l}, \blambda^{l}) \nonumber
     \\
    & =
    \calL_{\rho}^{\text{reg}, l+1} 
    +  \frac{\alpha_{\delta} \rho}{2} \| \mathrm{\vC}_g [\Delta \bnu]_{l} \|^2_2
    + \frac{\alpha_{\delta} (1- 2 \delta^2) }{2 \rho } \| [\Delta\blambda]_{l} \|_2^2 \nonumber
\end{align}
such that $ \tilde{\calL}^{l+1} - \tilde{\calL}^{l} = \Delta^l \calL_{\rho}^{\text{reg}} 
    + \alpha_{\delta} ( r_{\lambda}^l + r_{\nu}^l )$. Since $\calL_{\rho}^{\text{reg}, l+1} $ is lower bounded, the function $\tilde{\calL}^{l+1}$ is lower bounded.
\subsubsection{Convergence Analysis}
We can rewrite \eqref{Convergence analysis descent relation} in terms of the merit function $\tilde{\calL}$ as follows
\begin{equation}
\begin{aligned}
    [\Delta \tilde{\calL}]_{l}
    & \leq \alpha_{\delta} r_{c}^l + r_h^l
    - \frac{\rho}{2} \| \mathrm{\vC}_{y} [\Delta \by]_{l} \|_2^2
    \\
    &~~~~~~~
    - \frac{\rho}{2} \| \vC_g [ \Delta \bnu]_{l} \|_2^2
    - \hat{\alpha}_{\delta} \| [\Delta \blambda]_{l} \|_2^2
\end{aligned}
\label{convergence analysis descent relation in Ltilde}
\end{equation}
By summing \eqref{convergence analysis descent relation in Ltilde} over $l$, we get
\begin{align}
    & \sum_{l= 0}^{\infty}
    \bigg( \frac{\rho}{2} \| \mathrm{\vC}_{y} [\Delta \by]_{l} \|_2^2
    + \frac{\rho}{2} \| \vC_g [ \Delta \bnu]_{l} \|_2^2
    + \hat{\alpha}_{\delta} \| [\Delta \blambda]_{l} \|_2^2 \bigg) \nonumber
    \\ &~~~~~
    \leq
    \tilde{\calL}^0 - \tilde{\calL}^{\infty}
    + \sum_{l=0}^{\infty} ( \alpha_{\delta} r_{c}^l + r_h^l )
\end{align}
Assuming there exists $\bar{l} < \infty$ such that $\alpha_{\delta} r_{c}^l + r_h^l \leq 0$ for all $l \geq \bar{l}$. This assumption is based on the empirical data observed in all the experiments disclosed in the Section \ref{sec: Simulation}. Based on this assumption, we have $ \sum_{l=0}^{\infty} ( \alpha_{\delta} r_{c}^l + r_h^l ) < \infty$. Further, we have $\tilde{\calL}^0 < \infty$ and $\tilde{\calL}^{\infty}$ is lower-bounded which  gives $- \tilde{\calL}^{\infty} < \infty$. Therefore, we get
\begin{equation}
\begin{aligned}
    & \sum_{l= 0}^{\infty}
    \bigg( \frac{\rho}{2} \| \mathrm{\vC}_{y} [\Delta \by]_{l} \|_2^2
    + \frac{\rho}{2} \| \vC_g [ \Delta \bnu]_{l} \|_2^2
    + \hat{\alpha}_{\delta} \| [\Delta \blambda]_{l} \|_2^2 \bigg)
    < \infty \nonumber
\end{aligned}
\end{equation}
 On the LHS of the above, we have a summation of infinite number of non-negative terms. This summation can be upper bounded by a finite value only when $\| [\Delta \blambda]_{l} \|_2^2 \rightarrow 0$, $\| \mathrm{\vC}_{y}  [\Delta \by]_{l} \|^2_{2}  \rightarrow 0$, $\| \vC_g [ \Delta \bnu]_{l} \|_2^2 \rightarrow 0$ as $l \rightarrow \infty$. 

 Let $\by^*, \bnu^*, \blambda^*$ be the limit points (i.e., $\lim_{l \to \infty} \by^{l} = \by^{*}$, $\lim_{l \to \infty} \bnu^{l} = \bnu^{*}$, and $\lim_{l \to \infty} \blambda^{l} = \blambda^{*}$). We have that the limit points satisfy the KKT conditions (\eqref{convergence: KKT gradient condition y - main} - \eqref{convergence: KKT h constraint slackness condition - main}) of the original Problem \ref{Problem: DRO - compact form} except the consensus constraint \eqref{convergence: KKT consensus constraint satisfaction - main}. From the dual update step \eqref{convergence analysis: dual update}, we have $\mathrm{\vC}_{y} \by^{*} - \mathrm{\vC}_g \bnu^{*} = \delta \blambda^{*}/\rho$. By using the similar intuition as in Theorem 1 of \cite{yang2022proximal}, we have that points $\by^*, \bnu^*, (1+ \delta) \lambda^{*}$ are $\delta \rho^{-1} \|\blambda^{*}\|_2$- approximate stationary points.

\section*{References}
\bibliographystyle{IEEEtran}
\bibliography{IEEEabrv, references}

\end{document}